\documentclass[12pt,twoside,reqno,openany]{amsart}

\usepackage{amsbsy,amscd,amsfonts,amsmath,amssymb,amsthm,color,
fancybox,fancyhdr,footmisc,graphics,graphicx,ifthen,mathrsfs,
multicol,pdfpages,rotating,times,wasysym}
\usepackage[dvipsnames,svgnames,x11names]{xcolor}

\usepackage{pifont}

\usepackage[all]{xy}
\usepackage[utf8]{inputenc}
\usepackage[T1]{fontenc}
\usepackage{mathtools}
\newtagform{EngelLie}[\scriptstyle]{$}{$}
\makeatletter\newcommand{\leqnomode}{\tagsleft@true}
\newcommand{\reqnomode}{\tagsleft@false}\makeatother

\newtheorem{Theorem}[equation]{Theorem}

\newtheorem{Proposition}[equation]{Proposition}

\newtheorem{Lemma}[equation]{Lemma}

\newtheorem{Assertion}[equation]{Assertion}

\newtheorem{Observation}[equation]{Observation}

\theoremstyle{definition}

\newtheorem{Terminology}[equation]{Terminology}

\newtheorem{Problem}[equation]{Problem}

\newcommand{\C}{\mathbb{C}}

\newcommand{\N}{\mathbb{N}}

\renewcommand{\P}{\mathbb{P}}

\newcommand{\R}{\mathbb{R}}

\newcommand{\W}{\mathbb{W}}
\newcommand{\X}{\mathbb{X}}
\newcommand{\Y}{\mathbb{Y}}
\newcommand{\Z}{\mathbb{Z}}

\newcommand{\GG}{\text{\sc g}}

\newcommand{\KK}{\text{\sc k}}

\newcommand{\NN}{\text{\sc n}}

\newcommand{\RR}{\text{\sc r}}

\newcommand{\TT}{\text{\sc t}}

\newcommand{\daux}{{\text{\usefont{T1}{qcs}{m}{sl}d}}}

\newcommand{\CA}{{\text{\usefont{T1}{qcs}{m}{sl}CA}}}
\newcommand{\CR}{{\text{\usefont{T1}{qcs}{m}{sl}CR}}}
\newcommand{\CMR}{{\text{\usefont{T1}{qcs}{m}{sl}CMR}}}

\definecolor{blue}{cmyk}{1.,1.,0.,0.63}
\definecolor{red}{cmyk}{0.,1.,1.,0.63}
\definecolor{green}{cmyk}{1.,0.,1.,0.63}
\definecolor{black}{cmyk}{1.,1.,1.,1.}

\newcommand{\red}{\textcolor{red}}

\makeatletter
\renewcommand{\@fnsymbol}[1]
{\ensuremath{\ifcase#1\or $*$\or $**$\or $***$\or $****$\or $*****$
\else\@ctrerr\fi}}
\makeatother

\newcommand{\HEAD}[2]{%
\pagestyle{fancy}
\fancyhead[RO]{\tiny\sf\thepage}
\fancyhead[CO]{{\tiny\sf #1}}
\fancyhead[LE]{\tiny\sf\thepage}
\fancyhead[CE]{{\tiny\sf #2}}
\fancyfoot{}}

\numberwithin{equation}{section}

\newcommand{\Section}[1]{
\renewcommand{\thesection}{\bf\arabic{section}}
\section{#1}
\renewcommand{\thesection}{\arabic{section}}}

\newcommand{\style}[1]{\text{\footnotesize{\sf #1}}}

\newcommand{\stylesmall}[1]{{\sf #1}}

\newcommand{\constant}{\style{constant}}

\renewcommand{\cos}{\style{cos}}

\renewcommand{\exp}{\style{exp}}

\renewcommand{\lim}{\style{lim}}

\renewcommand{\log}{\style{log}}
\newcommand{\logsmall}{\stylesmall{log}}

\renewcommand{\max}{\style{max}}

\renewcommand{\min}{\style{min}}

\renewcommand{\sin}{\style{sin}}

\newcommand{\domaink}{\angle\!\square}

\newcommand{\explicationmath}[1]{\text{\scriptsize\sf [$#1$]}}

\newcommand{\explicationtext}[1]{\text{\scriptsize\sf [#1]}}

\newcommand{\linestop}{\medskip\centerline{\bf 
-----------------}\medskip}

\newcommand{\vf}{\vfill\end{document}}

\newcommand{\zero}[1]{\underline{#1}_{\red{\circ}}}

\begin{document}

\setcounter{section}{0}

$\:$


\begin{center}

{\large\bf 
Degrees
$d \geqslant \big( \sqrt{n}\, \log\, n\big)^n$ 
and 
$d \geqslant \big( n\, \log\, n\big)^n$}

\medskip

{\large\bf 
in the Conjectures of Green-Griffiths and of Kobayashi}
\label{constant-power-n}

\bigskip\bigskip

Jo\"el~{\sc Merker}\footnotemark[1]
and 
The-Anh {\sc Ta}\footnotemark[1]

\footnotetext[1]{Laboratoire de Mathématiques d'Orsay,
Université Paris-Sud, CNRS, Université Paris-Saclay, 
91405 Orsay Cedex, France.
{\bf joel.merker@math.u-psud.fr} and 
{\bf the-anh.ta@math.u-psud.fr}} 

\smallskip

{\large\footnotesize\sf D\'epartement de Math\'ematiques d'Orsay}

{\large\footnotesize\sf Université Paris-Sud, France}

\end{center}\bigskip

\begin{center}
\begin{minipage}[t]{12.5cm}
\parindent 0.53cm
\scriptsize
\noindent
{\sc Abstract}.
Once first answers in any dimension to the Green-Griffiths
and Kobayashi conjectures for generic algebraic hypersurfaces
$\mathbb{X}^{n-1} \subset \mathbb{P}^n(\mathbb{C})$ 
have been reached, the principal
goal is to decrease (to improve) the degree bounds, knowing
that the `celestial' horizon lies near $d \geqslant 2n$.

For Green-Griffiths algebraic degeneracy of entire holomorphic
curves, we obtain:
\[
d
\,\geqslant\,
\big(\sqrt{n}\,{\sf log}\,n\big)^n,
\]
and for Kobayashi-hyperbolicity (constancy of entire curves),
we obtain:
\[
d
\,\geqslant\,
\big(n\,{\sf log}\,n\big)^n.
\]
The latter improves $d \geqslant n^{2n}$ obtained by Merker in
{\tiny\sf arxiv.org/1807/11309/}. 

Admitting a certain technical conjecture
$I_0 \geqslant \widetilde{I}_0$, the method employed
(Diverio-Merker-Rousseau, B\'erczi, Darondeau)
conducts to constant power $n$, namely to:
\[
d
\,\geqslant\,
2^{5n}
\ \ \ \ \ \ \ \ \ \ \ \ \ \ \ \ \ \
\text{and, respectively, to:}
\ \ \ \ \ \ \ \ \ \ \ \ \ \ \ \ \ \
d
\,\geqslant\,
4^{5n}.
\]

In Spring 2019, a forthcoming prepublication based on intensive
computer explorations will present several subconjectures supporting
the belief that $I_0 \geqslant \widetilde{I}_0$, a conjecture which
will be established up to dimension $n = 50$.
\end{minipage}
\end{center}

\Section{\bf Introduction}
\label{introduction-constant-power-n}
\HEAD{{\ref{introduction-constant-power-n}}.~{\sf Introduction}
}{
Jo\"el {\sc Merker} and The-Anh Ta,
D\'epartement de Math\'ematiques d'Orsay, 
Universit\'e Paris-Sud, France}

The goal is to establish that generic algebraic hypersurfaces
of the projective space satisfy the Green-Griffiths
conjecture, as well as their complements,
with improvements on lower degree bounds. 

\begin{Theorem}
\label{Theorem-GG-sqrt-n-log-log-n}
For a generic hypersurface $\X^{n-1} \subset \P^n(\C)$ of degree:
\[
d
\,\geqslant\,
\big(
\sqrt{n}\,
\log\,n
\big)^n
\eqno
{\scriptstyle{(\forall\,n\,\geqslant\,\NN_{\GG\GG})}},
\]

\smallskip\noindent{\bf (1)}\,
there exists a proper subvariety $\Y \subset \P^n$ of codimension
$\geqslant 2$ such that all nonconstant entire holomorphic curves $f
\colon \C \longrightarrow \P^n \backslash \X$ are in fact contained in
$\Y \supset f(\C)$;

\smallskip\noindent{\bf (2)}\,
there exists a proper subvariety $\W \subset \X$ of codimension
$\geqslant 2$ such that all nonconstant entire holomorphic curves $f
\colon \C \longrightarrow \X$ are in fact contained in 
$\W \supset f(\C)$.

\end{Theorem}

This lower degree bound: 
\[
d
\,\geqslant\,
\daux_{\GG\GG}(n)
\,:=\, 
\big( \sqrt{n}\, \log\,n \big)^n
\]
improves $d \geqslant
(5n)^2\, n^n$ of~{\cite{Darondeau-IMRN-2016}} 
and improves $d \geqslant 2^{\,n^5}$
of~{\cite{Diverio-Merker-Rousseau-2010}}.
In the demonstrations, we will treat mainly
the details of the complement case {\small\bf (1)},
since the computations in the
compact case {\small\bf (2)} are essentially similar,
thanks to Darondeau's
works~{\cite{Darondeau-Fourier-2016,
Darondeau-MZ-2016, Darondeau-IMRN-2016}}.

By~{\cite{Riedl-Yang-2018}}, any solution to the Green-Griffiths
conjecture in all dimensions $n$ for hypersurfaces of degrees $d
\geqslant \daux_{\GG\GG} (n)$ implies a solution to the Kobayashi
conjecture in all dimensions $n$ for hypersurfaces of degrees:
\[
d
\,\geqslant\,
\daux_{\KK}(n) 
\,:=\,
\daux_{\GG\GG}(2n). 
\]
Rounding off a small technical improvement
of Theorem~{\ref{Theorem-GG-sqrt-n-log-log-n}} in order 
to present only an elegant 
degree bound, we obtain as a corollary the following

\begin{Theorem}
\label{Theorem-Kb-n-log-n}
For a generic hypersurface $\X^{n-1} \subset \P^n(\C)$ of degree
\[
d
\,\geqslant\,
\big(n\,\log\,n\big)^n
\eqno
{\scriptstyle{(\forall\,n\,\geqslant\,\NN_\KK)}}
\]

\smallskip\noindent{\bf (1)}\,
$\P^n \big\backslash \X^{n-1}$ is Kobayashi-hyperbolically
imbedded in $\P^n${\em ;}

\smallskip\noindent{\bf (2)}\,
$\X^{n-1}$ is Kobayashi-hyperbolic.

\end{Theorem}

An inspection of the end of Section~{\ref{final-minorations}} shows
that the dimensions $\NN_{\GG\GG}$ and $\NN_\KK$ at which these
statements begin to hold true can be made effective.

Theorem~{\ref{Theorem-Kb-n-log-n}}
improves the degree bound $d \geqslant n^{2n}$ 
obtained in~{\cite{Merker-2018}}.
For standard presentations of the research field,
and for up-to-date history, including 
degree bound comparisons, 
the reader is referred to the introductions 
of the articles~{\cite{Merker-2018, 
Riedl-Yang-2018,
Brotbek-Deng-2018,
Demailly-2018,
Huynh-2016,
Deng-2016,
Brotbek-2017,
Darondeau-IMRN-2016,
Darondeau-MZ-2016,
Darondeau-Fourier-2016,
Siu-2015,
Berczi-2018,
Diverio-Merker-Rousseau-2010}},
listed in chronological order of prepublication.

Under the technical assumption (or conjecture):
\[
I_0
\,\geqslant\,
\widetilde{I}_0,
\]
the explanation of which the reader will find in
Section~{\ref{preliminary-link-Darondeau}}, 
and which is equivalent to
Problem~{\ref{Problem-A-C}}, 
we obtain better results.

\begin{Theorem}
\label{Theorem-GG-constant-n}
If $I_0 \geqslant \widetilde{I}_0$ holds true, then 
for a generic hypersurface $\X^{n-1} \subset \P^n(\C)$ of degree:
\[
d
\,\geqslant\,
2^{\,5n}
\eqno
{\scriptstyle{(\forall\,n\,\geqslant\,10)}},
\]
the two conclusions {\small\bf (1)} and {\small\bf (2)}
of Theorem~{\ref{Theorem-GG-sqrt-n-log-log-n}} hold true.
\end{Theorem}

Similarly, we also obtain as corollary the

\begin{Theorem}
\label{Theorem-Kb-constant-n}
Under the same technical assumption $I_0 \geqslant
\widetilde{I}_0$, 
the conclusions {\small\bf (1)} and {\small\bf (2)}
of Theorem~{\ref{Theorem-Kb-n-log-n}} hold true
in degree: 
\[
d
\,\geqslant\,
4^{\,5n}
\eqno
{\scriptstyle{(\forall\,n\,\geqslant\,20)}}.
\]
\end{Theorem}

\smallskip\noindent{\bf Acknowledgments.}
In 2013, 2014, the first author exchanged
with Lionel Darondeau. 

\Section{\bf Preliminary: Link with Darondeau's Work}
\label{preliminary-link-Darondeau}
\HEAD{{\ref{preliminary-link-Darondeau}}.~{\sf Preliminary: 
Link with Darondeau's Work }
}{
Jo\"el {\sc Merker} and The-Anh Ta,
D\'epartement de Math\'ematiques d'Orsay, 
Universit\'e Paris-Sud, France}

This section continues~{\cite{Darondeau-IMRN-2016}},
and goes slightly beyond. The jet order $\kappa = n$
will be chosen equal to the dimension $n$, because
some reflections on the concerned estimates 
convince that any choice of $\kappa > n$ cannot improve
the degree bound anyway.

Let $n \geqslant 1$ be an integer. Let $t_1, \dots, t_n$
be formal variables. Introduce:
\[
C(t_1,\dots,t_n)
\,:=\,
\prod_{1\leqslant i<j\leqslant n}\,
\frac{t_j-t_i}{t_j-2\,t_i}\,
\prod_{2\leqslant i<j\leqslant n}\,
\frac{t_j-2\,t_i}{t_j-2\,t_i+t_{i+1}}.
\]
As explained in~{\cite{Darondeau-IMRN-2016}},
this rational expression possesses 
an iterated Laurent series at the origin as:
\[
C(t)
\,=\,
\sum_{
k_1,\dots,k_n\,\in\,\Z
\atop
k_1+\cdots+k_n\,=\,0}\,
C_{k_1,\dots,k_n}\,
t_1^{k_1}
\cdots
t_n^{k_n},
\]
for certain coefficients $C_{k_1, \dots, k_n}$; soon, this object
$C(t_1, t_2, \dots, t_n)$ will
be re-interpreted as a standard converging
power series $C(w_2, \dots, w_n)$  in
terms of alternative
new variables $(w_2, \dots, w_n)$, hence it is not necessary to recall
what an iterated Laurent series is.

For certain integer weights $a_1, \dots, a_n \in \N^\ast$, 
introduce also an expression which comes from an application
of the so-called {\sl holomorphic Morse inequalities:} 
\[
f_0(t)
\,:=\,
\big(
a_1t_1
+\cdots+
a_nt_n
\big)^{n^2}.
\]
It expands:
\[
f_0(t)
\,=\,
\sum_{
m_1,\dots,m_n\geqslant 0
\atop
m_1+\cdots+m_n=n^2}\,
\frac{(n^2)!}{m_1!\,\cdots\,m_n!}\,
\big(a_1\,t_1\big)^{m_1}
\cdots
\big(a_n\,t_n\big)^{m_n},
\]
by means of (integer) multinomial coefficients:
\[
M_{m_1,\dots,m_n}
\,:=\,
\frac{(n^2)!}{m_1!\,\cdots\,m_n!}.
\]

It is well known that the binomial $\binom{2n}{n}$ is the
unique largest one among all the $\binom{2n}{i}$ with $0 \leqslant
i \leqslant 2n$. In fact, an application of Stirling's 
asymptotic formula:
\[
n!
\underset{n\,\to\,\infty}{\,\,\,\sim\,\,\,}
\sqrt{2\pi\,n}\,
\Big(
\frac{n}{e}
\Big)^n\,
\bigg[
1
+
\frac{1}{12\,n}
+
\frac{1}{288\,n^2}
-
\frac{139}{51\,840\,n^3}
-
\frac{571}{2\,488\,320\,n^4}
+
{\rm O}
\big(
\frac{1}{n^5}
\big)
\bigg],
\]
shows that 
asymptotically as 
$n \longrightarrow \infty$:
\[
\binom{2\,n}{n}
\,\,\sim\,\,
\frac{2^{2\,n}}{\sqrt{\pi\,n}}\,
\bigg[
1
-
\frac{1}{8\,n}
+
\frac{1}{128\,n^2}
+
\frac{5}{1024\,n^3}
-
\frac{21}{32\,768\,n^4}
+
{\rm O}
\Big(
\frac{1}{n^5}
\Big)
\bigg].
\]
Similarly, the {\sl central multinomial
coefficient:}
\[
M_{n,\dots,n}
\,:=\,
\frac{(n^2)!}{n!\,\cdots\,n!}
\,=\,
\frac{(n^2)!}{(n!)^n},
\]
happens to be the unique largest one, as states the next observation
({\em see} also Lemma~{\ref{Lemma-M-less-1}}).

\begin{Lemma}
For all integers $m_1, \dots, m_n \geqslant 0$ with $m_1 + \cdots +
m_n = n^2$ and $(m_1, \dots, m_n) \neq (n, \dots, n)$, the
corresponding multinomial coefficients are smaller than the central
one:
\[
M_{m_1,\dots,m_n}
\,<\,
M_{n,\dots,n}.
\]
\end{Lemma}

\proof
This amounts to verify that:
\[
\frac{n!}{m_1!}
\,\cdots\,
\frac{n!}{m_i!}
\,\cdots\,
\frac{n!}{m_n!}
\overset{\text{\bf ?}}{\,\,<\,\,}
1.
\]
The $m_i = n$ are neutral, for $\frac{n!}{n!} = 1$.  By assumption, at
least one $m_i \neq n$.

\smallskip\noindent$\bullet$\,
When $m_i < n$, simplify:
\[
\frac{n!}{m_i!}
\,=\,
n\,(n-1)\cdots\big(m_i+1).
\]

\smallskip\noindent$\bullet$\,
When $m_i > n$, simplify:
\[
\frac{n!}{m_i!}
\,=\,
\frac{1}{(n+2)(n+1)\cdots m_i}.
\]
After these simplifications:
\[
\prod_{1\leqslant i\leqslant n}\,
\frac{n!}{m_i!}
\,=\,
\frac{
\prod_{m_i<n}\,\,\,\,\,\,\,
(n-0)\,(n-1)\cdots(m_i+1)
}{
\prod_{m_i>n}\,
(n+1)(n+2)\cdots m_i
}.
\]
Since $m_1 + \cdots + m_n = n^2$, 
the number of factors in the numerator is the same as that 
in the denominator, and since {\em each} factor upstairs
is $\leqslant n$, while {\em each} factor downstairs 
is $\geqslant n+1$, the result is indeed $< 1$. 
\endproof

A further application of Stirling's formula shows that, 
asymptotically as $n \longrightarrow \infty$:
\[
\frac{(n^2)!}{n!\,\cdots\,n!}
\,\,\sim\,\,
n^{\,n^2-\frac{n}{2}+1}\,
\frac{1}{(2\,\pi)^{\frac{n-1}{2}}}\,
\frac{1}{e^{\frac{1}{12}}}\,\,
\bigg[
1
+
\frac{31}{360\,n^2}
+
\frac{5287}{181\,440\,n^4}
+
{\rm O}
\Big(
\frac{1}{n^6}
\Big)
\bigg].
\]

\begin{Terminology}
Call the coefficient of $t_1^n \cdots t_n^n$ in $f_0(t)$:
\[
\aligned
\widetilde{I}_0
\,:=\,
&\,
\big[t_1^n\cdots t_n^n\big]
\big(f_0(t)\big)
\\
\,=\,
&\,
\frac{(n^2)!}{n!\,\cdots\,n!}\,
a_1^n\cdots a_n^n
\endaligned
\]
the {\sl central monomial}.
\end{Terminology}

Since $a_1, \dots, a_n \in \N^\ast$, this is a large integer.
The notation $\widetilde{I}_0$ is borrowed 
from~{\cite{Darondeau-IMRN-2016}}.

In fact, Appendices~1 and~2 of~{\cite{Darondeau-IMRN-2016}}
provided almost all the details to verify that the choice
of weights:
\[
a_i
\,:=\,
r^{n-i}
\eqno
{\scriptstyle{(1\,\leqslant\,i\,\leqslant\,n)}},
\]
for some constant $r$ independent of $n$, shall offer a degree bound
in the Green-Griffiths conjecture of the form:
\[
d
\,\geqslant\,
\constant^{\,n}.
\]
which would improve the current $d \gtrsim n^{n}$
obtained in~{\cite{Darondeau-IMRN-2016, Merker-2018}}.

For a certain nefness condition 
required to apply the holomorphic Morse inequalities, 
it is necessary to have at least:
\[
r
\,\geqslant\,
3.
\]
It is also allowed to take $r$ larger, for instance:
\[
r
\,=\,
9
\ \ \ \ \ \ \ \ \ \ \ \ \
\text{or}
\ \ \ \ \ \ \ \ \ \ \ \ \
r
\,=\,
12
\ \ \ \ \ \ \ \ \ \ \ \ \
\text{or}
\ \ \ \ \ \ \ \ \ \ \ \ \
r
\,=\,
20,
\]
but one should try {\em not} to choose 
$r$ increasing with $n$, like for instance $r = \sqrt{n}$,
since the final degree bound would otherwise be
(explanations will appear later): 
\[
d
\,\gtrsim\,
(\sqrt{n})^n
\,\gg\,
\constant^{\,n}.
\]
In~{\cite{Darondeau-IMRN-2016}}, the choice was $r := n$, and
this conducted to $d \gtrsim n^n$. 

With a fixed (bounded) constant $r \geqslant 3$, the
final degree bound for Green-Griffiths will be close to:
\[
d
\,\gtrsim\,
\big(r\,(1+\varepsilon(r))\big)^n
\,=\,
\constant^{\,n},
\]
as we will verify in details later.  The only remaining substantial
piece of work to be done is to solve the following

\begin{Problem}
\label{Problem-I-0}
{\sl 
With the choice of weights: 
\[
a_1
\,:=\,
r^{n-1},\ \ \
a_2
\,:=\,
r^{n-2},\ \ \
\dots\dots,\ \ \
a_{n-1}
\,:=\,
r,\ \ \
a_n
\,:=\,
1,
\]
to show that the coefficient of the monomial $t_1^n \cdots t_n^n$ in
the product $C(t) \cdot f_0(t)$, namely:}
\[
I_0
\,:=\,
\big[t_1^n\cdots t_n^n\big]
\Big(
C(t_1,\dots,t_n)
\cdot
f_0(t_1,\dots,t_n)
\Big)
\]
{\sl is at least equal to the central monomial:}
\[
\aligned
I_0
&
\overset{\text{\bf ?}}{\,\geqslant\,}
\widetilde{I}_0
\\
&
\,=\,
\frac{(n^2)!}{(n!)^n}\,
r^{\,n\frac{n(n-1)}{2}}.
\endaligned
\]
\end{Problem}

In fact, several computer experiments convince that
instead of $\frac{I_0}{\widetilde{I}_0} \geqslant 1$, 
a better inequality seems to hold: 
\[
\frac{I_0}{\widetilde{I}_0}
\,\gtrsim\, 
\big(\constant_r\big)^{\!n}, 
\]
for some $\constant_r > 1$ which depends on $r$, and is
closer and closer to $1$ when $r$ increases. So experimentally,
$I_0 \geqslant \widetilde{I}_0$ is more than true.
The goal is to set up a proof. 

\smallskip

We start in Section~{\ref{end-proof-under-hypothesis}}
by verifying that a proof of $\frac{I_0}{\widetilde{I}_0}
\geqslant 1$ implies a degree bound for Green-Griffiths 
of the announced form $d \geqslant \constant^{\,n}$; 
this task was already almost completely performed by Darondeau 
in~{\cite{Darondeau-IMRN-2016}}.

\smallskip

Then in subsequent sections, we study the product
$C(t_1, \dots, t_n)$ and we establish 
$I_0 \geqslant \widetilde{I}_0$.

\Section{\bf End of Proof of 
Theorem~{\ref{Theorem-GG-constant-n}}}
\label{end-proof-under-hypothesis}
\HEAD{{\ref{end-proof-under-hypothesis}}.~{\sf End of Proof 
of Theorem~{\ref{Theorem-GG-constant-n}}}
}{
Jo\"el {\sc Merker} and The-Anh Ta,
D\'epartement de Math\'ematiques d'Orsay, 
Universit\'e Paris-Sud, France}

It essentially suffices to read Appendices~1 and~2 
of~{\cite{Darondeau-IMRN-2016}}, with in mind that 
Darondeau's (simplifying) choice:
\[
a_i
\,:=\,
n^{n-i}
\eqno
{\scriptstyle{(1\,\leqslant\,i\,\leqslant\,n)}},
\]
should be replaced with the choice:
\[
a_i
\,:=\,
r^{n-i}
\eqno
{\scriptstyle{(1\,\leqslant\,i\,\leqslant\,n)}},
\]
where $r \geqslant 3$ is a fixed constant. Later, we will see that the
choice $r = 3$ might expose to some computational difficulties, while
as soon that $r \geqslant 9$, a serendipitous positivity property
occurs. In any case, the estimates of the mentioned Appendix~2
were prepared in advance to work for any
choice of $r = 3, 9, 12, 20, \log\, n, \sqrt{n}, n$, while they were
applied in~{\cite{Darondeau-IMRN-2016}} 
to $r = n$ by lack of a solution to
Problem~{\ref{Problem-I-0}}. Before solving this problem in the
next sections, let us admit
temporarily that it has a positive answer for a certain fixed:
\[
9
\,\leqslant\,
r
\,\leqslant\,
20
\eqno
{\scriptstyle{(\text{\rm hypothesis throughout})}}.
\]

\proof[End of proof of 
Theorem~{\ref{Theorem-GG-constant-n}}]
In the notations of~{\cite{Darondeau-IMRN-2016}},
the lower degree bound:
\[
d
\,\geqslant\,
\daux_{\rm GG}(n)
\]
is determined by the largest root of a certain polynomial equation:
\[
d^n\,I_0
+
d^{n-1}\,I_1
+\cdots+
d^{n-p}\,I_p
+\cdots+
I_n
\,\,=\,\,
0,
\]
with $I_0 > 0$. 
Of course, $I_0$ is the same 
as in Problem~{\ref{Problem-I-0}}, hence we assume temporarily 
not only that it is positive, but also that it is quite large:
\[
\frac{I_0}{\widetilde{I}_0}
\,\geqslant\,
1.
\]
We refer to~{\cite{Darondeau-IMRN-2016}} for a presentation
of the other coefficients $I_p$.

\begin{Proposition}
\label{Proposition-r-3-n}
The polynomial in the degree $d$ of a hypersurface
$\X^n \subset \P^{n-1}(\C)${\em :}
\[
d^n\,I_0
+
d^{n-1}\,I_1
+\cdots+
d^{n-p}\,I_p
+\cdots+
d\,I_{n-1}
+
I_n
\]
takes positive values for all degrees:
\[
\aligned
d
&
\,\geqslant\,
25\,n^2
\cdot
\big(
r+3
\big)^n.
\\
&
\,=:\,
\daux_{\GG\GG}(n,r)
\endaligned
\]
\end{Proposition}

In fact, a glance at the end of the proof
shows a slightly better, though more complicated:
\[
\daux_{\GG\GG}(n,r)
\,:=\,
\big(20\,n^2+4\,n\big)
\cdot
\frac{r^3}{(r-1)^3\,(r+3)}
\cdot
\big(
r+3
\big)^n.
\]

Theorem~{\ref{Theorem-GG-constant-n}} 
terminates by checking on a computer that:
\[
2^{\,5n}
\,\geqslant\,
\big(20\,n^2+4\,n\big)
\cdot
\frac{r^3}{(r-1)^3\,(r+3)}
\cdot
\big(
r+3
\big)^n
\eqno
{\scriptstyle{(\forall\,n\,\geqslant\,20)}},
\]
for any choice of $9 \leqslant r \leqslant 20$.
\endproof

\proof[Proof of Proposition~{\ref{Proposition-r-3-n}}]
In~{\cite{Darondeau-IMRN-2016}}, the 
p\^ole order of so-called
{\sl slanted vector fields} $c_n := n(n+2)$ is used.
But the article~{\cite{Darondeau-MZ-2016}} improves it to:
\[
c_n
\,:=\,
5\,n-2.
\]
Then with $c := c_n + 1$, the quantity $\frac{c+2}{2}$
appears several times in~{\cite{Darondeau-IMRN-2016}}, so we may read:
\[
\frac{c+2}{2}
\,=\,
\frac{5\,n+1}{2}.
\]

Next, with:
\[
a_i
\,:=\,
r^{n-i}
\eqno
{\scriptstyle{(1\,\leqslant\,i\,\leqslant\,n)}},
\]
set:
\[
\mu(a)
\,:=\,
1\,a_1
+
2\,a_2
+\cdots+
n\,a_n,
\]
and for all $1 \leqslant p \leqslant n$, set:
\[
\widetilde{I}_p
\,:=\,
\underbrace{
\frac{(n^2)!}{(n!)^n}\,
a_1^n\cdots a_n^n}_{
\text{\rm recognize}\,\,\widetilde{I}_0}\,\,
\big(
2\,n\,\mu(a)
\big)^p\!
\sum_{1\leqslant i_1<\cdots<i_p\leqslant n}\,
\frac{1}{a_{i_1}}\,
\cdots\,
\frac{1}{a_{i_p}},
\]
Importantly, 
Lemma~A.6 on page~1919 of Appendix~2 shows that:
\[
\frac{\vert I_p\vert}{\widetilde{I}_p}
\,\,\leqslant\,\,
\frac{5n+1}{2}
\cdot
\vert B\vert
\Big(
\frac{2n\mu(a)h}{a_1},\dots,\frac{2n\mu(a)h}{a_n}
\Big)
\cdot
\vert C\vert
\Big(
\frac{1}{a_1},\dots,\frac{1}{a_n}
\Big)
\eqno
{\scriptstyle{(1\,\leqslant\,p\,\leqslant\,n)}}.
\]
It is not necessary to dwell into details about the middle quantity 
$\vert B\vert$, since Lemma~A.7 on page~1920 shows that for any
choice of weights $a_1, \dots, a_n$:
\reqnomode\usetagform{EngelLie}
\begin{align}
\vert B\vert
\Big(
\frac{2n\mu(a)h}{a_1},\dots,\frac{2n\mu(a)h}{a_n}
\Big)
&
\,\,\leqslant\,\,
\Big(
\frac{2\,n}{2\,n-1}
\Big)^{\!n+1}
\notag
\\
&
\,\,\leqslant\,\,
2
\tag{(\forall\,n\,\geqslant\,4\text{\rm \,\,--\,\,exercise}).}
\end{align}
Consequently, we get:
\[
\frac{\vert I_p\vert}{\widetilde{I}_p}
\,\,\leqslant\,\,
\big(5\,n+1\big)
\cdot
\vert C\vert
\Big(
\frac{1}{a_1},\dots,\frac{1}{a_n}
\Big)
\eqno
{\scriptstyle{(1\,\leqslant\,p\,\leqslant\,n)}}.
\]

Next, page~1914 uses the control:
\[
\vert C\vert
\Big(
\frac{1}{a_1},\dots,\frac{1}{a_n}
\Big)
\,\,\leqslant\,\,
\widehat{C}
\Big(
\frac{1}{a_1},\dots,\frac{1}{a_n}
\Big),
\]
by the `majorant' series:
\[
\widehat{C}\big(t_1,\dots,t_n\big)
\,:=\,
\prod_{1\leqslant i<j\leqslant n}\,
\frac{t_j-t_i}{t_j-2\,t_i}\,\,
\prod_{2\leqslant i<j\leqslant n}\,
\frac{t_j-2\,t_i}{t_j-2\,t_i-t_{i-1}}.
\]
Replacing the formal variables by the
inverses of the weights, we get:
\[
\widehat{C}
\Big(
\frac{1}{a_1},\dots,\frac{1}{a_{n-1}},\frac{1}{a_n}
\Big)
\,=\,
\prod_{1\leqslant i<j\leqslant n}\,
\frac{a_i/a_j-1}{a_i/a_j-2}\,\,
\prod_{2\leqslant i<j\leqslant n}\,
\frac{a_i/a_j-2}{a_i/a_j-2-a_i/a_{i-1}}.
\]
Since $a_i = r^{n-i}$ for $1 \leqslant i \leqslant n$, this
rewrites as:
\[
\aligned
\widehat{C}
\Big(
\frac{1}{r^{n-1}},\dots,\frac{1}{r},\frac{1}{1}
\Big)
&
\,\,=\,\,
\prod_{1\leqslant i<j\leqslant n}\,
\frac{r^{j-i}-1}{r^{j-i}-2}\,\,
\prod_{2\leqslant i<j\leqslant n}\,
\frac{r^{j-i}-2}{r^{j-i}-2-\frac{1}{r}}
\\
\explicationtext{Extract $i=1$}
\ \ \ \ \ \ \ \ \ \ \ \ \ \ \ \ \ \ \ \ \ \ \ \ \ \
&
\,\,=\,\,
\prod_{2\leqslant j\leqslant n}\,
\frac{r^{j-1}-1}{r^{j-1}-2}\,\,
\prod_{2\leqslant i<j\leqslant n}\,
\bigg[\,\,
\frac{r^{j-i}-1}{\zero{r^{j-i}-2}}\,\,
\frac{\zero{r^{j-i}-2}}{r^{j-i}-2-\frac{1}{r}}
\,\,\bigg]
\\
\explicationtext{Simplify}
\ \ \ \ \ \ \ \ \ \ \ \ \ \ \ \ \ \ \ \ \ \ \ \ \ \
&
\,\,=\,\,
\prod_{2\leqslant j\leqslant n}\,
\frac{r^{j-1}-1}{r^{j-1}-2}\,\,
\prod_{2\leqslant i<j\leqslant n}\,
\frac{r^{j-i+1}-r}{r^{j-i+1}-2\,r-1}
\,\,\bigg]
\\
\explicationtext{Rename indices}
\ \ \ \ \ \ \ \ \ \ \ \ \ \ \ \ \ \ \ \ \ \ \ \ \ \
&
\,\,=\,\,
\prod_{1\leqslant k\leqslant n-1}\,
\frac{r^k-1}{r^k-2}\,\,
\prod_{2\leqslant\ell\leqslant n-1}\,
\Big(
\frac{r^\ell-r}{r^\ell-2\,r-1}
\Big)^{n-\ell}.
\endaligned
\]

Using inequalities valid as soon as $r \geqslant 4$ hence for 
$r \geqslant 9$:
\[
\frac{1}{r^k-2}
\,\leqslant\,
\frac{2}{r^k}
\eqno
{\scriptstyle{(\forall\,k\,\geqslant\,1)}},
\]
the first product is bounded by a universal constant,
and even by a constant which decreases as $r$ increases:
\[
\aligned
\prod_{1\leqslant k\leqslant n-1}\,
\Big(
1
+
\frac{1}{r^k-2}
\Big)
&
\,\leqslant\,
\prod_{k=1}^\infty\,
\Big(
1
+
\frac{2}{r^k}
\Big)
\\
&
\,=\,
\exp
\bigg(
\sum_{k=1}^\infty\,
\log\,
\Big(
1
+
\frac{2}{r^k}
\Big)
\bigg)
\\
\explicationmath{\logsmall\,(1+\varepsilon)\leqslant 1+\varepsilon}
\ \ \ \ \ \ \ \ \ \ \ \ \ \ \ \ \ \ \ \ \ \ \ \ \ \
&
\,\leqslant\,
\exp\,
\bigg(
\frac{2}{r}\,
\sum_{k=0}^\infty\,
\frac{1}{r^k}
\bigg)
\\
&
\,=\,
\exp\,
\Big(
\frac{2}{r-1}
\Big)
\\
\explicationtext{Computer help}
\ \ \ \ \ \ \ \ \ \ \ \ \ \ \ \ \ \ \ \ \ \ \ \ \ \
&
\,\leqslant\,
1
+
\frac{3}{r}.
\endaligned
\]

The second product is bounded by a constant power $n-2$:
\[
\aligned
\prod_{2\leqslant\ell\leqslant n-1}\,
\Big(
\frac{r^\ell-r}{r^\ell-2\,r-1}
\Big)^{n-\ell}
&
\,\,=\,\,
\prod_{2\leqslant\ell\leqslant n-1}\,
\Big(
1
+
\frac{r+1}{r^\ell-2\,r-1}
\Big)^{n-\ell}
\notag
\\
&
\,\,\leqslant\,\,
\bigg(
\prod_{2\leqslant\ell\leqslant n-1}\,
\Big(
1
+
\frac{r+1}{r^\ell-2\,r-1}
\Big)
\bigg)^{n-2}
\notag
\\
&
\,\,\leqslant\,\,
\bigg(
\prod_{\ell=2}^\infty\,
\Big(
1
+
\frac{r+1}{r^\ell-2\,r-1}
\Big)
\bigg)^{n-2}.
\endaligned
\]
Let us estimate this constant, which depends on $r$:
\[
\alpha(r)
\,:=\,
\prod_{\ell=2}^\infty\,
\Big(
1
+
\frac{r+1}{r^\ell-2\,r-1}
\Big).
\]

\begin{Lemma}
For all $r \geqslant 6$, one has:
\[
1
+
\frac{r+1}{r^\ell-2\,r-1}
\,\,\leqslant\,\,
1
+
\frac{2}{r^{\ell-1}}
\eqno
{\scriptstyle{(\forall\,\ell\,\geqslant\,2)}}.
\]
\end{Lemma}

\proof
This is equivalent to:
\[
4\,r+2
\,\leqslant\,
r^\ell
-
r^{\ell-1}
\eqno
{\scriptstyle{(\forall\,r\,\geqslant\,6,\,\,
\forall\,\ell\,\geqslant\,2)}},
\]
which is easily checked, on a computer, to be true.
\endproof

Hence we can majorize still assuming $r \geqslant 9$
throughout:
\[
\aligned
\alpha(r)
&
\,\leqslant\,
\prod_{\ell=2}^\infty\,
\Big(
1
+
\frac{2}{r^{\ell-1}}
\Big)
\\
&
\,=\,
\exp\,
\bigg(
\sum_{\ell=0}^\infty\,
\log\,
\Big(
1
+
\frac{2}{r^{\ell+1}}
\Big)
\bigg)
\\
\explicationmath{\log\,(1+\varepsilon)\leqslant 1+\varepsilon}
\ \ \ \ \ \ \ \ \ \ \ \ \ \ \ \ \ \ \ \ \ \ \ \ \ \
&
\,\leqslant\,
\exp\,
\bigg(
\frac{2}{r}\,
\sum_{\ell=0}^\infty\,
\frac{1}{r^\ell}
\bigg)
\\
&
\,=\,
\exp\,
\Big(
\frac{2}{r-1}
\Big)
\\
\explicationtext{Computer help}
\ \ \ \ \ \ \ \ \ \ \ \ \ \ \ \ \ \ \ \ \ \ \ \ \ \
&
\,\leqslant\,
1
+
\frac{3}{r}.
\endaligned
\]

In summary, we have shown that:
\[
\frac{\vert I_p\vert}{\widetilde{I}_p}
\,\leqslant\,
\big(5\,n+1\big)
\cdot
\Big(
1
+
\frac{3}{r}
\Big)
\cdot
\Big(
1
+
\frac{3}{r}
\Big)^{n-2}
\eqno
{\scriptstyle{(\forall\,n\,\geqslant\,2)}}.
\]

Next, we estimate, still with $a_i = r^{n-i}$ for $i=1, \dots, n$:
\[
\aligned
\mu(a)
&
\,=\,
1\cdot a_1
+
2\cdot a_2
+\cdots+
(n-1)\,a_{n-1}
+
n\,a_n
\\
&
\,=\,
1\,r^{n-1}
+
2\,r^{n-2}
+\cdots+
(n-1)\,r^1
+
n\,r^0
\\
&
\,=\,
(n+1)\,
\big[
r^{n-1}
+
r^{n-2}
+\cdots+
r^1
+
r^0
\big]
\\
&
\ \ \ \ \
-\,
n\,r^{n-1}
-
(n-1)\,r^{n-2}
-\cdots-
2\,r^1
-
1\,r^0
\\
&
\,=\,
(n+1)\,
\frac{r^n-1}{r-1}
\\
&
\ \ \ \ \
-\,
\frac{n\,r^{n+1}-(n+1)\,r^n+1}{(r-1)^2},
\endaligned
\]
the result in this last line being obtained simply by differentiating
with respect to $r$ the classical:
\[
r^n
+
r^{n-1}
+\cdots+
r^2
+
r
+
1
\,\,=\,\,
\frac{r^{n+1}-1}{r-1}.
\]
A reduction to the same denominator contracts:
\[
\aligned
\mu(a)
&
\,=\,
\frac{r^{n+1}-(n+1)\,r+n}{(r-1)^2}
\\
&
\,\leqslant\,
\frac{r^{n+1}}{(r-1)^2}
\\
&
\,=\,
\frac{r}{(r-1)^2}\,
r^n.
\endaligned
\]

Next, consider generally a polynomial of degree 
$n \geqslant 1$ with complex coefficients $c_p \in \C$:
\[
c_0\,z^n
+
c_1\,z^{n-1}
+\cdots+
c_{n-1}\,z^1
+
c_n
\eqno
{\scriptstyle{(c_0\,\neq\,0)}}.
\]
Abbreviate:
\[
\KK_n
\,:=\,
\text{\rm unique positive zero of}\,\,
z^n-z^{n-1}-\cdots-z-1,
\]
which satisfies:
\[
1
\,<\,
\KK_n
\,<\,
2
\eqno
{\scriptstyle{(\text{\rm close to}\,2)}}.
\]

\begin{Theorem}
{\bf [Fujiwara]}
The moduli of all roots of $c_0 z^n + c_1 z^{n-1} + \cdots + 
c_n$ are bounded by:
\[
\max\,
\big\vert
{\sf roots}
\big\vert
\,\,\leqslant\,\,
\underbrace{\KK_n}_{<\,\,2}\,\,
\underset{1\leqslant p\leqslant n}{\max}\,
\sqrt[p]{\,
\frac{\vert c_p\vert}{\vert c_0\vert}\,}.
\eqno\qed
\]
\end{Theorem}

Now, come back to the polynomial $d^n\, I_0 + d^{n-1}\, I_1 + 
\cdots + I_n$ of Proposition~{\ref{Proposition-r-3-n}}. 
Thanks to Fujiwara:
\[
\aligned
\max\,
\big\vert
{\sf roots}
\big\vert
&
\,\,\leqslant\,\,
2\,\,
\underset{1\leqslant p\leqslant n}{\max}\,
\sqrt[p]{\,
\frac{\vert I_p\vert}{I_0}}
\\
\explicationmath{I_0\geqslant\widetilde{I}_0}
\ \ \ \ \ \ \ \ \ \ \ \ \ \ \ \ \ \ \ \ \ \ \ \ \ \
&
\,\,\leqslant\,\,
2\,\,
\underset{1\leqslant p\leqslant n}{\max}\,
\sqrt[p]{\,
\frac{\vert I_p\vert}{\widetilde{I}_0}}
\\
&
\,\,=\,\,
2\,\,
\underset{1\leqslant p\leqslant n}{\max}\,
\sqrt[p]{\,
\frac{\vert I_p\vert}{\widetilde{I}_p}
\cdot
\frac{\widetilde{I}_p}{\widetilde{I}_0}\,}
\\
\explicationtext{Seen above}
\ \ \ \ \ \ \ \ \ \ \ \ \ \ \ \ \ \ \ \ \ \ \ \ \ \
&
\,\,\leqslant\,\,
2\,\,
\underset{1\leqslant p\leqslant n}{\max}\,
\sqrt[p]{\,
\big(
5\,n+1
\big)\,
\Big(
1
+
\frac{3}{r}
\Big)^{n-1}
\cdot
\frac{\widetilde{I}_p}{\widetilde{I}_0}\,}
\\
&
\,\,\leqslant\,\,
2\,\,
\underset{1\leqslant p\leqslant n}{\max}\,
\sqrt[p]{\,
\big(
5\,n+1
\big)\,
\Big(
1+\frac{3}{r}
\Big)^{n-1}}
\cdot\,
\underset{1\leqslant p\leqslant n}{\max}\,
\sqrt[p]{\,
\frac{\widetilde{I}_p}{\widetilde{I}_0}}\,\,
\\
&
\,\,=\,\,
2\,\,
\big(
5\,n+1
\big)\,
\Big(
1
+
\frac{3}{r}
\Big)^{n-1}
\cdot
\underset{1\leqslant p\leqslant n}{\max}\,
\sqrt[p]{\,
\frac{\widetilde{I}_p}{\widetilde{I}_0}}\,\,
\cdot\,
\endaligned
\]

Next, coming back to the definition of $\widetilde{I}_p$, 
it remains to estimate the $p$-th roots of the quotients:
\[
\frac{\widetilde{I}_p}{\widetilde{I}_0}
\,=\,
\big(
2n\,\mu(a)
\big)^p\,\,
\underbrace{
\sum_{1\leqslant i_1<\cdots<i_p\leqslant n}\,
\frac{1}{a_{i_1}}\,\cdots\,\frac{1}{a_{i_p}}}_{=:\,\,
\sigma_p(\frac{1}{a_1},\dots,\frac{1}{a_n})},
\]
which incorporate the $p$-th symmetric functions $\sigma_p$ of
the weight inverses $\frac{1}{a_i}$. 
We start by extracting the $p$-th root
of $\big(2n\,\mu(a)\big)^p$ easily:
\[
\aligned
\underset{1\leqslant p\leqslant n}{\max}\,
\sqrt[p]{\,
\frac{\widetilde{I}_p}{\widetilde{I}_0}}
&
\,\,=\,\,
2n\,\mu(a)
\cdot
\underset{1\leqslant p\leqslant n}{\max}\,
\sqrt[p]{
\sigma_p
\Big(
\frac{1}{a_1},\dots,\frac{1}{a_n}
\Big)}
\\
\explicationtext{Seen above}
\ \ \ \ \ \ \ \ \ \ \ \ \ \ \ \ \ \ \ \ \ \ \ \ \ \
&
\,\,\leqslant\,\,
2n\,
\frac{r}{(r-1)^2}\,
r^n
\cdot
\underset{1\leqslant p\leqslant n}{\max}\,
\sqrt[p]{
\sigma_p
\Big(
\frac{1}{a_1},\dots,\frac{1}{a_n}
\Big)}.
\endaligned
\]

\begin{Lemma}
One has:
\[
\underset{1\leqslant p\leqslant n}{\max}\,
\sqrt[p]{
\sigma_p
\Big(
\frac{1}{a_1},\dots,\frac{1}{a_n}
\Big)}
\,\,=\,\,
\sigma_1
\Big(
\frac{1}{a_1},\dots,\frac{1}{a_n}
\Big).
\]
\end{Lemma}

\proof
For positive real numbers $b_1, \dots, b_n > 0$, the renormalized
symmetric functions:
\[
\aligned
s_p
\big(b_1,\dots,b_n\big)
\,:=\,
&\,
\frac{1}{\binom{n}{p}}\,
\sum_{1\leqslant i_1<\cdots<i_p\leqslant n}\,
b_{i_1}\cdots b_{i_p}
\\
\,=\,
&\,
\frac{1}{\binom{n}{p}}\,
\sigma_p\,
\big(
b_1,\dots,b_p
\big),
\endaligned
\]
satisfy the classical {\sl Mac Laurin inequality:}
\[
s_1
\,\,\geqslant\,\,
\sqrt[2]{s_2}
\,\,\geqslant\,\,
\sqrt[3]{s_3}
\,\,\geqslant\,\,
\cdots\cdots
\,\,\geqslant\,\,
\sqrt[n]{s_n}.
\]
A modified version, useful to us, is:

\begin{Assertion}
A similar, less fine, inequality, holds before renormalization:
\[
\sigma_1
\,\,\geqslant\,\,
\sqrt[2]{\sigma_2}
\,\,\geqslant\,\,
\cdots\cdots
\,\,\geqslant\,\,
\sqrt[p]{\sigma_p}
\,\,\geqslant\,\,
\sqrt[p+1]{\sigma_{p+1}}
\,\,\geqslant\,\,
\cdots\cdots
\,\,\geqslant\,\,
\sqrt[n]{\sigma_n}.
\]
\end{Assertion}

\proof
For $1 \leqslant p \leqslant n-1$, we would deduce from Mac Laurin
what we want:
\[
\Bigg(
\sqrt[p]{\frac{\sigma_p}{\binom{n}{p}}}
\overset{\sf known}{\,\,\geqslant\,\,}
\sqrt[p+1]{\frac{\sigma_{p+1}}{\binom{n}{p+1}}}
\Bigg)
\ \ \ \ \
\Longrightarrow
\ \ \ \ \
\Big(
\sqrt[p]{\sigma_p}
\overset{\text{\bf ?}}{\,\,\geqslant}
\sqrt[p+1]{\sigma_{p+1}}
\Big),
\] 
provided it would be true that:
\[
\frac{\sqrt[p]{\binom{n}{p}}}{
\sqrt[p+1]{\binom{n}{p+1}}}
\overset{\text{\bf ?}}{\,\,\geqslant\,\,}
1
\eqno
{\scriptstyle{(\forall\,1\,\leqslant\,p\,\leqslant\,n-1)}}.
\]
We claim that such numerical
inequalities hold true. Indeed, from
the two visible minorations:
\[
\aligned
n\,(n-1)\,\cdots\,(n-p+1)
&
\,\,\geqslant\,\,
\big(n-p\big)^p,
\\
\big(p+1\big)^p
&
\,\,\geqslant\,\,
1\cdot 2\cdot\ldots\cdot p,
\endaligned
\]
comes:
\[
\frac{n\,(n-1)\,\cdots\,(n-p+1)}{1\cdot 2\cdot\ldots\cdot p}
\,\,\geqslant\,\,
\frac{(n-p)^p}{(p+1)^p},
\]
whence:
\[
\bigg[
\frac{n\,(n-1)\,\cdots\,(n-p+1)}{1\cdot 2\cdot\ldots\cdot p}
\bigg]^{p+1}
\,\,\geqslant\,\,
\bigg[
\frac{n\,(n-1)\,\cdots\,(n-p+1)}{
1\cdot 2\cdot\ldots\cdot p}
\cdot
\frac{(n-p)}{(p+1)}
\bigg]^{p},
\]
and this is exactly what we wanted:
\[
\binom{n}{p}^{\!p+1}
\,\,\geqslant\,\,
\binom{n}{p+1}^{\!p}.
\qedhere
\]
\endproof

Lastly, with $b_1 := \frac{1}{a_1}$, \dots, $b_n := \frac{1}{a_n}$,
we get:
\[
\sigma_1
\Big(
\frac{1}{a_1},\dots,\frac{1}{a_n}
\Big)
\,\,\geqslant\,\,
\sqrt[p]{\sigma_p
\Big(
\frac{1}{a_1},\dots,\frac{1}{a_n}
\Big)}
\eqno
{\scriptstyle{(\forall\,1\,\leqslant\,p\,\leqslant\,n)}},
\]
which forces the maximum to be attained precisely when $p = 1$.
\endproof

So we obtain:
\[
\underset{1\leqslant p\leqslant n}{\max}\,
\sqrt[p]{\,
\frac{\widetilde{I}_p}{\widetilde{I}_0}}
\,\,\leqslant\,\,
2n\,
\frac{r}{(r-1)^2}\,
r^n
\cdot
\sigma_1
\Big(
\frac{1}{a_1},\dots,\frac{1}{a_n}
\Big).
\]
and it only remains to estimate:
\[
\aligned
\sigma_1
\Big(
\frac{1}{a_1},\dots,\frac{1}{a_{n-1}},\,\frac{1}{a_n}
\Big)
&
\,\,\leqslant\,\,
\frac{1}{r^{n-1}}
+\cdots+
\frac{1}{r}
+
1
\\
&
\,\,\leqslant\,\,
\frac{r}{r-1},
\endaligned
\]
in order to finish the proof of Proposition~{\ref{Proposition-r-3-n}}:
\begin{align}
\max\,
\big\vert
{\sf roots}
\big\vert
&
\,\,\leqslant\,\,
\big(
10\,n+2
\big)\,
\Big(
1+\frac{3}{r}
\Big)^{n-1}
\cdot
\underset{1\leqslant p\leqslant n}{\max}\,
\sqrt[p]{\,
\frac{\widetilde{I}_p}{\widetilde{I}_0}}\,\,
\notag
\\
&
\,\,\leqslant\,\,
\big(
10\,n+2
\big)\,
\Big(
1+\frac{3}{r}
\Big)^{n-1}
\cdot
2n\,
\frac{r}{(r-1)^2}\,
r^n
\cdot
\sigma_1
\Big(
\frac{1}{a_1},\dots,\frac{1}{a_n}
\Big)
\notag
\\
&
\,\,\leqslant\,\,
\big(
10\,n+2
\big)\,
\frac{(r+3)^{n-1}}{r^{n-1}}
\cdot
2n\,
\frac{r}{(r-1)^2}\,
r^n
\cdot
\frac{r}{r-1}
\notag
\\
&
\,\,=\,\,
\big(20\,n^2+4\,n\big)
\cdot
\frac{r^3}{(r-1)^3\,(r+3)}
\cdot
\big(
r+3
\big)^n
\notag
\\
&
\,\,\leqslant\,\,
25\,n^2
\cdot
\big(r+3\big)^n.
\qedhere
\end{align}
\endproof

\Section{\bf From Coordinates $(t_1, t_2, \dots, t_n)$ 
to Coordinates $(w_2, \dots, w_n)$}
\label{coordinates-t-coordinates-w}
\HEAD{{\ref{coordinates-t-coordinates-w}}.~{\sf From Coordinates 
$(t_1, t_2, \dots, t_n)$ to Coordinates $(w_2, \dots, w_n)$}
}{
Jo\"el {\sc Merker} and The-Anh Ta,
D\'epartement de Math\'ematiques d'Orsay, 
Universit\'e Paris-Sud, France}

The goal of this section is to transform 
both the product $C(t)$ and the $n^2$-power $f_0(t)$ 
into more tractable expressions, by introducing the formal variables:
\[
w_2
\,:=\,
\frac{t_1}{t_2},\ \ \ \ \
w_3
\,:=\,
\frac{t_2}{t_3},\ \ \ \ \
\dots\dots,\ \ \ \ \
w_n
\,:=\,
\frac{t_{n-1}}{t_n}.
\]

To enhance intuition, 
start by expanding the writing of the factors of two types
in the considered double big product:
\[
\aligned
C(t_1,\dots,t_n)
\,\,=\,\,
\frac{t_2-t_1}{t_2-2\,t_1}\,
\frac{t_3-t_1}{t_3-2\,t_1}\,
\cdots\,
\frac{t_n-t_1}{t_n-2\,t_1}
\\
\frac{t_3-t_2}{t_3-2\,t_2}
\cdots
\frac{t_n-t_2}{t_n-2\,t_2}
\\
\ddots
\ \ \ \ \ \ 
\vdots
\ \ \ \ \ \ \
\\
\frac{t_n-t_{n-1}}{t_n-2\,t_{n-1}}
\\
\frac{t_3-2\,t_2}{t_3-2\,t_2+t_1}
\cdots
\frac{t_n-2\,t_2}{t_n-2\,t_2+t_1}
\\
\ddots
\ \ \ \ \ \ \ \ \ \ \ \ \ \ 
\vdots
\ \ \ \ \ \ \ \ \
\\
\frac{t_n-2\,t_{n-1}}{t_n-2\,t_{n-1}+t_{n-2}}.
\endaligned
\]
To pass to the new variables, compute first for instance:
\[
\aligned
\frac{t_2-t_1}{t_2-2\,t_1}
&
\,=\,
\frac{1-\frac{t_1}{t_2}}{1-2\,\frac{t_1}{t_2}}
\,=\,
\frac{1-w_2}{1-2\,w_2},
\\
\frac{t_3-t_1}{t_3-2\,t_1}
&
\,=\,
\frac{1-\frac{t_1}{t_3}}{1-2\,\frac{t_1}{t_3}}
\,=\,
\frac{1-\frac{t_1}{t_2}\,\frac{t_2}{t_3}}{
1-2\,\frac{t_1}{t_2}\,\frac{t_2}{t_3}}
\,=\,
\frac{1-w_2w_3}{1-2\,w_2w_3},
\\
\frac{t_5-2\,t_2}{t_5-2\,t_2+t_1}
&
\,=\,
\frac{1-2\,\frac{t_2}{t_3}\,\frac{t_3}{t_4}\,\frac{t_4}{t_5}}{
1-2\,\frac{t_2}{t_3}\,\frac{t_3}{t_4}\,\frac{t_4}{t_5}
+\frac{t_1}{t_2}\,\frac{t_2}{t_3}\,\frac{t_3}{t_4}\,\frac{t_4}{t_5}}
\,=\,
\frac{1-2\,w_3w_4w_5}{1-2\,w_3w_4w_5+w_2w_3w_4w_5}.
\endaligned
\]
Generally, with as above:
\[
w_i
\,:=\,
\frac{t_{i-1}}{t_i}
\eqno
{\scriptstyle{(2\,\leqslant\,i\,\leqslant\,n)}},
\]
we can transform all the factors of first type, for indices
$2 \leqslant i \leqslant j \leqslant n$\,\,---\,\,mind
the shift $i \longmapsto i-1$ from the original definition of
$C(t)$:
\[
\aligned
E_{i,j}(t)
\,:=\,
\frac{t_j-t_{i-1}}{t_j-2\,t_{i-1}}
\,=\,
\frac{1-\frac{t_{i-1}}{t_j}}{1-2\,\frac{t_{i-1}}{t_j}}
&
\,=\,
\frac{1-\frac{t_{i-1}}{t_i}\cdots\frac{t_{j-1}}{t_j}}{
1-2\,\frac{t_{i-1}}{t_i}\cdots\frac{t_{j-1}}{t_j}}
\\
&
\,=\,
\frac{1-w_i\cdots w_j}{1-2\,w_i\cdots w_j}
\\
&
\,=:\,
E_{i,j}(w),
\endaligned
\]
Similarly, for $3 \leqslant i \leqslant j \leqslant n$, again
with the shift $i \longmapsto i-1$:
\[
\aligned
F_{i,j}(t)
\,:=\,
\frac{t_j-2\,t_{i-1}}{t_j-2\,t_{i-1}+t_{i-2}}
\,=\,
\frac{1-2\,\frac{t_{i-1}}{t_j}}{1-2\,\frac{t_{i-1}}{t_j}+
\frac{t_{i-2}}{t_j}}
&
\,=\,
\frac{1-2\,\frac{t_{i-1}}{t_i}\cdots\frac{t_{j-1}}{t_j}}{
1-2\,\frac{t_{i-1}}{t_i}\cdots\frac{t_{j-1}}{t_j}
+\frac{t_{i-2}}{t_{i-1}}\frac{t_{i-1}}{t_i}\cdots\frac{t_{j-1}}{t_j}}
\\
&
\,=\,
\frac{1-2\,w_i\cdots w_j}{
1-2\,w_i\cdots w_j+w_{i-1}w_i\cdots w_j}
\\
&
\,=:\,
F_{i,j}(w).
\endaligned
\]
Consequently:
\[
\aligned
C(t_1,t_2,\dots,t_n)
\,=\,
C(w_2,\dots,w_n)
\,:=\,
\frac{1-w_2}{1-2\,w_2}\,
\frac{1-w_2w_3}{1-2\,w_2w_3}\,
\cdots\,
\frac{1-w_2w_3\cdots w_n}{1-2\,w_2w_3\cdots w_n}
\\
\ \ \ \ \ \ \ \ 
\frac{1-w_3}{1-2\,w_3}\,
\cdots\,
\frac{1-w_3\cdots w_n}{1-2\,w_3\cdots w_n}
\\
\ddots
\ \ \ \ \ \ \ \ \ \ \
\vdots
\ \ \ \ \ \ \ \ \ \
\\
\frac{1-w_n}{1-2\,w_n}
\\
\frac{1-2\,w_3}{1-2\,w_3+w_2w_3}\,
\cdots\,
\frac{1-2\,w_3\cdots w_n}{
1-2\,w_3\cdots w_n+w_2w_3\cdots w_n}
\\
\ddots
\ \ \ \ \ \ \ \ \ \ \
\vdots
\ \ \ \ \ \ \ \ \ \
\\
\frac{1-2\,w_n}{1-2\,w_n+w_{n-1}w_n}.
\endaligned
\]
This can be abbreviated as:
\[
\aligned
C(w)
&
\,=\,
\prod_{2\leqslant i\leqslant j\leqslant n}\,
\frac{1-w_i\cdots w_j}{1-2\,w_i\cdots w_j}\,\,
\prod_{3\leqslant i\leqslant j\leqslant n}\,
\frac{1-2\,w_i\cdots w_j}{1-2\,w_i\cdots w_j+w_{i-1}w_i\cdots w_j}
\\
&
\,=\,
\prod_{2\leqslant i\leqslant j\leqslant n}\,
E_{i,j}
\prod_{3\leqslant i\leqslant j\leqslant n}\,
F_{i,j}.
\endaligned
\]

As is visible\,\,---\,\,and as was already visible before in variables
$(t_1, \dots, t_n)$\,\,---, 
the terms $1 - 2\, w_i \cdots w_j$ 
that appear in the denominators of the $E_{i,j}$ 
cancel out 
with the same terms appearing in the numerators of the $F_{i,j}$,
though only for $3 \leqslant i \leqslant j \leqslant n$. 
These simplifications conduct to the shorter representation:
\[
\aligned
C(w_2,\dots,w_n)
\,:=\,
\frac{1-w_2}{1-2\,w_2}\,
\ \ \ \ \ \
\frac{1-w_2w_3}{1-2\,w_2w_3}\,
\ \ \ \ \ \
\cdots\cdots\cdots\cdots\cdots\,
\ \ \ \ \ \
\frac{1-w_2w_3\cdots w_n}{1-2\,w_2w_3\cdots w_n}
\\
\frac{1-w_3}{1-2\,w_3+w_2w_3}\,
\ \
\cdots\cdots\,
\ \
\frac{1-w_3\cdots w_n}{
1-2\,w_3\cdots w_n+w_2w_3\cdots w_n}
\\
\ddots
\ \ \ \ \ \ \ \ \ \ \ \ \ \ \ \ \ \ \ \ \ \ \ \ \ \ \ \ \ \ \ \ \
\ \ \ \ \ \ 
\vdots
\ \ \ \ \ \ \ \ \ \ \ \ \ \ \ \ \ \
\\
\ddots
\ \ \ \ \ \ \ \ \ \ \ \ \ \ \ \ \ \ \ \ \ \ \ \ \ \ \ \ \ 
\vdots
\ \ \ \ \ \ \ \ \ \ \ \ \ \ \ \ \ \
\\
\frac{1-w_n}{1-2\,w_n+w_{n-1}w_n},
\endaligned
\]
which can be abbreviated as:
\[
\aligned
C(w_2,\dots,w_n)
\,=\,
E_2'(w_2)\,
E_3'(w_2,w_3)\,
\cdots
E_n'(w_2,w_3,\dots,w_n)
\\
F_{3,3}'(w_2,w_3)
\cdots
F_{3,n}'(w_2,w_3,\dots,w_n)
\\
\ddots
\ \ \ \ \ \ \ \ \ \ \ \ \ \ \ \ \ \ \ \ \ \ \ 
\vdots
\ \ \ \ \ \ \ \ \ \ \ \ \ \ \
\\
F_{n,n}'(w_{n-1},w_n),
\endaligned
\]
that is to say:
\[
C(w_2,\dots,w_n)
\,=\,
\prod_{2\leqslant j\leqslant n}\,
\underbrace{
\frac{1-w_2\cdots w_j}{1-2\,w_2\cdots w_j}}_{=:\,\,E_j'}\,\,
\prod_{3\leqslant i\leqslant j\leqslant n}\,
\underbrace{
\frac{1-w_i\cdots w_j}{
1-2\,w_i\cdots w_j+w_{i-1}w_i\cdots w_j}
}_{=:\,\,F_{i,j}'}.
\]

Next, let us re-express in the $w_i$ variables:
\[
\aligned
f_0(t)
&
\,=\,
\Big(
a_1t_1
+\cdots+
a_{n-2}t_{n-2}
+
a_{n-1}t_{n-1}
+
a_nt_n
\Big)^{n^2}
\\
&
\,=\,
\Big(
a_1
{\textstyle{\frac{t_1}{t_n}}}
+\cdots+
a_{n-2}
{\textstyle{\frac{t_{n-2}}{t_n}}}
+
a_{n-1}
{\textstyle{\frac{t_{n-1}}{t_n}}}
+
a_n
{\textstyle{\frac{t_n}{t_n}}}
\Big)^{n^2}\,
\big(
t_n
\big)^{n^2}
\\
&
\,=\,
\Big(
r^{n-1}w_2\cdots w_n
+\cdots+
r^2w_{n-1}w_n
+
r\,w_n
+
1
\Big)^{n^2}\,
\big(t_n\big)^{n^2}.
\endaligned
\]
To yet transform $t_n^{n^2}$ at the end, observe that:
\[
\aligned
\frac{1}{(w_2)^n(w_3)^{2n}\cdots (w_{n-1})^{n^2-2n}(w_n)^{n^2-n}}
&
\,=\,
\frac{1}{
\big(
\frac{t_1}{t_2}
\big)^n\,
\big(
\frac{t_2}{t_3}
\big)^{2n}
\cdots
\big(
\frac{t_{n-2}}{t_{n-1}}
\big)^{n^2-2n}
\big(
\frac{t_{n-1}}{t_n}
\big)^{n^2-n}}
\\
&
\,=\,
\frac{1}{
t_1^n\,t_2^n\,\cdots t_{n-2}^n\,t_{n-1}^n\,
\frac{1}{t_n^{n^2-n}}}
\\
&
\,=\,
\frac{t_n^{n^2}}{
t_1^n\,t_2^n\,\cdots t_{n-2}^n\,t_{n-1}^n\,t_n^n},
\endaligned
\]
whence:
\[
f_0(t)
\,=\,
\frac{
\big(
r^{n-1}w_2\cdots w_n
+\cdots+
r^2w_{n-1}w_n
+
r\,w_n
+
1
\big)^{n^2}}{
w_2^n\,w_3^{2n}\,\cdots\,
w_{n-1}^{n^2-2n}\,w_n^{n^2-n}}\,\,
\frac{1}{t_1^n\,t_2^n\,\cdots\,t_{n-1}^n\,t_n^n}.
\]

Consequently, in Problem~{\ref{Problem-I-0}},
the coefficient $I_0$ of the monomial $t_1^n \cdots t_n^n$
in the product $C(t) \cdot f_0(t)$ identifies with the
{\em constant} term, namely the coefficient of
$w_2^0 \cdots w_n^0 = 1$, in the product:
\[
\footnotesize
\!\!\!\!\!\!\!\!\!\!\!\!\!\!\!\!\!\!\!\!
\aligned
\frac{
\big(
r^{n-1}w_2\cdots w_n
+\cdots+
r^2w_{n-1}w_n
+
r\,w_n
+
1
\big)^{n^2}}{
w_2^n\,w_3^{2n}\,\cdots\,
w_{n-1}^{n^2-2n}\,w_n^{n^2-n}}
\,\,\cdot\,\,
\frac{1-w_2}{1-2\,w_2}\,
\ \ \ \ \ \
\frac{1-w_2w_3}{1-2\,w_2w_3}\,
\ \ \ \ \ \
\cdots\cdots\cdots\cdots\cdots\,
\ \ \ \ \ \
\frac{1-w_2w_3\cdots w_n}{1-2\,w_2w_3\cdots w_n}
\\
\frac{1-w_3}{1-2\,w_3+w_2w_3}\,
\ \
\cdots\cdots\,
\ \
\frac{1-w_3\cdots w_n}{
1-2\,w_3\cdots w_n+w_2w_3\cdots w_n}
\\
\ddots
\ \ \ \ \ \ \ \ \ \ \ \ \ \ \ \ \ \ \ \ \ \ \ \ \ \ \ \ \ \ \ \ \
\ \ \ \ \ \ 
\vdots
\ \ \ \ \ \ \ \ \ \ \ \ \ \ \ \ \ \
\\
\ddots
\ \ \ \ \ \ \ \ \ \ \ \ \ \ \ \ \ \ \ \ \ \ \ \ \ \ \ \ \ 
\vdots
\ \ \ \ \ \ \ \ \ \ \ \ \ \ \ \ \ \
\\
\frac{1-w_n}{1-2\,w_n+w_{n-1}w_n}.
\endaligned
\]

It is now appropriate to expand the $n^2$ power in the numerator
above plainly as:
\[
\frac{1}{(w_2)^n\cdots (w_n)^{n^2-2n}}\,
\sum_{
i_1,\dots,i_n\geqslant 0
\atop
i_1+\cdots+i_n=n^2}\,
(1)^{i_1}\,(rw_n)^{i_2}
\big(r^2w_{n-1}w_n\big)^{i_3}
\cdots
\big(r^{n-1}w_2\cdots w_n\big)^{i_n}\,
\frac{(n^2)!}{i_1!\,i_2!\,i_3!\,\cdots i_n!}
\]

Next, we would like to point out 
that $C(w_2, \dots, w_n)$ is a product of
rational expressions which expand all in converging power series at
the origin. More precisely, using the trivial expansion:
\[
\aligned
E(x)
\,:=\,
\frac{1-x}{1-2\,x}
&
\,=\,
1
+
\frac{x}{1-2\,x}
\\
&
\,=\,
1
+
\sum_{i=1}^\infty\,
2^{i-1}\,
x^i,
\endaligned
\]
together with the
expansion of 
Lemma~{\ref{Lemma-expansion-F-F-hat}}\,\,---\,\,with 
the convention that $\binom{\ell-1}{-1} = 0 = 
\binom{\ell-1}{\ell}$\,\,---\,\,:
\[
\aligned
F(x,y)
\,:=\,
\frac{1-y}{1-2\,y+x\,y}
&
\,=\,
1
+
\frac{y-x\,y}{1-2\,y+x\,y}
\\
&\,
\,=\,
1
+
\sum_{\ell=1}^\infty\,
y^\ell\,
\sum_{0\leqslant k\leqslant\ell}\,
(-1)^k\,
x^k\,
\Big[
2^{\ell-1-k}\,
{\textstyle{\binom{\ell-1}{k}}}
+
2^{\ell-k}\,
{\textstyle{\binom{\ell-1}{k-1}}}
\Big],
\endaligned
\]
and re-expressing:
\[
C(w_2,\dots,w_n)
\,=\,
\prod_{2\leqslant j\leqslant n}\,
E\big(w_2\cdots w_j\big)\,\,
\prod_{3\leqslant i\leqslant j\leqslant n}\,
F\big(w_{i-1},\,w_i\cdots w_j\big),
\]
and lastly, multiplying all the obtained converging power series,
one can in principle receive an expansion: 
\[
C(w_2,\dots,w_n)
\,=\,
\sum_{k_2=0}^\infty
\cdots
\sum_{k_n=0}^\infty\,
C_{k_2,\dots,k_n}\,
(w_2)^{k_2}
\cdots
(w_n)^{k_n},
\]
which is holomorphic in a neighborhood of the origin. However, it is
very delicate to reach closed explicit expressions for these integer
Taylor coefficients $C_{k_2, \dots, k_n}$, a difficulty which lies at
the very core of Problem~{\ref{Problem-I-0}}.

In summary, the quantity $I_0$ we want to determine, in order to show
that it satisfies $I_0 \geqslant \widetilde{I}_0$, is the coefficient
of the constant term $w_2^0 \cdots w_n^0$ in a product consisting of
$2$ rows:
\[
\aligned
{}
&
\sum_{k_2=0}^\infty
\cdots
\sum_{k_n=0}^\infty\,
C_{k_2,\dots,k_n}\,
(w_2)^{k_2}
\cdots
(w_n)^{k_n}
\,\cdot
\\
&
\cdot
\frac{1}{(w_2)^n\cdots (w_n)^{n^2-n}}\,
\sum_{
i_1,\dots,i_n\geqslant 0
\atop
i_1+\cdots+i_n=n^2}\,
(1)^{i_1}\,(rw_n)^{i_2}
\big(r^2w_{n-1}w_n\big)^{i_3}
\cdots
\big(r^{n-1}w_2\cdots w_n\big)^{i_n}\,
\frac{(n^2)!}{i_1!\,i_2!\,i_3!\,\cdots i_n!}.
\endaligned
\]
Clearly, the second row becomes, after reorganization, a Laurent
series of the form:
\[
\sum_{-n\leqslant\ell_2}\,
\cdots
\sum_{-(n^2-n)\leqslant\ell_n}\,
J_{\ell_2,\dots,\ell_n}\,
(w_2)^{\ell_2}
\cdots
(w_n)^{\ell_n}.
\]
But because in the first row one always has $k_2, \dots, k_n \geqslant
0$, all Laurent monomials $(w_2)^{\ell_2} \cdots (w_n)^{\ell_n}$ in
the second row for which $\ell_i \geqslant 1$ for some $1 \leqslant i
\leqslant n$ do {\em not} contribute to the determination of the
desired constant term $w_2^0 \cdots w_n^0$. So the summation in the
second row can be truncated to:
\[
\sum_{-n\leqslant\ell_2\leqslant 0}\,
\cdots
\sum_{-(n^2-n)\leqslant\ell_2\leqslant 0}\,
J_{\ell_2,\dots,\ell_n}\,
(w_2)^{\ell_2}
\cdots
(w_n)^{\ell_n}.
\] 
A supplementary change of indices
followed by a reorganization conducts to 
an appropriate reformulation of what is $I_0$:
the following statement will then constitute the 
very starting point of our further explorations.

\begin{Proposition} 
One has:
\[
I_0
\,=\,
\big[w_2^0\cdots w_n^0\big]\,
\Big(
A(w_2,\dots,w_n)
\cdot
C(w_2,\dots,w_n)
\Big),
\]
where:
\[
\aligned
\!\!\!\!\!\!\!\!\!\!\!\!\!\!\!\!\!\!\!\!
A(w_2,\dots,w_n)
\,:=\,
\!\!\!\!\!\!\!\!\!\!\!\!\!\!\!
\sum_{
\substack{
0\leqslant k_2\leqslant n
\\
\ \ \ \ \ 
0\leqslant k_3\leqslant n+k_2
\\
\ \ \ \ \ \
\cdots\cdots\cdots\cdots\cdots\cdots\cdots
\\
\ \ \ \ \ \ \ \ \ \ \ 
0\leqslant k_{n-1}\leqslant n+k_{n-2}
\\
\ \ \ \ \ \ \ \ \ \ \
0\leqslant\,k_n\,\,\,\,\,\leqslant n+k_{n-1}
}}
\!\!\!\!\!
&
\frac{(n^2)!}{(n-k_2)!(n+k_2-k_3)!\cdots 
(n+k_{n-2}-k_{n-1})!(n+k_{n-1}-k_n)!(n+k_n)!}
\,\cdot
\\
&
\cdot
\frac{r^{\,n\frac{n(n-1)}{2}}}{r^{k_2+\cdots+k_n}}\,
\frac{1}{(w_2)^{k_2}\cdots (w_n)^{k_n}},
\endaligned
\]
and where $C(w_2,\dots,w_n)$ is as before.
\end{Proposition}

\proof
We therefore rewrite:
\[
\footnotesize
\!\!\!\!\!\!\!\!\!\!\!\!\!\!\!
\aligned
{}
&
\frac{1}{w_2^nw_3^{2n}\cdots w_{n-1}^{(n-2)n}w_n^{(n-1)n}}
\!\!\!
\sum_{
i_1,\dots,i_n\geqslant 0
\atop
i_1+\cdots+i_n=n^2}\,
(1)^{i_1}\,(r\,w_n)^{i_2}
\big(r^2w_{n-1}w_n\big)^{i_3}
\cdots
\big(r^{n-2}w_3\cdots w_{n-1}w_n\big)^{i_{n-1}}
\big(r^{n-1}w_2w_3\cdots w_{n-1}w_n\big)^{i_n}
\cdot
\\
&
\ \ \ \ \ \ \ \ \ \ \ \ \ \ \ \ \ \ \ \ \ \ \ \ \ \ \ \ \ \ \ \ \ \
\ \ \ \ \ \ \ \ \ \ \ \ \ \ \ \ \ \ \ \ \ \ \ \ \ \ \ \ \ \ \  
\cdot
\frac{(n^2)!}{i_1!\,i_2!\,i_3!\,\cdots i_{n-1}!i_n!}
\,\,=
\\
&
=:\,\,
\sum_{-n\leqslant\ell_2}\,
\sum_{-2n\leqslant\ell_3}\,
\cdots
\sum_{-(n-2)n\leqslant\ell_{n-1}}\,
\sum_{-(n-1)n\leqslant\ell_n}\,
J_{\ell_2,\ell_3,\dots,\ell_{n-1},\ell_n}\,
(w_2)^{\ell_2}
(w_3)^{\ell_3}
\cdots
(w_{n-1})^{\ell_{n-1}}
(w_n)^{\ell_n},
\endaligned
\]
so that the correspondence between exponents is:
\[
\aligned
-\,n
\,\leqslant\,
\ell_2
&
\,=\,
i_n
-
n,
\\
-\,2\,n
\,\leqslant\,
\ell_3
&
\,=\,
i_n
+
i_{n-1}
-
2\,n,
\\
\cdots
&
\cdots\cdots\cdots\cdots\cdots\cdots
\\
-\,(n-2)\,n
\,\leqslant\,
\ell_{n-1}
&
\,=\,
i_n+i_{n-1}
+\cdots+
i_3
-
(n-2)\,n,
\\
-\,(n-1)\,n
\,\leqslant\,
\ell_n
&
\,=\,
i_n+i_{n-1}
+\cdots+
i_3
+
i_2
-
(n-1)\,n.
\endaligned
\]
Performing the harmless truncations $\ell_2 \leqslant 0$, \dots,
$\ell_n \leqslant 0$ leads then to the inequalities:
\[
\aligned
0
&
\,\leqslant\,
i_n
\,\leqslant\,
n,
\\
0
&
\,\leqslant\,
i_n+i_{n-1}
\,\leqslant\,
2\,n,
\\
\cdots
&
\cdots\cdots\cdots\cdots\cdots\cdots\cdots\cdots\cdots
\\
0
&
\,\leqslant\,
i_n+i_{n-1}+\cdots+i_3
\,\leqslant\,
(n-2)\,n,
\\
0
&
\,\leqslant\,
i_n+i_{n-1}+\cdots+i_3+i_2
\,\leqslant\,
(n-1)\,n,
\endaligned
\]
so that it suffices to consider, before multiplying by $C(w_2,
\dots, w_n)$, the truncated series:
\[
\aligned
A
\,:=\,
\sum_{
\substack{
i_1+\cdots+i_n=n^2
\\
i_1\geqslant 0,\,\dots,\,i_n\geqslant 0
\\
0\leqslant i_n\leqslant n
\\
0
\leqslant i_n+i_{n-1}\leqslant 2\,n
\\
\cdots\cdots\cdots\cdots\cdots\cdots\cdots\cdots\cdots\cdots\cdots
\\
0
\leqslant i_n+i_{n-1}+\cdots+i_3\leqslant (n-2)\,n
\\
0
\leqslant i_n+i_{n-1}+\cdots+i_3+i_2\leqslant (n-1)\,n}}
&
\!\!\!\!\!
\frac{1}{w_2^{n-i_n}\,w_3^{2n-i_{n-1}-i_n}
\cdots
w_{n-1}^{(n-2)n-i_n-i_{n-1}-\cdots-i_3}
w_n^{(n-1)n-i_n-i_{n-1}-\cdots-i_3-i_2}}
\,\cdot\,
\\
&
\,\cdot\,
r^{i_2+2\,i_3+\cdots+(n-2)i_{n-1}+(n-1)i_n}
\,\cdot\,
\frac{(n^2)!}{i_1!\,i_2!\,i_3!\,\cdots\,i_{n-1}!\,i_n!}.
\endaligned
\]

To reach the expression shown by the proposition, introduce
the new nonnegative integer indices:
\reqnomode\usetagform{EngelLie}
\begin{align}
k_2
&
\,:=\,
n-i_n
\tag{(k_2\,\geqslant\,0),}
\\
k_3
&
\,:=\,
2\,n-i_{n-1}-i_n
\tag{(k_3\,\geqslant\,0),}
\\
\cdots
&
\cdots\cdots\cdots\cdots\cdots\cdots\cdots\cdots
\tag{(\cdots\cdots\cdots\!),}
\\
k_{n-1}
&
\,:=\,
(n-2)\,n
-i_3-\cdots-i_{n-1}-i_n
\tag{(k_{n-1}\,\geqslant\,0),}
\\
k_n
&
\,:=\,
(n-1)\,n
-i_2-i_3-\cdots-i_{n-1}-i_n
\tag{(k_n\,\geqslant\,0).}
\end{align}
To finish, three explanations are needed.

Firstly, one has the inequalities:
\[
\aligned
0
&
\,\leqslant\,
k_2
\,\leqslant\,
n,
\\
0
&
\,\leqslant\,
k_3
\,\leqslant\,
n+k_2,
\endaligned
\]
because $i_n \geqslant 0$ and because:
\[
k_3
\,=\,
n-i_{n-1}+n-i_n
\,=\,
\underbrace{n-i_{n-1}}_{i_{n-1}\geqslant 0}
+
k_2
\,\,\leqslant\,\,
n+k_2.
\]
Similarly:
\[
0
\,\leqslant\,
k_4
\,=\,
3\,n-i_{n-2}-i_{n-1}-i_n
\,=\,
n-i_{n-2}
+
k_3
\,\,\leqslant\,\,
n+k_3,
\]
and so on up to:
\[
\aligned
0
&
\,\leqslant\,
k_{n-1}
\,\leqslant\,
n+k_{n-2},
\\
0
&
\,\leqslant\,
k_n
\,\leqslant\,
n+k_{n-1}.
\endaligned
\]

Secondly, since:
\[
k_2+k_3+\cdots+k_{n-1}+k_n
\,=\,
n\,\big(
1+2+\cdots+(n-2)+(n-1)
\big)
-
i_2-2\,i_3-\cdots-(n-2)\,i_{n-1}-(n-1)\,i_n,
\]
the exponent of $r$ becomes:
\[
i_2+2\,i_3+\cdots+(n-2)\,i_{n-1}+(n-1)\,i_n
\,=\,
n\,
{\textstyle{\frac{n(n-1)}{2}}}
-k_2-k_3-\cdots-k_{n-1}-k_n.
\]

Thirdly and lastly, the factorials become:
\[
\aligned
i_n!
&
\,=\,
(n-k_2)!,
\\
i_{n-1}!
&
\,=\,
\big(n+n-i_n-k_3\big)!
\,=\,
\big(n+k_2-k_3\big)!,
\\
\cdot\cdots
&
\cdots\cdots\cdots\cdots\cdots\cdots\cdots\cdots\cdots
\cdots\cdots\cdots\cdots\cdots\cdots
\\
i_2!
&
\,=\,
\big(
-\,k_n
+
(n-1)\,n
-
i_3-\cdots-i_n
\big)!
\,=\,
\big(n+k_{n-1}-k_n\big)!,
\\
i_1!
&
\,=\,
\big(
n\,n
-
(n+k_{n-1}-k_n)
-\cdots-
(n+k_2-k_3)
-
(n-k_2)
\big)
\,=\,
\big(n+k_n\big)!.
\endaligned
\]
These three explanations yield the expression of
$A(w_2, \dots, w_n)$ stated by the proposition.
\endproof

Next, because only the quotient $\frac{I_0}{\widetilde{I}_0}$ must be
studied in order to reach the minoration $I_0 \geqslant
\widetilde{I}_0$, we can divide everything in advance by
the central monomial: 
\[
\widetilde{I}_0
\,=\,
\frac{(n^2)!}{n!\cdots\,n!}\,
r^{\,n\frac{n(n-1)}{2}}.
\]
Equivalently, we factor:
\[
\aligned
A
\,=\,
&\,
\frac{(n^2)!}{n!\cdots\,n!}\,
r^{\,n\frac{n(n-1)}{2}}
\,\cdot
\\
&\,
\cdot\,
\sum_{
\substack{
0\leqslant k_2\leqslant n
\\
\ \ \ \ \ 
0\leqslant k_3\leqslant n+k_2
\\
\ \ \ \ \ \
\cdots\cdots\cdots\cdots\cdots\cdots\cdots
\\
\ \ \ \ \ \ \ \ \ \ \ 
0\leqslant k_{n-1}\leqslant n+k_{n-2}
\\
\ \ \ \ \ \ \ \ \ \ \
0\leqslant\,k_n\,\,\,\,\,\leqslant n+k_{n-1}
}}
\!\!\!\!\!\!\!\!\!\!\!\!\!\!\!
\frac{n!}{(n-k_2)!}\,
\frac{n!}{(n+k_2-k_3)!}\,
\cdots\,
\frac{n!}{(n+k_{n-2}-k_{n-1})!}\,
\frac{n!}{(n+k_{n-1}-k_n)!(n+k_n)!}
\,\cdot
\\
&
\cdot
\frac{1}{r^{k_2+k_3+\cdots+k_{n-1}+k_n}}\,
\frac{1}{w_2^{k_2}w_3^{k_3}\cdots w_{n-1}^{k_{n-1}}w_n^{k_n}},
\endaligned
\]
we keep the same name $A$ after eliminating the
factor $\widetilde{I}_0$ on the first line, 
and we reformulate our goal as a more precise

\begin{Problem}
\label{Problem-A-C}
{\sl For some specific choice of a fixed constant $r \geqslant 3$,
to show that for any $n \geqslant 2$, 
the coefficient of the constant
monomial $w_2^0 \cdots w_n^0$ in the product
$C(w) \cdot A(w)$ is at least equal to $1$, namely:}
\[
1
\,\,\leqslant\,\,
\big[w_2^0\cdots w_n^0\big]
\Big(
C(w_2,\dots,w_n)
\cdot
A(w_2,\dots,w_n)
\Big),
\]
{\sl where:}
\[
\aligned
C(w_2,\dots,w_n)
\,:=\,
\frac{1-w_2}{1-2\,w_2}\,
\frac{1-w_2w_3}{1-2\,w_2w_3}\,
\cdots\,
\frac{1-w_2w_3\cdots w_n}{1-2\,w_2w_3\cdots w_n}
\\
\ \ \ \ \ \ \ \ 
\frac{1-w_3}{1-2\,w_3}\,
\cdots\,
\frac{1-w_3\cdots w_n}{1-2\,w_3\cdots w_n}
\\
\ddots
\ \ \ \ \ \ \ \ \ \ \
\vdots
\ \ \ \ \ \ \ \ \ \
\\
\frac{1-w_n}{1-2\,w_n}
\\
\frac{1-2\,w_3}{1-2\,w_3+w_2w_3}\,
\cdots\,
\frac{1-2\,w_3\cdots w_n}{
1-2\,w_3\cdots w_n+w_2w_3\cdots w_n}
\\
\ddots
\ \ \ \ \ \ \ \ \ \ \
\vdots
\ \ \ \ \ \ \ \ \ \
\\
\frac{1-2\,w_n}{1-2\,w_n+w_{n-1}w_n}.
\endaligned
\]
{\sl and where:}
\[
\aligned
\!\!\!\!\!\!\!\!\!\!\!\!\!\!\!
A(w_2,\dots,w_n)
\,:=\,
&\,
\!\!\!\!\!\!\!\!\!\!\!\!\!\!\!
\sum_{
\substack{
0\leqslant k_2\leqslant n
\\
\ \ \ \ \ 
0\leqslant k_3\leqslant n+k_2
\\
\ \ \ \ \ \
\cdots\cdots\cdots\cdots\cdots\cdots\cdots
\\
\ \ \ \ \ \ \ \ \ \ \ 
0\leqslant k_{n-1}\leqslant n+k_{n-2}
\\
\ \ \ \ \ \ \ \ \ \ \
0\leqslant\,k_n\,\,\,\,\,\leqslant n+k_{n-1}
}}
\!\!\!\!\!\!\!\!\!\!
\frac{n!}{(n-k_2)!}\,
\frac{n!}{(n+k_2-k_3)!}\,
\cdots\,
\frac{n!}{(n+k_{n-2}-k_{n-1})!}\,
\frac{n!}{(n+k_{n-1}-k_n)!(n+k_n)!}
\,\cdot
\\
&
\ \ \ \ \ \ \ \ \ \ \ \ \ \ \ \ \ \ \ \ \ \ \ \ \ 
\cdot
\frac{1}{r^{k_2+k_3+\cdots+k_{n-1}+k_n}}\,
\frac{1}{w_2^{k_2}w_3^{k_3}\cdots w_{n-1}^{k_{n-1}}w_n^{k_n}}.
\endaligned
\]
\end{Problem}

Of course, under the hypothesis that the power series expansion
of $C(w)$ is known:
\[
C(w_2,\dots,w_n)
\,=\,
\sum_{k_2=0}^\infty
\cdots
\sum_{k_n=0}^\infty\,
C_{k_2,\dots,k_n}\,
(w_2)^{k_2}
\cdots
(w_n)^{k_n},
\]
the coefficient in question writes up as the sum:
\[
\aligned
\CA_{n-1}^n
\,:=\,
&\,
\!\!\!\!\!\!\!\!\!\!\!\!\!\!\!
\sum_{
\substack{
0\leqslant k_2\leqslant n
\\
\ \ \ \ \ 
0\leqslant k_3\leqslant n+k_2
\\
\ \ \ \ \ \
\cdots\cdots\cdots\cdots\cdots\cdots\cdots
\\
\ \ \ \ \ \ \ \ \ \ \ 
0\leqslant k_{n-1}\leqslant n+k_{n-2}
\\
\ \ \ \ \ \ \ \ \ \ \
0\leqslant\,k_n\,\,\,\,\,\leqslant n+k_{n-1}
}}
\!\!\!\!\!\!\!\!\!\!
\frac{n!}{(n-k_2)!}\,
\frac{n!}{(n+k_2-k_3)!}\,
\cdots\,
\frac{n!}{(n+k_{n-2}-k_{n-1})!}\,
\frac{n!}{(n+k_{n-1}-k_n)!(n+k_n)!}
\,\cdot
\\
&
\ \ \ \ \ \ \ \ \ \ \ \ \ \ \ \ \ \ \ \ \ \ \ \ \ 
\cdot
\frac{1}{r^{k_2+k_3+\cdots+k_{n-1}+k_n}}\,
C_{k_2,k_3,\dots,k_{n-1},k_n},
\endaligned
\]
which should satisfy:
\[
\CA_{n-1}^n
\overset{\text{\bf ?}}{\,\,\geqslant\,\,}
1
\eqno
{\scriptstyle{(\forall\,n\,\geqslant\,2)}}.
\]

\Section{\bf Approximations of multinomial quotients 
$M_{k_2, \dots, k_n}^n$}
\label{approximations-M-k-2-k_n}
\HEAD{{\ref{approximations-M-k-2-k_n}}.~Approximations 
of multinomial quotients $M_{k_2, \dots, k_n}^n$
}{
Jo\"el {\sc Merker} and The-Anh Ta,
D\'epartement de Math\'ematiques d'Orsay, 
Universit\'e Paris-Sud, France}

Let us attribute a name to the quotients of multinomial
coefficients which have appeared above:
\[
\aligned
M_{k_2,k_3,\dots,k_{n-1},k_n}^n
\,:=\,
&\,
\frac{n!}{(n-k_2)!}\,\,
\frac{n!}{(n+k_2-k_3)!}\,\,
\cdots\cdots\,\,
\frac{n!}{(n+k_{n-1}-k_n)!}\,\,
\frac{n!}{(n+k_n)!}
\\
\,=\,
&\,
\frac{
\frac{(n^2)!}{
(n-k_2)!\,(n+k_2-k_3)!\,\cdots\,(n+k_{n-1}-k_n)!\,(n+k_n)!}
}{
\frac{(n^2)!}{n!\,n!\,\cdots\,n!\,n!}
}.
\endaligned
\]
When $k_2 = k_3 = \cdots = k_{n-1} = k_n = 0$, this is just:
\[
M_{0,0,\dots,0,0}^n
\,=\,
1.
\]

\begin{Lemma}
\label{Lemma-M-less-1}
For all indices $(k_2, k_3, \dots, k_{n-1}, k_n) \neq
(0, 0, \dots, 0, 0)$ in the domain:
\[
\aligned
0
&
\,\leqslant\,
k_2
\,\leqslant\,
n,
\\
0
&
\,\leqslant\,
k_3
\,\leqslant\,
n+k_2,
\\
\cdot\cdot
&
\cdots\cdots\cdots\cdots\cdots\cdots
\\
0
&
\,\leqslant\,
k_{n-1}
\,\leqslant\,
n+k_{n-2},
\\
0
&
\,\leqslant\,
k_n
\ \ \
\,\leqslant\,
n+k_{n-1},
\endaligned
\]
there are strict inequalities:
\[
\big(
0
\,\leqslant\,
\big)
\ \ \ \ \ \ \ \ \
M_{k_2,k_3,\dots,k_{n-1},k_n}^n
\,\,<\,\,
1, 
\]
with equality $=1$ only when $k_2 = k_3 = \cdots = k_{n-1} = k_n
= 0$.
\end{Lemma}

\proof
Coming back to the old (nonnegative) indices:
\[
\aligned
i_n
&
\,=\,
n
-
k_2,
\\
i_{n-1}
&
\,=\,
n+k_2-k_3,
\\
\cdots
&
\cdots\cdots\cdots\cdots\cdots\cdot
\\
i_2
&
\,=\,
n+k_{n-1}-k_n,
\\
i_1
&
\,=\,
n+k_n,
\endaligned
\]
which satisfy $i_1 + i_2 + \cdots + i_{n-1} + i_n = n\,n = n^2$
and are not all equal to $n$\,\,---\,\,otherwise all
$k_\lambda = 0$\,\,---, we have to explain the inequalities:
\[
\frac{n!}{i_1!}\,\,
\frac{n!}{i_2!}\,\,
\cdots\cdots\,\,
\frac{n!}{i_{n-1}!}\,\,
\frac{n!}{i_n!}
\overset{\text{\bf ?}}{\,\,\,<\,\,}
1.
\]

After a reordering, we can assume that:
\[
i_1
\neq
n,\,\,\,
\dots\dots,\,\,\,
i_\kappa
\neq
n,\,\,\,
i_{\kappa+1}
=
n,\,\,\,
\dots\dots,\,\,\,
i_n
=
n,
\]
for a certain integer $1 \leqslant \kappa \leqslant n$. 
Since the factors $\frac{n!}{n!} = 1$ have no effect,
we are led to ask whether:
\[
\frac{n!}{i_1!}\,\,
\cdots\cdots\,\,
\frac{n!}{i_\lambda!}\,\,
\cdots\cdots\,\,
\frac{n!}{i_\kappa!}
\overset{\text{\bf ?}}{\,\,\,<\,\,}
1.
\]

Observing that:
\[
i_1
+\cdots+
i_\lambda
+\cdots+
i_\kappa
\,=\,
\kappa\,n,
\]
let us distinguish two cases about these $i_\lambda$
for every $1 \leqslant \lambda \leqslant \kappa$:
\[
i_\lambda
\,<\,
n
\ \ \ \ \ \ \ \ \ \ \ \ \ \ \ \ \ \
\text{or}
\ \ \ \ \ \ \ \ \ \ \ \ \ \ \ \ \ \
i_\lambda
\,>\,
n.
\]
When $i_\lambda < n$, we simplify:
\[
\frac{n!}{i_\lambda!}
\,=\,
n\,\big(n-1\big)
\cdots
\big(i_\lambda+1\big),
\]
and when $i_\lambda > n$, we simplify:
\[
\frac{n!}{i_\lambda!}
\,=\,
\frac{1}{i_\lambda\,(i_\lambda-1)\cdots(n+1)},
\]
so that:
\[
\frac{n!}{i_1!}\,\,
\cdots\cdots\,\,
\frac{n!}{i_\lambda!}\,\,
\cdots\cdots\,\,
\frac{n!}{i_\kappa!}
\,\,=\,\,
\frac{
\prod_{i_\lambda<n}\,
n\,(n-1)\,\cdots\,(i_\lambda+1)
}{
\prod_{i_\lambda>n}\,
i_\lambda(i_\lambda-1)\cdots(n+1)}.
\]

Now, we observe that in this
fraction the number of integer factors at numerator
place is equal to the number of integer factors at 
denominateur place, because the equality above:
\[
\kappa\,n
\,=\,
\sum_{1\leqslant\lambda\leqslant\kappa}\,
i_\lambda
\,=\,
\sum_{i_\lambda<n}\,
i_\lambda
+
\sum_{i_\lambda>n}\,
i_\lambda
\]
can be rewritten as:
\[
\sum_{i_\lambda<n}\,
\big(n-i_\lambda\big)
\,=\,
\sum_{i_\lambda>n}\,
\big(
i_\lambda-n
\big).
\]
But {\em each} integer factor at denominator place is 
{\em larger} than
{\em all} integer factors at numerator place, 
so the fraction must be $<1$.
\endproof

Visibly, in the quantity under study:
\[
M_{k_2,k_3,\dots,k_{n-1},k_n}^n
\,=\,
\frac{n!}{(n-k_2)!}\,\,
\frac{n!}{(n+k_2-k_3)!}\,\,
\cdots\cdots\,\,
\frac{n!}{(n+k_{n-1}-k_n)!}\,\,
\frac{n!}{(n+k_n)!},
\]
there are two types of quotients:
\[
\frac{n!}{(n-k)!}
\ \ \ \
\text{\rm with}\,\,\,
k\geqslant 0
\ \ \ \ \ \ \ \ \ \ \ \ \ \ \ \ \ \
\text{and}
\ \ \ \ \ \ \ \ \ \ \ \ \ \ \ \ \ \
\frac{n!}{(n+\ell)!}
\ \ \ \
\text{\rm with}\,\,\,
\ell\geqslant 0.
\]
We can simplify, factorize, and rewrite the first type
quotients as:
\[
\aligned
\frac{n!}{(n-k)!}
\,=\,
\frac{n\,(n-1)\,\cdots\,(n-k+1)}{1}
&
\,=\,
n^k\,
\Big(
1
-
\frac{0}{n}
\Big)\,
\Big(
1
-
\frac{1}{n}
\Big)\,
\cdots
\Big(
1
-
\frac{k-1}{n}
\Big)
\\
&
\,=\,
n^k\,
\prod_{0\leqslant i\leqslant k-1}\,
\Big(
1
-
\frac{i}{n}
\Big),
\endaligned
\]
and the second type quotients as:
\[
\aligned
\frac{n!}{(n+\ell)!}
\,=\,
\frac{1}{(n+\ell)\cdots(n+1)}
&
\,=\,
\frac{1}{
\big(1+\frac{\ell}{n}\big)\,
\cdots
\big(1+\frac{1}{n}\big)\,n^\ell}
\\
&
\,=\,
n^{-\ell}\,
\prod_{1\leqslant j\leqslant\ell}\,
\Big(
1
+
\frac{j}{\ell}
\Big)^{-1}.
\endaligned
\]

In order to estimate the proximity to $1$ of these products,
let us take their logarithms:
\[
\aligned
\log\,
\prod_{0\leqslant i\leqslant k-1}\,
\Big(
1
-
\frac{i}{n}
\Big)
&
\,=\,
\log\,
\Big(
1
-
\frac{0}{n}
\Big)
+
\log\,
\Big(
1
-
\frac{1}{n}
\Big)
+\cdots+
\log\,
\Big(
1
-
\frac{k-1}{n}
\Big)
\\
&
\,\leqslant\,
0,
\endaligned
\]
and:
\[
\aligned
\log\,
\prod_{1\leqslant j\leqslant\ell}\,
\Big(
1
+
\frac{j}{n}
\Big)^{-1}
&
\,=\,
-\,
\log\,
\Big(
1
+
\frac{1}{n}
\Big)
-
\log\,
\Big(
1
+
\frac{2}{n}
\Big)
-\cdots-
\log\,
\Big(
1
+
\frac{\ell}{n}
\Big)
\\
&
\,\leqslant\,
0.
\endaligned
\]

Certainly, we have already seen implicitly in the proof
of the previous
Lemma~{\ref{Lemma-M-less-1}}
that all the logarithms of these products
are $\leqslant 0$. But we are now searching for 
a {\em minoration} of these coefficients:
\[
M_{k_2,\dots,k_n}^n
\,\,\geqslant\,\,
{\sf what}\text{\bf ?}
\]
For a reason that will become transparent 
just after a preliminary lemma, we will soon
restrict ourselves to suppose that:
\[
k_2+k_3+\cdots+k_{n-1}+k_n
\,\,\leqslant\,\,
\sqrt{n}.
\]

\begin{Lemma}
For all $0 \leqslant \delta \leqslant 3/5$:
\[
\log\,\big(1-\delta\big)
\,\,\geqslant\,\,
-\,\delta
-
\delta^2,
\]
and for all $\varepsilon \geqslant 0$:
\[
-\,\log\,\big(1+\varepsilon\big)
\,\,\geqslant\,\,
-\,\varepsilon.
\]
\end{Lemma}

\proof
The first inequality\,\,---\,\,which 
is in fact true for $0 \leqslant \delta
\leqslant 0, 683$ as can be seen with the
help of a computer\,\,---:
\[
-\,\delta
-
\frac{\delta^2}{2}
-
\frac{\delta^3}{3}
-
\frac{\delta^4}{4}
-
\frac{\delta^5}{5}
-\cdots
\overset{\text{\bf ?}}{\,\,\geqslant\,\,}
-\,\delta
-
\delta^2,
\]
is equivalent to:
\[
\frac{\delta^2}{2}
\overset{\text{\bf ?}}{\,\,\geqslant\,\,}
\frac{\delta^3}{3}
+
\frac{\delta^4}{4}
+
\frac{\delta^5}{5}
+\cdots.
\]
In this inequality under questioning, let us 
insert a computable infinite sum:
\[
1
\overset{\text{\bf ?}}{\,\,\geqslant\,\,}
\frac{2}{3}\,\delta\,
\big(
1
+
\delta
+
\delta^2
+
\cdots
\big)
\,\,\geqslant\,\,
\frac{2}{3}\,\delta
+
\frac{2}{4}\,\delta^2
+
\frac{2}{5}\,\delta^3
+\cdots,
\]
in order to come to an elementary minoration:
\[
1
\overset{\text{\bf ?}}{\,\,\geqslant\,\,}
\frac{2}{3}\,\delta\,
\frac{1}{1-\delta}
\ \ \ \ \
\Longleftrightarrow
\ \ \ \ \
3
-
3\,\delta
\overset{\text{\sf yes}}{\,\,\geqslant\,\,}
2\,\delta.
\]

The second inequality $\log\, (1+\varepsilon) \leqslant \varepsilon$
is well known.
\endproof

Now, let us suppose that:
\[
k
\,\leqslant\,
\sqrt{n},
\]
whence as soon as $n \geqslant 4$:
\[
\frac{k-1}{n}
\,<\,
\frac{1}{\sqrt{n}}
\,\leqslant\,
\frac{1}{2}
\,<\,
\frac{3}{5}.
\]
Then:
\[
\aligned
\sum_{i=0}^{k-1}\,
\log\,
\Big(
1
-
\frac{i}{n}
\Big)
&
\,\,\geqslant\,\,
-\,
\sum_{i=0}^{k-1}\,
\frac{i}{n}
-
\sum_{i=0}^{k-1}\,
\frac{i^2}{n^2}
\\
&
\,\,=\,\,
-\,\frac{(k-1)\,k}{2\,n}
-
\frac{(k-1)\,k\,(2\,k-1)}{6\,n^2}
\\
&
\,\,=\,\,
-\,
\frac{k^2}{2\,n}
+
\underline{\frac{k}{2\,n}}
-
\frac{k^3}{3\,n^2}
+
\underline{\frac{k^2}{2\,n^2}}
-
\underline{\frac{k}{6\,n^2}}.
\endaligned
\]

The three terms here underlined have a positive contribution
and we can even neglect the second of them:
\[
\frac{k}{2\,n}
+
\underline{\frac{k^2}{2\,n^2}}
-
\frac{k}{6\,n^2}
\,\,\geqslant\,\,
\frac{k}{2\,n}\,
\Big(
1
-
\frac{1}{3\,n}
\Big)
\,\,>\,\,
0.
\]
Therefore, we obtain a useful minoration:
\[
\aligned
\sum_{i=0}^{k-1}\,
\log\,
\Big(
1
-
\frac{i}{n}
\Big)
\,\,\geqslant\,\,
-\,\frac{k^2}{2\,n}
-
\frac{k^3}{3\,n^2}
&
\,\,=\,\,
-\,
\frac{k^2}{2\,n}
-
\frac{k^2}{3\,n}\,
\frac{k}{n}
\\
\explicationmath{{\textstyle{\frac{k}{n}}}\leqslant 1}
\ \ \ \ \ \ \ \ \ \ \ \ \ \ \ \ \ \ \ \ \ \ \ \ \ \
\ \ \ \ \ \ \ \ \ \ \ \ \ \ \ \ \ \ \ \ \ \ \ \ \ \
\ \ \ \ \ \ \ \ \ \ \ \ \ \ \ \ \ \ \ \ \ \ \ \ \ \
&
\,\,\geqslant\,\,
-\,\frac{k^2}{2\,n}
-
\frac{k^2}{3\,n}
\\
&
\,\,\geqslant\,\,
-\,\frac{k^2}{n}.
\endaligned
\]

Next, for the quotients of second type which are present in
the various $M_{k_2, \dots, k_n}^n$, the minoration work is easier:
\[
\aligned
-\,
\sum_{j=1}^\ell\,
\log\,
\Big(
1
+
\frac{j}{n}
\Big)
&
\,\,\geqslant\,\,
-\,\sum_{j=1}^\ell\,
\frac{j}{n}
\\
&
\,\,=\,\,
-\,\frac{\ell^2}{2\,n}
-
\frac{\ell}{2\,n}
\\
\explicationmath{\ell\leqslant\ell^2}
\ \ \ \ \ \ \ \ \ \ \ \ \ \ \ \ \ \ \ \ \ \ \ \ \ \
\ \ \ \ \ \ \ \ \ \ \ \ \ \ \ \ \ \ \ \ \ \ \ \ \ \
\ \ \ \ \ \ \ \ \ \ \ \ \ \ \ \ \ \ \ \ \ \ \ \ \ \
&
\,\,\geqslant\,\,
-\,\frac{\ell^2}{n}.
\endaligned
\]

Without forgetting the powers $n^k$ and $n^{-\ell}$,
these estimates can now be summarized as the following

\begin{Lemma}
For all $0 \leqslant k \leqslant \frac{3}{5}\, n$: 
\[
\frac{n!}{(n-k)!}
\,\,\geqslant\,\,
n^k\,
e^{-\,\frac{k^2}{n}}
\]
and for all $0 \leqslant
\ell$:
\[
\frac{n!}{(n+\ell)!}
\,\,\geqslant\,\,
n^{-\ell}\,
e^{-\,\frac{\ell^2}{n}}.
\eqno\qed
\]
\end{Lemma}

Importantly, we point out that there is a {\em uniform}
minoration:
\[
\frac{n!}{(n+m)!}
\,\,\geqslant\,\,
n^{-m}\,
e^{-\frac{m^2}{n}},
\]
valid for all integers $m \in \Z$, positive or negative, 
in the range:
\[
-\,{\textstyle{\frac{3}{5}}}\,
n
\,\,\leqslant\,\,
m
\,\,<\,\,
\infty.
\]
Notice that the exponential factor is always $\leqslant 1$. 

Next, thanks to all this, we will assume from now on that
the range of the integers $k_2, \dots, k_n$ is restricted
to:
\[
0
\,\leqslant\,
k_2+k_3+\cdots+k_{n-1}+k_n
\,\,\leqslant\,\,
\frac{\sqrt{n}}{c(n)},
\]
for some function $c(n) \underset{n\to\infty}{\longrightarrow} \infty$
that will be chosen later\,\,---\,\,think for instance
$c(n) := \log\, \log\, \log\, n$. In particular, this implies
that:
\[
0
\,\leqslant\,
k_2
\,\leqslant\,
{\textstyle{\frac{3}{5}}}\,
n,
\ \ \ \ \ \ \ \ \ \ \ \ \ \ \ \ \ \
\big\vert
k_2
-
k_3
\big\vert
\,\leqslant\,
{\textstyle{\frac{\sqrt{n}}{c(n)}}}
\,\leqslant\,
{\textstyle{\frac{3}{5}}},\,\,
\dots\dots,\,\,
\big\vert
k_{n-1}
-
k_n
\big\vert
\,\leqslant\,
{\textstyle{\frac{\sqrt{n}}{c(n)}}}
\,\leqslant\,
{\textstyle{\frac{3}{5}}},
\]
so that the lemma applies to minorize:
\[
\aligned
M_{k_2,k_3,\dots,k_{n-1},k_n}^n
&
\,=\,
\frac{n!}{(n-k_2)!}\,\,
\frac{n!}{(n+k_2-k_3)!}\,\,
\cdots\cdots\,\,
\frac{n!}{(n+k_{n-1}-k_n)!}\,\,
\frac{n!}{(n+k_n)!}
\\
&
\,\geqslant\,
\zero{n^{k_2}}\,
e^{-\frac{k_2^2}{n}}\,
\zero{n^{-k_2+k_3}}\,
e^{-\frac{(k_2-k_3)^2}{n}}\,\,
\cdots\cdots\,\,
\zero{n^{-k_{n-1}+k_n}}\,
e^{-\frac{(k_{n-1}-k_n)^2}{n}}\,
\zero{n^{-k_n}}\,
e^{-\frac{k_n^2}{n}}
\\
&
\,=\,
e^{-\frac{1}{n}\,
\left[
k_2^2+(k_2-k_3)^2+\cdots+(k_{n-1}-k_n)^2+k_n^2
\right]}
\\
&
\,=\,
e^{-\frac{1}{n}\,
\left[
2\,k_2^2-2\,k_2k_3+k_3^2
+\cdots+
2\,k_{n-1}^2-2\,k_{n-1}k_n+2\,k_n^2
\right]}
\\
&
\,\geqslant\,
e^{-\frac{1}{n}\,
\left[
2\,
\left(
k_2+k_3+\cdots+k_{n-1}+k_n
\right)^2
\right]}
\\
&
\,\geqslant\,
e^{-\frac{1}{n}\,2\,\frac{n}{c(n)^2}}
\\
&
\,=\,
e^{-\frac{2}{c(n)^2}}
\underset{n\to\infty}{\,\,\,\longrightarrow\,\,\,}
1.
\endaligned
\]
We thus have proved the key

\begin{Proposition}
\label{Proposition-minoration-M}
For any choice of function $c(n) \overset{n\to\infty}{\longrightarrow}
\infty$, the quantities:
\[
M_{k_2,k_3,\dots,k_{n-1},k_n}^n
\,=\,
\frac{n!}{(n-k_2)!}\,\,
\frac{n!}{(n+k_2-k_3)!}\,\,
\cdots\cdots\,\,
\frac{n!}{(n+k_{n-1}-k_n)!}\,\,
\frac{n!}{(n+k_n)!}
\]
enjoy the inequalities:
\[
e^{-\frac{2}{c(n)^2}}
\,\,\leqslant\,\,
M_{k_2,k_3,\dots,k_{n-1},k_n}^n
\,\,\leqslant\,\,
1,
\]
when their indices range in the set:
\reqnomode\usetagform{EngelLie}
\begin{align}
\Big\{
&
\big(k_2,k_3,\dots,k_{n-1},k_n\big)
\in
\N^n
\colon\,\,
0
\,\leqslant\,
k_2
\,\leqslant\,
n,
\notag
\\
&
\ \ \ \ \ \ \ \ \ \ \ \ \ \ \ \ \ \ \ \ \ \ \ \ \ \ \ \ \ \ \ \ \ \ \
\ \ \ \ \ \ \ \ \ \ \ \ \ 
0
\,\leqslant\,
k_3
\,\leqslant\,
n+k_2,
\notag
\\
&
\ \ \ \ \ \ \ \ \ \ \ \ \ \ \ \ \ \ \ \ \ \ \ \ \ \ \ \ \ \ \ \ \ \ \
\ \ \ \ \ \ \ \ \ \ \ \ \ 
\cdots\cdots\cdots\cdots\cdots\cdots
\notag
\\
&
\ \ \ \ \ \ \ \ \ \ \ \ \ \ \ \ \ \ \ \ \ \ \ \ \ \ \ \ \ \ \ \ \ \ \
\ \ \ \ \ \ \ \ \ \ \ \ \ 
0
\,\leqslant\,
k_{n-1}
\,\leqslant\,
n+k_{n-2},
\notag
\\
&
\ \ \ \ \ \ \ \ \ \ \ \ \ \ \ \ \ \ \ \ \ \ \ \ \ \ \ \ \ \ \ \ \ \ \
\ \ \ \ \ \ \ \ \ \ \ \ \ 
0
\,\leqslant\,
k_n
\ \ \
\,\leqslant\,
n+k_{n-1},
\notag
\\
&
k_2+k_3+\cdots+k_{n-1}+k_n
\,\leqslant\,
\frac{\sqrt{n}}{c(n)}
\Big\}.
\tag{\qed}
\end{align}
\end{Proposition}

\Section{\bf Majorant power series $\widehat{C}(w_2,\dots,w_n)$
and its diagonalization $\widehat{C}(x,\dots,x)$}
\label{majorant-widehat-C}
\HEAD{{\ref{majorant-widehat-C}}.~Majorant power series 
$\widehat{C}(w_2,\dots,w_n)$
and its diagonalization $\widehat{C}(x,\dots,x)$
}{
Jo\"el {\sc Merker} and The-Anh Ta,
D\'epartement de Math\'ematiques d'Orsay, 
Universit\'e Paris-Sud, France}

Now, come back to:
\[
F(x,y)
\,=\,
\frac{1-y}{1-2\,y+x\,y}
\]
and observe that for all $3 \leqslant i \leqslant n-1$:
\[
\frac{1-x^{i-1}}{1-2\,x^{i-1}+x^i}
\,\,=\,\,
F\big(x,\,x^{i-1}\big).
\]
Its expansion:
\[
\aligned
F(x,y)
&
\,=\,
1
+
\frac{y-x\,y}{1-2\,y+x\,y}
\\
&
\,=\,
1
+
\sum_{\ell=1}^\infty\,
y^\ell\,
\sum_{k=0}^\infty\,
x^k\,
F_{k,\ell}
\endaligned
\]
will be (easily) computed soon. With $F(x,y)$, introduce 
also\,\,---\,\,notice the single sign change in the denominator:
\[
\aligned
\widehat{F}(x,y)
\,=\,
\frac{1-y}{1-2\,y-x\,y}
&
\,=\,
1
+
\frac{y+x\,y}{1-2\,y-x\,y}
\\
&
\,=\,
1
+
\sum_{\ell=1}^\infty\,
y^\ell\,
\sum_{k=0}^\infty\,
x^k\,\widehat{F}_{k,\ell},
\endaligned
\]
a new function which
will act as a {\sl majorant series}, in the sense that:
\[
\big\vert
F_{k,\ell}
\big\vert
\,\,\leqslant\,\,
\widehat{F}_{k,\ell}
\eqno
{\scriptstyle{(\forall\,k\,\geqslant\,0,\,\,
\forall\,\ell\,\geqslant\,0)}}.
\]
Such inequalities are made transparent from
the following clear explicit
expressions, in which just a factor $(-1)^k$ drops.

\begin{Lemma}
\label{Lemma-expansion-F-F-hat}
With the convention that $\binom{\ell-1}{-1} = 0 = 
\binom{\ell-1}{\ell}$, the power series expansions are:
\[
\aligned
F(x,y)
&
\,=\,
1
+
\sum_{\ell=1}^\infty\,
y^\ell\,
\sum_{0\leqslant k\leqslant\ell}\,
(-1)^k\,
x^k\,
\Big[
2^{\ell-1-k}\,
{\textstyle{\binom{\ell-1}{k}}}
+
2^{\ell-k}\,
{\textstyle{\binom{\ell-1}{k-1}}}
\Big],
\\
\widehat{F}(x,y)
&
\,=\,
1
+
\sum_{\ell=1}^\infty\,
y^\ell\,
\sum_{0\leqslant k\leqslant\ell}\,
\ \ \ \ \ \ \ \ \ \
x^k\,
\Big[
2^{\ell-1-k}\,
{\textstyle{\binom{\ell-1}{k}}}
+
2^{\ell-k}\,
{\textstyle{\binom{\ell-1}{k-1}}}
\Big].
\endaligned
\]
\end{Lemma}

\proof
Expand:
\[
\aligned
F(x,y)
&
\,=\,
1
+
\frac{y-x\,y}{1-y\,(2-x)}
\\
&
\,=\,
1
+
\big(
y
-
x\,y
\big)\,
\sum_{h=0}^\infty\,
y^h\,
\big(
2
-
x
\big)^h
\\
&
\,=\,
1
+
\big(
y
-
x\,y
\big)\,
\sum_{h=0}^\infty\,
y^h\,
\sum_{0\leqslant m\leqslant h}\,
(-1)^m\,x^m\,2^{h-m}\,
\binom{h}{m}.
\endaligned
\]
Two double sums must be reorganized. In the first one,
replace $h = \ell-1$ and $m = k$:
\[
\sum_{h=0}^\infty\,
y^{h+1}\,
\sum_{0\leqslant k\leqslant h}\,
(-1)^m\,x^m\,
2^{h-m}\,
\binom{h}{m}
\,\,=\,\,
\sum_{\ell=1}^\infty\,
y^\ell\,
\sum_{0\leqslant k\leqslant\ell-1}\,
(-1)^k\,x^k\,
2^{\ell-1-k}\,
\binom{\ell-1}{k},
\]
and observe that the last sum can be extended to the range
$0 \leqslant k \leqslant \ell$, thanks to the convention.
In the second one, replace  $h = \ell-1$ and $m = k -1$:
\[
-\,
\sum_{h=0}^\infty\,
y^{h+1}\,
\sum_{0\leqslant m\leqslant h}\,
(-1)^m\,
x^{m+1}\,
2^{h-m}\,
\binom{h}{m}
\,\,=\,\,
-\,
\sum_{\ell=1}^\infty\,
\sum_{1\leqslant k\leqslant\ell}\,
(-1)^{k-1}\,x^k\,
2^{\ell-k}\,
\binom{\ell-1}{k-1},
\]
and observe that the term $k = 0$ in the sum can be included,
thanks to the convention. Adding these two expressions yield 
the stated power expansion of $F(x,y)$.

Next, for what concerns:
\[
\aligned
\widehat{F}(x,y)
&
\,=\,
1
+
\frac{y+x\,y}{1-2\,y-x\,y}
\\
&
\,=\,
1
+
\big(y+x\,y\big)\,
\sum_{h=0}^\infty\,
y^h\,
\big(2+x\big)^h
\\
&
\,=\,
1
+
\big(y+x\,y\big)\,
\sum_{h=0}^\infty\,
y^h\,
\sum_{0\leqslant m\leqslant h}\,
x^m\,2^{h-m}\,
\binom{h}{m},
\endaligned
\]
exactly the same transformations work, except that the
$(-1)^m$ factor has disappeared.
\endproof

Next, our goal is to introduce a majorant power series
$\widehat{C} (w_2, \dots, w_n)$ for the power series
$C (w_2, \dots, w_n)$. As anticipated above, it is now clear by
means of the triangle inequality that:
\leqnomode\usetagform{default}
\begin{align}
\label{F-smaller-F-hat}
\big\vert
F_{k,\ell}
\big\vert
\,\,\leqslant\,\,
\widehat{F}_{k,\ell},
\end{align}
for all $k \geqslant 0$ and all $\ell \geqslant 0$.
In terms of $F(x,y)$ and of the already seen power series:
\[
E(x)
\,:=\,
\frac{1-x}{1-2\,x}
\,=\,
\sum_{k=0}^\infty\,
E_k\,x^k,
\]
having positive coefficients $E_0 = 1$ and $E_k = 2^{k-1}$
for $k \geqslant 1$, we can write:
\[
\aligned
C\big(w_2,\dots,w_n\big)
&
\,=\,
\prod_{2\leqslant i\leqslant n}\,
\frac{1-w_2\cdots w_i}{1-2\,w_2\cdots w_i}
\\
&
\ \ \ \ \
\prod_{2\leqslant i<j\leqslant n}\,
\frac{1-w_{i+1}\cdots w_j}{
1-2\,w_{i+1}\cdots w_j+w_iw_{i+1}\cdots w_j}
\\
&
\,=:\,
\prod_{2\leqslant i\leqslant n}\,
E\big(w_2\cdots w_i\big)
\\
&
\ \ \ \ \
\prod_{2\leqslant i<j\leqslant n}\,
F\big(w_i,\,\,w_{i+1}\cdots w_j\big).
\endaligned
\]
Hence we may introduce similarly:
\[
\aligned
\widehat{C}\big(w_2,\dots,w_n\big)
&
\,:=\,
\prod_{2\leqslant i\leqslant n}\,
E\big(w_2\cdots w_i\big)
\\
&
\ \ \ \ \
\prod_{2\leqslant i<j\leqslant n}\,
\widehat{F}
\big(w_i,\,\,w_{i+1}\cdots w_j\big).
\endaligned
\]

The expansions of the factors of the first product show as:
\[
\aligned
C\big(w_2,\dots,w_n\big)
&
\,=\,
\prod_{2\leqslant i\leqslant n}\,
\bigg(
\sum_{k=0}^\infty\,
E_k\,
\big(
w_2\cdots w_i
\big)^k
\bigg)
\\
&
\ \ \ \ \
\prod_{2\leqslant i<j\leqslant n}\,
\bigg(
\sum_{k=0}^\infty\,
\sum_{\ell=0}^\infty\,
F_{k,\ell}\,
\big(w_i\big)^k\,
\big(
w_{i+1}\cdots w_j
\big)^\ell
\bigg)
\\
&
\,=:\,
\sum_{k_2,\dots,k_n\,\geqslant\,0}\,
C_{k_2,\dots,k_n}\,
\big(w_2\big)^{k_2}
\cdots
\big(w_n\big)^{k_n},
\endaligned
\]
and similarly:
\[
\aligned
\widehat{C}\big(w_2,\dots,w_n\big)
&
\,=\,
\prod_{2\leqslant i\leqslant n}\,
\bigg(
\sum_{k=0}^\infty\,
E_k\,
\big(
w_2\cdots w_i
\big)^k
\bigg)
\\
&
\ \ \ \ \
\prod_{2\leqslant i<j\leqslant n}\,
\bigg(
\sum_{k=0}^\infty\,
\sum_{\ell=0}^\infty\,
\widehat{F}_{k,\ell}\,
\big(w_i\big)^k\,
\big(
w_{i+1}\cdots w_j
\big)^\ell
\bigg)
\\
&
\,=:\,
\sum_{k_2,\dots,k_n\,\geqslant\,0}\,
\widehat{C}_{k_2,\dots,k_n}\,
\big(w_2\big)^{k_2}
\cdots
\big(w_n\big)^{k_n}.
\endaligned
\]

Since all $E_k \geqslant 0$ and all $\widehat{F}_{k,\ell} \geqslant
0$, we have all $\widehat{C}_{k_2,\dots,k_n} \geqslant 0$ as 
well\,\,---\,\,however, many $C_{k_2,\dots,k_n}$ are $\leqslant -1$.

Thanks to~({\ref{F-smaller-F-hat}}) and to the triangle inequality
in expansions, we obtain:
\leqnomode\usetagform{default}
\begin{align}
\label{C-k-C-hat-k}
\big\vert
C_{k_2,\dots,k_n}
\big\vert
\,\,\leqslant\,\,
\widehat{C}_{k_2,\dots,k_n},
\end{align}
for all $k_2, \dots, k_n \geqslant 0$, which means
that $\widehat{C}$ is a {\sl majorant power series} for $C$.
Notice that:
\[
C_{k_2,\dots,k_n}
\,\in\,
\Z
\ \ \ \ \ \ \ \ \ \ \ \ \ \ \ \ \ \
\text{and}
\ \ \ \ \ \ \ \ \ \ \ \ \ \ \ \ \ \
\widehat{C}_{k_2,\dots,k_n}
\,\in\,
\N.
\]

Now, passing to the diagonal:
\[
\big\{
w_2
\,=\,\cdots\,=\,
w_n
\,=:\,
x
\big\},
\]
we deduce for every $k \in \N$, again by means of the triangle
inequality:
\[
\aligned
\big\vert
C_h
\big\vert
&
\,=\,
\bigg\vert
\sum_{k_2+\cdots+k_n=h}\,
C_{k_2,\dots,k_n}
\bigg\vert
\\
&
\,\leqslant\,
\sum_{k_2+\cdots+k_n=h}\,
\big\vert
C_{k_2,\dots,k_n}
\big\vert
\\
\explicationtext{({\ref{C-k-C-hat-k}})}
\ \ \ \ \ \ \ \ \ \ \ \ \ \ \ \ \ \ \ \ \ \ \ \ \ \
&
\,\leqslant\,
\sum_{k_2+\cdots+k_n=h}\,
\widehat{C}_{k_2,\dots,k_n}
\,\,\,=:\,\,
\widehat{C}_h.
\endaligned
\]

In fact, these integers $\widehat{C}_h \geqslant 0$ express
as {\em coefficients} of the diagonal majorant series:
\[
\aligned
\widehat{C}^{n-1}(x)
\,:=\,
&\,
\widehat{C}(x,\dots,x)
\\
\,=\,
&\,
\prod_{2\leqslant i\leqslant n}\,
\frac{1-x^{i-1}}{1-2\,x^{i-1}}\,
\prod_{2\leqslant i<j\leqslant n}\,
\frac{1-x^{j-i}}{1-2\,x^{j-i}-x^{j-i+1}}
\\
&
\,=\,
\sum_{k_2,\dots,k_n\geqslant 0}\,
\widehat{C}_{k_2,\dots,k_n}\,
x^{k_2}\cdots x^{k_n}
\\
&
\,=\,
\sum_{h=0}^\infty\,
\bigg(
\sum_{k_2+\cdots+k_n=h}\,
\widehat{C}_{k_2,\dots,k_n}
\bigg)\,
x^h
\,\,\,=:\,\,
\sum_{h=0}^\infty\,
\widehat{C}_h\,x^h.
\endaligned
\]
Let us therefore state these observations as a

\begin{Lemma}
\label{Lemma-expressions-C-C-hat}
The $1$-variable products\big/series:
\[
\aligned
C^{n-1}(x)
&
\,:=\,
\prod_{i=1}^{n-1}\,
\frac{1-x^i}{1-2\,x^i}\,\,
\prod_{i=2}^{n-1}\,
\Big(
\frac{1-x^{i-1}}{1-2\,x^{i-1}+x^i}
\Big)^{n-i}
\,\,=\,\,
\sum_{h=0}^\infty\,
C_h^{n-1}\,
x^h,
\\
\widehat{C}^{n-1}(x)
&
\,:=\,
\prod_{i=1}^{n-1}\,
\frac{1-x^i}{1-2\,x^i}\,\,
\prod_{i=2}^{n-1}\,
\Big(
\frac{1-x^{i-1}}{1-2\,x^{i-1}-x^i}
\Big)^{n-1}
\,\,=\,\,
\sum_{h=0}^\infty\,
\widehat{C}_h^{n-i}\,
x^h,
\endaligned
\]
have coefficients satisfying the inequalities:
\[
\big\vert
C_h^{n-1}
\big\vert
\,\,\leqslant\,\,
\widehat{C}_h^{n-1}
\eqno
{\scriptstyle{(\forall\,h\,\geqslant\,0)}}.\,\qed
\]
\end{Lemma}

\Section{\bf Positivity of diagonal sums coefficients $C_h^{n-1}$}
\label{positivity-diagonal-sums-coefficients}
\HEAD{{\ref{positivity-diagonal-sums-coefficients}}.~Positivity 
of diagonal sums coefficients $C_h^{n-1}$
}{
Jo\"el {\sc Merker} and The-Anh Ta,
D\'epartement de Math\'ematiques d'Orsay, 
Universit\'e Paris-Sud, France}

Now, study the power series $C(w_2, \dots, w_n)$ along the diagonal:
\[
\big\{
w_2
\,=\,\cdots\,=\,
w_n
\,=:\,
x
\big\},
\]
that is to say, introduce:
\[
\aligned
C^{n-1}(x)
\,:=\,
&\,
C\big(x,\dots,x\big)
\\
\,=\,
&\,
\sum_{k_2=0}^\infty\,
\cdots
\sum_{k_n=0}^\infty\,
C_{k_2,\dots,k_n}\,
x^{k_2}\cdots x^{k_n}
\\
\,=\,
&\,
\sum_{h=0}^\infty\,
\bigg(
\sum_{k_2+\cdots+k_n=h}\,
C_{k_2,\dots,k_n}
\bigg)\,
x^h
\\
\,=:\,
&\,
\sum_{h=0}^\infty\,
C_h^{n-1}\,
x^h,
\endaligned
\]
in terms of certain integer coefficients $C_h^{n-1} \in \Z$. 
In fact, coming back to the product expression
of $C(w_2, \dots, w_n)$, we realize that in:
\[
\aligned
C^{n-1}(x)
\,=\,
\frac{1-x}{1-2\,x}\,\,
\frac{1-x^2}{1-2\,x^2}\,
\,\,\cdots\cdots\,\,
\frac{1-x^{n-1}}{1-2\,x^{n-1}}
\\
\frac{1-x}{\zero{1-2\,x}}\,\,
\,\,\cdots\cdots\,\,
\frac{1-x^{n-2}}{\zero{1-2\,x^{n-2}}}
\\
\ddots
\ \ \ \ \ \ \ \ \ \ \ \ 
\vdots
\ \ \ \ \ \ \ 
\\
\frac{1-x}{\zero{1-2\,x}}
\\
\frac{\zero{1-2\,x}}{1-2\,x+x^2}
\cdots
\frac{\zero{1-2\,x^{n-2}}}{1-2\,x^{n-2}+x^{n-1}}
\\
\ddots
\ \ \ \ \ \ \ \ \ \ \ \ \ \ \ \ \ 
\vdots
\ \ \ \ \ \ \ \ \ \ \ \ \ 
\\
\frac{\zero{1-2\,x}}{1-2\,x+x^2},
\endaligned
\]
some simplifications indicated by underlinings conduct us to:
\[
\aligned
C^{n-1}(x)
&
\,=\,
\frac{1-x}{1-2\,x}\,\,
\frac{1-x^2}{1-2\,x^2}\,\,
\frac{1-x^3}{1-2\,x^3}\,\,
\cdots\cdots\cdots\cdots\cdots\cdots\,\,
\frac{1-x^{n-1}}{1-2\,x^{n-1}}
\\
&
\ \ \ \ \ \ \ \ \ \
\bigg(
\frac{1-x}{1-2\,x+x^2}
\bigg)^{n-2}
\bigg(
\frac{1-x^2}{1-2\,x^2+x^3}
\bigg)^{n-3}
\!\!\!\!\!\cdots\cdots
\bigg(
\frac{1-x^{n-2}}{1-2\,x^{n-2}+x^{n-1}}
\bigg)^1.
\endaligned
\]

Furthermore, 
on the second line, the first fraction 
to the power $(\cdot)^{n-2}$ trivially 
simplifies as:
\[
\frac{1-x}{(1-x)^2}
\,=\,
\frac{1}{1-x},
\]
whence:
\[
\aligned
C^{n-1}(x)
&
\,=\,
\frac{1-x}{1-2\,x}\,\,
\frac{1-x^2}{1-2\,x^2}\,\,
\frac{1-x^3}{1-2\,x^3}\,\,
\cdots\cdots\cdots\cdots\cdots\cdots\,\,
\frac{1-x^{n-1}}{1-2\,x^{n-1}}
\\
&
\ \ \ \ \ \ \ \ \ \
\bigg(
\frac{1}{1-x}
\bigg)^{n-2}
\bigg(
\frac{1-x^2}{1-2\,x^2+x^3}
\bigg)^{n-3}
\!\!\!\!\!\cdots\cdots
\bigg(
\frac{1-x^{n-2}}{1-2\,x^{n-2}+x^{n-1}}
\bigg)^1.
\endaligned
\]

Let us focus on the second line, which we now call:
\[
\aligned
P^{n-1}(x)
\,:=\,
&\,
\bigg(
\frac{1}{1-x}
\bigg)^{n-2}
\bigg(
\frac{1-x^2}{1-2\,x^2+x^3}
\bigg)^{n-3}
\!\!\!\!\!\cdots\cdots
\bigg(
\frac{1-x^{n-2}}{1-2\,x^{n-2}+x^{n-1}}
\bigg)^1
\\
\,=:\,
&\,
\sum_{h=0}^\infty\,
P_h^{n-1}\,
x^h.
\endaligned
\]
We believe that all the coefficients of the full product
$C^{n-1}(x)$ are positive,
but a restricted statement will be enough for our purposes.

\begin{Lemma}
For all indices $h$ in the range:
\[
0
\,\leqslant\,
h
\,\leqslant\,
\lfloor
\sqrt{n}
\rfloor
\]
one has:
\[
\aligned
P_h^{n-1}
&
\,\geqslant\,
1,
\\
C_h^{n-1}
&
\,\geqslant\,
2^h.
\endaligned
\]
\end{Lemma}

\proof
First, we make the following transformation for each term in the
product $P^{n-1}(x)$:
\reqnomode\usetagform{EngelLie}
\begin{align}
\bigg(
\frac{1-x^k}{1-2\,x^k+x^{k+1}}
\bigg)^{n-k-1}
&=
\bigg(
\frac{1-x^k}{1-x^k - (x^k-x^{k+1})}
\bigg)^{n-k-1} 
\notag
\\
&=
\bigg(
\frac{1}{1 - \frac{x^k-x^{k+1}}{1-x^k}}
\bigg)^{n-k-1}
\notag
\\
&=
\bigg(
\frac{1}{1 - \frac{x^k}{1+x+ \cdots +x^{k-1}}}
\bigg)^{n-k-1}
\tag{(1\,\leqslant\,k\,\leqslant\,n-2).}
\end{align}

Using the expansion and factorization:
\[
\frac{1}{1-T} 
=
\Big(
1+T
\Big)
\,\sum\limits_{i=0}^{\infty}
T^{2i}
= 
\Big(
1+T
\Big)
\Big(
1+T^2+T^4+T^6+\cdots
\Big) 
\]
and substituting $T = \frac{x^k}{1+x+ \cdots +x^{k-1}}$ gives us:
\[
\aligned
\bigg(
\frac{1}{1 - \frac{x^k}{1+x+ \cdots +x^{k-1}}}
\bigg)^1
&=
\bigg(
1+\frac{x^k}{1+x+ \cdots +x^{k-1}}
\bigg)
\bigg(
\sum\limits_{i=0}^{\infty}
\Big(
\frac{x^k}{1+x+ \cdots +x^{k-1}}
\Big)^{2i}
\bigg) \\
&=
\bigg(
\frac{1+x+ \cdots +x^k}{1+x+ \cdots +x^{k-1}}
\bigg)
\bigg(
\sum\limits_{i=0}^{\infty}
\Big(
\frac{x^k}{1+x+ \cdots +x^{k-1}}
\Big)^{2i}
\bigg).
\endaligned
\]

We then put together these expansions of terms in the product $P^{n-1}(x)$ to obtain: 
\[
\aligned
P^{n-1}(x)
=
&\Big(
1+x
\Big)^{n-2}
\Big(
1+x^2+x^4+x^6+\cdots
\Big)^{n-2} \\
&.\bigg(
\frac{1+x+x^2}{1+x}
\bigg)^{n-3}
\bigg(
1
+\Big(\frac{x^2}{1+x} \Big)^2
+\Big(\frac{x^2}{1+x} \Big)^4
+\cdots
\bigg)^{n-3} \\
& \!\!\;\;\;\;\;\;\;\;\;\;\;\;\;\;\;\;\;\;\;\cdots\cdots\cdots \\
&.\bigg(
\frac{1+x+\cdots+x^{k-1}+x^k}{1+x+\cdots+x^{k-1}}
\bigg)^{n-k-1}
\bigg(
\sum\limits_{i=0}^{\infty}
\Big(
\frac{x^k}{1+x+\cdots+x^{k-1}} 
\Big)^{2i}
\bigg)^{n-k-1} \\
& \!\!\;\;\;\;\;\;\;\;\;\;\;\;\;\;\;\;\;\;\;\cdots\cdots\cdots \\
&.\bigg(
\frac{1+x+\cdots+x^{n-3}+x^{n-2}}{1+x+\cdots+x^{n-3}}
\bigg)^{1}
\bigg(
1
+\Big(
\frac{x^{n-2}}{1+x+\cdots+x^{n-3}} 
\Big)^2
+\cdots
\bigg)^{1}.
\endaligned
\]

Notice that the product of the first terms in all lines admits simplification as follows:
\[
\aligned
\Big(
1+x
\Big)^{n-2}
&\bigg(
\frac{1+x+x^2}{1+x}
\bigg)^{n-3}
\cdots
\bigg(
\frac{1+\cdots+x^{k-1}+x^k}{1+\cdots+x^{k-1}}
\bigg)^{n-k-1}
\cdots
\bigg(
\frac{1+\cdots+x^{n-3}+x^{n-2}}{1+\cdots+x^{n-3}}
\bigg)^{1} \\
&=
(1+x)(1+x+x^2) \cdots (1+x+\cdots+x^k) \cdots (1+x+\cdots+x^{n-2}),
\endaligned
\]
while the other terms can be expanded using
\[
\aligned
\Big(
\sum\limits_{i=0}^{\infty} 
T^{2i}
\Big)^m
&=
\Big(
1+T^2+T^4+\cdots +T^{2j}+\cdots
\Big)^m \\
&=
\sum\limits_{j=0}^{\infty} 
\binom{m+j-1}{j} 
T^{2j} \\
&= 
1
+
\binom{m}{1} T^2 
+
\binom{m+1}{2} T^4 
+
\cdots 
+
\binom{m+j-1}{j}T^{2j} 
+ 
\cdots
.
\endaligned
\]

The expansion of $P^{n-1}(x)$ now becomes
\[
\aligned
P^{n-1}(x)
=
&\Big(
1+x
\Big)
\Big(
1+x+x^2
\Big) 
\cdots 
\Big(
1+x+\cdots+x^k
\Big) 
\cdots 
\Big(
1+x+\cdots+x^{n-2}
\Big) \\
&.\bigg(
1
+\binom{n-2}{1} x^2
+\binom{n-1}{2}x^4
+\binom{n}{3}x^6
+\cdots
\bigg) \\
&.\bigg(
1
+\binom{n-3}{1} \Big(\frac{x^2}{1+x} \Big)^2
+\binom{n-2}{2} \Big(\frac{x^2}{1+x} \Big)^4
+\cdots
\bigg) \\
& \!\!\;\;\;\;\;\;\;\;\;\;\;\;\;\;\;\;\;\;\;\cdots\cdots\cdots \\
&.\bigg(
\sum\limits_{i=0}^{\infty}
\binom{n-k-1 + i-1}{i}
\Big(
\frac{x^k}{1+x+\cdots+x^{k-1}} 
\Big)^{2i}
\bigg) \\
& \!\!\;\;\;\;\;\;\;\;\;\;\;\;\;\;\;\;\;\;\;\cdots\cdots\cdots \\
&.\bigg(
1
+\Big(\frac{x^{n-2}}{1+x+\cdots+x^{n-3}} \Big)^2
+\Big(\frac{x^{n-2}}{1+x+\cdots+x^{n-3}} \Big)^4
+\cdots
\bigg).
\endaligned
\]

Since we are only interested in the coefficients $P^{n-1}_h$ with $ 0 \,\leqslant\, h \,\leqslant\, \lfloor \sqrt{n} \rfloor $,  we will ignore the terms $\big(\frac{x^k}{1+x+\cdots+x^{k-1}} \big)^{2i}$ with $k \cdot 2i > n$, i.e. with $i > \frac{n}{2k}$. The first $\lfloor \sqrt{n} \rfloor$ coefficients in the power series expansion of $P^{n-1}(x)$ are the same as those of
\[
\aligned
\Big(
&1+x
\Big)
\Big(
1+x+x^2
\Big)
\cdots 
\Big(
1+x+\cdots+x^k
\Big)
\cdots 
\Big(
1+x+\cdots+x^{n-2}
\Big) \\
&.\bigg(
1
+\binom{n-2}{1} x^2
+\binom{n-1}{2}x^4
+\cdots
+\binom{n - 2 + \lfloor \frac{\sqrt{n}}{2} \rfloor -1}{\lfloor \frac{\sqrt{n}}{2} \rfloor}
x^{2 \lfloor \frac{\sqrt{n}}{2} \rfloor}
\bigg) \\
&.\bigg(
1
+\binom{n-3}{1}
\Big(
\frac{x^2}{1+x} 
\Big)^2
+\cdots
+
\binom{n- 3+ \lfloor \frac{\sqrt{n}}{4} \rfloor -1}{\lfloor \frac{\sqrt{n}}{4} \rfloor}
\Big(
\frac{x^2}{1+x} 
\Big)^{2\lfloor \frac{\sqrt{n}}{4} \rfloor}
\bigg) \\
& \!\!\;\;\;\;\;\;\;\;\;\;\;\;\;\;\;\;\;\;\;\cdots\cdots\cdots \\
&.\bigg(
1
+\cdots
+\binom{n-k-1 + \lfloor \frac{\sqrt{n}}{2k} \rfloor -1}{\lfloor \frac{\sqrt{n}}{2k} \rfloor}
\Big(
\frac{x^k}{1+x+\cdots+x^{k-1}} 
\Big)^{2 \lfloor \frac{\sqrt{n}}{2k} \rfloor}
\bigg) \\
& \!\!\;\;\;\;\;\;\;\;\;\;\;\;\;\;\;\;\;\;\;\cdots\cdots\cdots \\
&.\bigg(
1
+
\binom{n-\lfloor \frac{\sqrt{n}}{2} \rfloor -1}{1}
\Big(
\frac{x^{
\lfloor \frac{\sqrt{n}}{2} \rfloor}}
{1+x+\cdots+x^{\lfloor \frac{\sqrt{n}}{2} \rfloor -1}} 
\Big)^2
\bigg).
\endaligned
\]

Now it is clear that in order to show the positivity of $P^{n-1}_h$ for all $ 0 \,\leqslant\, h \,\leqslant\, \lfloor \sqrt{n} \rfloor $, it suffices to prove that the product
\[
\Big(
1+x
\Big)^{2 \lfloor \frac{\sqrt{n}}{4} \rfloor}
\Big(
1+x+x^2
\Big)^{2 \lfloor \frac{\sqrt{n}}{6} \rfloor}
\!\!\cdots
\Big(
1+x+\cdots+x^{\lfloor \frac{\sqrt{n}}{2} \rfloor -1}
\Big)^2
\]
is divisible by 
\[
\Big(
1+x
\Big)
\Big(
1+x+x^2
\Big)
\cdots 
\Big(
1+x+\cdots+x^k
\Big)
\cdots 
\Big(
1+x+\cdots+x^{n-2}
\Big),
\]
and at the same time that the quotient also has nonnegative coefficients.

Note that for any $j \geqslant 0$, one has 
\[
1+x+\cdots + x^{kj+k-1} 
= 
\Big(
1+x+\cdots+x^{k-1}
\Big)
\Big(
1+x^k+x^{2k}+\cdots+x^{kj}
\Big),
\]
that is 
$
1+x+\cdots+x^{k-1}
$ 
is divisible by 
$
1+x+\cdots + x^{kj+k-1}
$ 
with quotient having nonnegative coefficients. 

Now, we divide the set of indices 
$
\Big\{ 
1,2, \dots, (\lfloor \sqrt{n} \rfloor)^2-1
\Big\}
$ 
into 
$
\lfloor \sqrt{n} \rfloor
$ 
disjoint sets:
\[
\Big\{
\lfloor \sqrt{n} \rfloor j+1,
\lfloor \sqrt{n} \rfloor j+2,
\dots,
\lfloor \sqrt{n} \rfloor j+\lfloor \sqrt{n} \rfloor
\Big\}
\]
for 
$
j=0,1, \dots, \lfloor \sqrt{n} \rfloor -1
$. 
Then, for each index $k$, the number of integers of the form $kj+k-1$ in the interval 
$
\Big\{ 
\lfloor \sqrt{n} \rfloor k+1,
\lfloor \sqrt{n} \rfloor k+2,
\dots,
\lfloor \sqrt{n} \rfloor k+\lfloor \sqrt{n} \rfloor \Big\}
$ 
is at least 
$
\frac{\lfloor \sqrt{n} \rfloor}{k}
$.

Since 
$ 
2 \lfloor \frac{\sqrt{n}}{2k} \rfloor
\leqslant
\frac{\lfloor \sqrt{n} \rfloor}{k}
$, 
the polynomial 
$
\Big(
1+x+\cdots+x^{k-1}
\Big)^{2 \lfloor \frac{\sqrt{n}}{2k} \rfloor}
$ 
is divisible by 
\[
\prod\limits_{i= \lfloor \sqrt{n} \rfloor k+1}^{\lfloor\sqrt{n} \rfloor (k+1)}
\Big(
1+x+\cdots+x^i
\Big),
\]
with quotient having nonnegative coefficients. 
Taking in account all the values of 
$
k = 1,2, \dots, \lfloor \frac{\sqrt{n}}{2} \rfloor -1
$
, and making the product of all the 
$
\lfloor \frac{\sqrt{n}}{2} \rfloor -1
$ 
terms gives us the desired divisibility. 

At this point, notice further that the set
$
\Big\{ 
1,
2,
\dots,
\lfloor \sqrt{n} \rfloor \Big\},
$ 
corresponding to $k=0$, has not been used in obtaining the above divisibility. Thus, the first 
$
\lfloor \sqrt{n} \rfloor
$
coefficients of $P^{n-1}(x)$ are those of the product between 
\[
\Big(
1+x
\Big)
\Big(
1+x+x^2
\Big)
\cdots 
\Big(
1+x+\cdots+x^{\lfloor \sqrt{n} \rfloor}
\Big)
\]
and
a power series having constant coefficient 1 and the first 
$
\lfloor \sqrt{n} \rfloor
$
coefficients nonnegative.
This clearly implies the positivity of $P^{n-1}_h$ for all 
$ 
0 
\,\leqslant\, 
h 
\,\leqslant\, 
\lfloor \sqrt{n} \rfloor.
$

For the first 
$
\lfloor \sqrt{n} \rfloor
$ 
coefficients in the power series expansion of $C^{n-1}(x)$, it is enough to consider the product 
\[
\bigg( 
\frac{1-x}{1-2x} 
\bigg)
P^{n-1}(x),
\]
since all the remaining terms in the product $C^{n-1}(x)$ have power series expansions with nonnegative coefficients and constant coefficient $1.$ Now using the expansion
\[
 \frac{1-x}{1-2x} 
 =
 1+ \sum\limits_{i=1}^{\infty} 2^{i-1}x^i,
\]
we get  
\[
C^{n-1}_h 
= 
P^{n-1}_h
+
\sum\limits_{i=1}^{h} 2^{i-1} P^{n-1}_{h-i}.
\]

Since we have already showed that 
$
P^{n-1}_h 
\geqslant 1
$ 
for all 
$ 
0 
\,\leqslant\, 
h 
\,\leqslant\, 
\lfloor \sqrt{n} \rfloor 
$
, it follows that
\[
C^{n-1}_h 
\geqslant
1
+
\sum\limits_{i=1}^{h} 2^{i-1}
= 2^h
\]
for all $ 0 \,\leqslant\, h \,\leqslant\, \lfloor \sqrt{n} 
\rfloor$. This finishes our proof of the lemma.
\endproof

\Section{\bf Cauchy inequalities}
\label{Cauchy-inequalities}
\HEAD{{\ref{Cauchy-inequalities}}.~Cauchy inequalities
}{
Jo\"el {\sc Merker} and The-Anh Ta,
D\'epartement de Math\'ematiques d'Orsay, 
Universit\'e Paris-Sud, France}

Next, we will set up 
a useful (and trivial) version of the Cauchy inequalities 
for power series
having nonnegative coefficients. We start by determining
the radius of convergence $\RR > 0$ of $C^{n-1}(x)$ and
the one $\widehat{\RR} > 0$ of $\widehat{C}^{n-1}(x)$,
where, from Lemma~{\ref{Lemma-expressions-C-C-hat}}:
\[
\aligned
C^{n-1}(x)
&
\,:=\,
\prod_{i=1}^{n-1}\,
\frac{1-x^i}{1-2\,x^i}\,\,
\prod_{i=2}^{n-1}\,
\Big(
\frac{1-x^{i-1}}{1-2\,x^{i-1}+x^i}
\Big)^{n-i},
\\
\widehat{C}^{n-1}(x)
&
\,:=\,
\prod_{i=1}^{n-1}\,
\frac{1-x^i}{1-2\,x^i}\,\,
\prod_{i=2}^{n-1}\,
\Big(
\frac{1-x^{i-1}}{1-2\,x^{i-1}-x^i}
\Big)^{n-1},
\endaligned
\]

\begin{Lemma}
The smallest moduli of poles of the
two rational functions $C^{n-1}(x)$ and
$\widehat{C}^{n-1}(x)$ are:
\[
\RR
\,:=\,
{\textstyle{\frac{1}{2}}}
\,=\,
0.5
\ \ \ \ \ \ \ \ \ \ \ \ \ \ \ \ \ \
\text{and}
\ \ \ \ \ \ \ \ \ \ \ \ \ \ \ \ \ \
\widehat{\RR}
\,:=\,
\sqrt{2}-1
\,\approx\,
0.414\cdots.
\]
\end{Lemma}

\proof
The moduli of the roots of the denominator of the first
product $\prod_{1\leqslant i\leqslant n-1}\,
\frac{\ast}{1-2\,x^i}$ appearing in $C^{n-1}(x)$ are
$\frac{1}{2}$, $\frac{1}{\sqrt[2]{2}}$, 
$\frac{1}{\sqrt[3]{2}}$, \dots, $\frac{1}{\sqrt[n-1]{2}}$,
and the smallest among them is $\frac{1}{2}$.
But then in the disc $\big\{ x \in \C \colon\,
\vert x\vert < \frac{1}{2} \big\}$, we assert that
all denominators in the second product constituting
$C^{n-1}(x)$ are nowhere vanishing. Indeed, as already
observed above, taking account of
the simplification for $i = 2$:
\[
\frac{1-x^{2-1}}{1-2\,x^{2-1}+x^2}
\,=\,
\frac{1}{1-x},
\]
this second product writes as:
\[
\Big(
\frac{1}{1-x}
\Big)^{n-2}\,
\prod_{i=3}^{n-1}\,
\Big(
\frac{\ast}{1-2\,x^{i-1}+x^i}
\Big)^{n-i}.
\]
Then the root $1$ is certainly $> \frac{1}{2}$, while the subsequent
denominators for $3 \leqslant i \leqslant n$ are nonvanishing
when $\vert x\vert \leqslant \frac{1}{2}$, because:
\[
\aligned
\big\vert
1-2\,x^{i-1}+x^i
\big\vert
&
\,\geqslant\,
1
-
2\,\vert x\vert^{i-1}
-
\vert x\vert^i
\\
&
\,\geqslant\,
1
-
2\,
\big(
{\textstyle{\frac{1}{2}}}
\big)^{i-1}
-
\big(
{\textstyle{\frac{1}{2}}}
\big)^i
\\
&
\,\geqslant\,
1
-
2\,
\big(
{\textstyle{\frac{1}{2}}}
\big)^2
-
\big(
{\textstyle{\frac{1}{2}}}
\big)^3
\,\,\,=\,
{\textstyle{\frac{3}{8}}}
\,\,>\,
0.
\endaligned
\]

On the other hand, while the first product constituting 
$\widehat{C}^{n-1}(x)$ is exactly the same, such a simplification
in the second product does not occur, and in fact, in:
\[
\Big(
\frac{\ast}{1-2\,x-x^2}
\Big)^{n-2}\,\,
\prod_{i=3}^{n-1}\,
\Big(
\frac{\ast}{1-2\,x^{i-1}-x^i}
\Big)^{n-i},
\]
the same minoration for $3 \leqslant i\leqslant n-1$ applies:
\[
\big\vert
1
-
2\,x^{i-1}
-
x^i
\big\vert
\,\,\geqslant\,\,
1
-
2\,\vert x\vert^{i-1}
-
\vert x\vert^i
\,\,\geqslant\,\,
{\textstyle{\frac{3}{8}}},
\]
whereas the positive root $\sqrt{2}-1$ of $1-2\,x-x^2 = 0$
is smaller than $\frac{1}{2}$, and the other root $-1-\sqrt{2}$
has (much) larger modulus.
\endproof

Let therefore $0 < \rho < \sqrt{2}-1$ be any radius in these
convergence discs. A trivial version of the Cauchy inequalities
for power series having nonnegative coefficients is as follows.
From:
\[
\widehat{C}^{n-1}(\rho)
\,=\,
\sum_{h=0}^\infty\,
\widehat{C}_h^{n-1}\,
\rho^h,
\]
it comes for any $h \in \N$ fixed, since all terms are $\geqslant 0$:
\[
\widehat{C}^{n-1}(\rho)
\,\,\geqslant\,\,
\widehat{C}_h^{n-1}\,
\rho^h.
\]
Soon, we will take $\rho = \rho(n) \
\overset{>}{\underset{n\to\infty}{\longrightarrow}} 0$, in fact:
\[
\rho 
\,:=\,
\frac{1}{\sqrt{n}}
\eqno
{\scriptstyle{(\text{later})}}.
\]

\begin{Observation}
\label{Observation-Cauchy-inequalities}
For any $0 < \rho < \sqrt{2}-1$ and every $h \in \N$:
\[
\widehat{C}_h
\,\,\leqslant\,\,
\frac{1}{\rho^h}\,
\widehat{C}(\rho).
\eqno\qed
\]
\end{Observation}

Section~{\ref{vert-z-vert-rho}} provides an exploration of the way 
moduli of the elementary constituents $\frac{1-x^k}{1-2\,x^k}$ and
$\frac{1-x^\ell}{1-2\,x^\ell-x^{\ell+1}}$ vary
with drastic oscillations on circles
$\{ \vert x \vert = \rho\}$.

Thanks to these basic Cauchy inequalities, we can now start to 
control the growth of $\widehat{C}^{n-1}(\rho)$.

\Section{\bf Estimations of $\widehat{C}(\frac{1}{r})$ and of
$C(\frac{1}{r})$}
\label{estimations-C-1-r-C-hat-1-r}
\HEAD{{\ref{estimations-C-1-r-C-hat-1-r}}.~Estimations of 
$\widehat{C}(\frac{1}{r})$ and of
$C(\frac{1}{r})$
}{
Jo\"el {\sc Merker} and The-Anh Ta,
D\'epartement de Math\'ematiques d'Orsay, 
Universit\'e Paris-Sud, France}

At first, we reorganize $\widehat{C}^{n-1}(x)$
from Lemma~{\ref{Lemma-expressions-C-C-hat}}, writing its
second product up to $i = n$ included instead of $i = n-1$,
using $(\ast)^{n-n}=1$:
\[
\aligned
\widehat{C}^{n-1}(x)
&
\,=\,
\prod_{i=1}^{n-1}\,
\frac{1-x^i}{1-2\,x^i}\,\,
\prod_{i=2}^n\,
\Big(
\frac{1-x^{i-1}}{1-2\,x^{i-1}-x^i}
\Big)^{n-i}
\\
&
\,=\,
\underbrace{
\prod_{i=1}^{n-1}\,
\frac{1}{1-2\,x^i}}_{i\,=:\,k}
\,\,
\underbrace{
\prod_{i=1}^{n-1}\,
\big(
1-x^i
\big)}_{i\,=:\,k}
\,\,
\underbrace{
\prod_{i=2}^n\,
\big(
1-x^{i-1}
\big)^{n-i}}_{i\,=:\,k+1}
\,\,
\underbrace{
\prod_{i=2}^n\,
\frac{1}{(1-2\,x^{i-1}-x^i)^{n-i}}}_{
i\,=:\,\,k+1}
\\
&
\,=\,
\prod_{k=1}^{n-1}\,
\frac{1}{1-2\,x^k}\,\,
\underline{
\prod_{k=1}^{n-1}\,
\big(1-x^k\big)^1\,\,
\prod_{k=1}^{n-1}\,
\big(
1-x^k
\big)^{n-k-1}}\,\,
\prod_{k=1}^{n-1}\,
\frac{1}{(1-2\,x^k-x^{k+1})^{n-k-1}}
\\
&
\,=\,
\prod_{k=1}^{n-1}\,
\frac{1}{1-2\,x^k}\,\,
\underline{
\prod_{k=1}^{n-1}\,
\big(1-x^k\big)^{n-k}}\,\,
\prod_{k=1}^{n-1}\,
\frac{1}{(1-2\,x^k-x^{k+1})^{n-k-1}}.
\endaligned
\]

In order to set up a general statement, we will take:
\[
x
\,:=\,
\frac{1}{r},
\]
with $r = r(n) \underset{n\to\infty}{\longrightarrow}
\infty$, always with $0 < \frac{1}{r} < \sqrt{2}-1$.
In fact, to fix ideas, we shall assume at least $r \geqslant 10$.

\begin{Lemma}
\label{Lemma-majoration-C-hat}
One has:
\[
\widehat{C}
\big(
{\textstyle{\frac{1}{r}}}
\big)
\,\,\leqslant\,\,
e^{\,\frac{n}{r}+12\,\frac{n}{r^2}}.
\]
\end{Lemma}

\proof
Take logarithm:
\[
\aligned
\log\,
\widehat{C}
\big(
{\textstyle{\frac{1}{r}}}
\big)
&
\,=\,
\sum_{k=1}^{n-1}\,
\bigg(
-\,\log\,
\Big(
1
-
\frac{2}{r^k}
\Big)
+
(n-k)\,
\log\,
\Big(
1
-
\frac{1}{r^k}
\Big)
-
(n-k-1)\,
\log\,
\Big(
1
-
\frac{2}{r^k}
-
\frac{1}{r^{k+1}}
\Big)
\bigg)
\\
&
\,=\,
-\,\log\,
\Big(
1
-
\frac{2}{r}
\Big)
+
(n-1)\,\log\,
\Big(
1
-
\frac{1}{r}
\Big)
-
(n-2)\,
\log\,
\Big(
1
-
\frac{2}{r}
-
\frac{1}{r^2}
\Big)
\,-
\\
&
\ \ \ \ \
-\,\log\,
\Big(
1
-
\frac{2}{r^2}
\Big)
+
(n-2)\,\log\,
\Big(
1
-
\frac{1}{r^2}
\Big)
-
(n-3)\,
\log\,
\Big(
1
-
\frac{2}{r^2}
-
\frac{1}{r^3}
\Big)
\,+
\\
&
\ \ \ \ \
+
\sum_{k=3}^{n-1}\,
\bigg\{
-\,\log\,
\Big(
1
-
\frac{2}{r^k}
\Big)
+
(n-k)\,
\log\,
\Big(
1
-
\frac{1}{r^k}
\Big)
-
(n-k-1)\,
\log\,
\Big(
1
-
\frac{2}{r^k}
-
\frac{1}{r^{k+1}}
\Big)
\bigg\}.
\endaligned
\]
Now, employ the majorations valuable for $0 \leqslant \delta \leqslant
0.5$:
\[
\aligned
\log\,\big(1-\delta\big)
&
\,\,\leqslant\,\,
-\,\delta
-
{\textstyle{\frac{1}{2}}}\,
\delta^2,
\\
-\,\log\,\big(1-\varepsilon\big)
&
\,\,\leqslant\,\,
\varepsilon
+
\varepsilon^2,
\endaligned
\]
to get using the assumption $r \geqslant 10$:
\begin{align}
\!\!\!\!\!\!\!\!\!\!\!\!\!\!\!
\log\,
\widehat{C}
\big(
{\textstyle{\frac{1}{r}}}
\big)
&
\,\leqslant\,
\frac{2}{r}
+
\frac{4}{r^2}
+
(n-1)\,
\bigg[
-\frac{1}{r}
-
\frac{1}{2}\,\frac{1}{r^2}
\bigg]
+
(n-2)\,
\bigg[
\Big(
\frac{2}{r}
+
\frac{1}{r^2}
\Big)
+
\underbrace{
\Big(
\frac{2}{r}
+
\frac{1}{r^2}
\Big)^2}_{\leqslant\,\frac{5}{r^2}}
\bigg]
\notag
\,+
\\
&
\ \ \ \ \
+
\underbrace{
\frac{2}{r^2}
+
\frac{4}{r^4}}_{\leqslant\,\frac{3}{r^2}}
+
(n-2)\,
\bigg[
-\,\frac{1}{r^2}
\underbrace{
-
\frac{1}{2}\,
\frac{1}{r^4}}_{<\,0}
\bigg]
+
(n-3)\,
\bigg[
\underbrace{
\Big(
\frac{2}{r^2}
+
\frac{1}{r^3}
\Big)
+
\Big(
\frac{2}{r^2}
+
\frac{1}{r^3}
\Big)^2}_{
\leqslant\,\frac{4}{r^2}}
\bigg]
\notag
\,+
\\
&
\ \ \ \ \
+
\sum_{k=3}^{n-1}\,
\bigg\{
\underbrace{
\frac{2}{r^k}
+
\frac{4}{r^{2k}}}_{\leqslant\,\frac{3}{r^k}}
+
\underbrace{
(n-k)\,
\log\,
\Big(
1
-
\frac{1}{r^k}
\Big)}_{<\,0}
+
(n-k-1)\,
\bigg[
\underbrace{
\Big(
\frac{2}{r^k}
+
\frac{1}{r^{k+1}}
\Big)
+
\Big(
\frac{2}{r^k}
+
\frac{1}{r^{k+1}}
\Big)^2}_{\leqslant\,\frac{3}{r^k}}
\bigg]
\bigg\}
\notag
\\
&
\,\leqslant\,
\frac{1}{r}\,
\big[
2-n+1+2\,n-4
\big]
+
\frac{1}{r^2}\,
\big[
4
-
{\textstyle{\frac{n}{2}}}
+
{\textstyle{\frac{1}{2}}}
+
n
-
2
+
5\,n
-
10
+
3
-
n
+
2
+
0
+
4\,n
-
12
\big]
\notag
\\
&
\ \ \ \ \
+
\sum_{k=3}^{n-1}\,
\frac{3}{r^k}\,
\big[
1
+
0
+
n-k-1
\big]
\notag
\\
&
\,\leqslant\,
\frac{1}{r}\,
\big[
n-1
\big]
+
\frac{1}{r^2}\,
\big[
{\textstyle{\frac{17}{2}}}\,n
-
{\textstyle{\frac{29}{2}}}
\big]
\,+
\notag
\\
&
\ \ \ \ \
+
\frac{3}{r^3}\,
\frac{1}{1-\frac{1}{r}}\,
\big[
n
\big]
\notag
\\
&
\,\leqslant\,
\frac{n}{r}
+
9\,\frac{n}{r^2}
+
3\,\frac{n}{r^2}.
\qedhere
\end{align}
\endproof

Similarly to the expression:
\[
\widehat{C}(x)
\,=\,
\prod_{k=1}^{n-1}\,
\frac{1}{1-2\,x^k}\,\,
\prod_{k=1}^{n-1}\,
\big(
1
-
x^k
\big)^{n-k}\,\,
\prod_{k=1}^{n-1}\,
\frac{1}{(1-2\,x^k-x^{k+1})^{n-k-1}},
\]
we obtain by simply changing the last sign $-$ to the sign $+$
in the denominator of the third product:
\[
C(x)
\,=\,
\prod_{k=1}^{n-1}\,
\frac{1}{1-2\,x^k}\,\,
\prod_{k=1}^{n-1}\,
\big(
1
-
x^k
\big)^{n-k}\,\,
\prod_{k=1}^{n-1}\,
\frac{1}{(1-2\,x^k+x^{k+1})^{n-k-1}}.
\]

\begin{Lemma}
\label{Lemma-majoration-C-hat-over-C}
One has:
\[
\frac{
\widehat{C}
\big(
\frac{1}{r}
\big)}{
C
\big(
\frac{1}{r}
\big)}
\,\,\leqslant\,\,
e^{\,17\,\frac{n}{r^2}}.
\]
\end{Lemma}

\proof
This quotient writes as:
\[
\frac{
\widehat{C}
\big(
\frac{1}{r}
\big)}{
C
\big(
\frac{1}{r}
\big)}
\,\,=\,\,
\prod_{k=1}^{n-2}\,
\Bigg(
\frac{1-2\,\frac{1}{r^k}
+\frac{1}{r^{k+1}}}{
1-2\,\frac{1}{r^k}-\frac{1}{r^{k+1}}}
\Bigg)^{n-k-1},
\]
since the terms for $k = n-1$ drop. Take logarithm
and use the above majorations:
\[
\!\!\!\!\!\!\!\!\!\!\!\!\!\!\!
\aligned
\log\,
\frac{
\widehat{C}
\big(
\frac{1}{r}
\big)}{
C
\big(
\frac{1}{r}
\big)}
&
\,\,=\,\,
\sum_{k=1}^{n-2}\,
\big(n-k-1\big)\,
\bigg[
\log\,
\Big(
1
-
\Big(
\frac{2}{r^k}
-
\frac{1}{r^{k+1}}
\Big)
\Big)
-
\log\,
\Big(
1
-
\Big(
\frac{2}{r^k}
+
\frac{1}{r^{k+1}}
\Big)
\Big)
\bigg]
\\
&
\,\,=\,\,
(n-2)\,
\bigg[
\log\,
\Big(
1
-
\Big(
\frac{2}{r}
-
\frac{1}{r^2}
\Big)
\Big)
-
\log\,
\Big(
1
-
\Big(
\frac{2}{r}
+
\frac{1}{r^2}
\Big)
\Big)
\bigg]
\,+
\\
&
\ \ \ \ \ 
+
(n-3)\,
\bigg[
\log\,
\Big(
1
-
\Big(
\frac{2}{r^2}
-
\frac{1}{r^3}
\Big)
\Big)
-
\log\,
\Big(
1
-
\Big(
\frac{2}{r^2}
+
\frac{1}{r^3}
\Big)
\Big)
\bigg]
\,+
\\
&
\ \ \ \ \
+
\sum_{k=3}^{n-2}\,
\big(n-k-1\big)\,
\bigg[
\underbrace{
\log\,
\Big(
1
-
\Big(
\frac{2}{r^k}
-
\frac{1}{r^{k+1}}
\Big)
\Big)}_{<\,0}
-
\log\,
\Big(
1
-
\Big(
\frac{2}{r^k}
+
\frac{1}{r^{k+1}}
\Big)
\Big)
\bigg]
\\
&
\,\,\leqslant\,\,
\big(n-2\big)\,
\bigg[
\zero{
-\,
\Big(
\frac{2}{r}}
-
\frac{1}{r^2}
\big)
-
\underbrace{
\frac{1}{2}\,
\Big(
\frac{2}{r}
-
\frac{1}{r^2}
\Big)^2}_{<\,0}
+
\zero{
\Big(
\frac{2}{r}}
+
\frac{1}{r^2}
\Big)
+
\underbrace{
\Big(
\frac{2}{r}
+
\frac{1}{r^2}
\Big)^2}_{\leqslant\,\frac{5}{r^2}}
\bigg]
\,+
\\
&
\ \ \ \ \
+
\big(n-3\big)\,
\bigg[
\underbrace{
-\,
\Big(
\frac{2}{r^2}
-
\frac{1}{r^3}
\big)}_{<\,0}
\underbrace{
-
\frac{1}{2}\,
\Big(
\frac{2}{r^2}
-
\frac{1}{r^3}
\Big)^2}_{<\,0}
+
\underbrace{
\Big(
\frac{2}{r^2}
+
\frac{1}{r^3}
\Big)
+
\Big(
\frac{2}{r^2}
+
\frac{1}{r^3}
\Big)^2}_{\leqslant\,\frac{4}{r^2}}
\bigg]
\,+
\\
\explicationmath{\varepsilon+\varepsilon^2\leqslant2\,\varepsilon}
\ \ \ \ \ \ \ \ \ \ \ \ \ \ \ \ \ \ \ \ \ \ \ \ \ \
&
\ \ \ \ \
+
\sum_{k=3}^{n-2}\,
\big(n-k-1\big)\,
\bigg[
0
+
2\,
\Big(
\underbrace{
\frac{2}{r^k}
+
\frac{1}{r^{k+1}}}_{\leqslant\,\frac{3}{r^k}}
\Big)
\bigg],
\endaligned
\]
and notice, {\em importantly}, that the $\frac{1}{r}$-terms
disappear, so that at the end:
\begin{align}
\log\,
\frac{
\widehat{C}
\big(
\frac{1}{r}
\big)}{
C
\big(
\frac{1}{r}
\big)}
&
\,\leqslant\,
\big(n-2\big)\,
\Big[
\frac{7}{r^2}
\Big]
+
\big(n-3\big)\,
\Big[
\frac{4}{r^2}
\Big]
+
6\,n\,
\sum_{k=3}^\infty\,
\frac{1}{r^k}
\notag
\\
&
\,\leqslant\,
11\,n\,\frac{1}{r^2}
+
6\,n\,
\frac{1}{r^2}\,
\frac{1}{r-1}.
\qedhere
\end{align}
\endproof

Lastly, making the choice:
\[
r
\,:=\,
\sqrt{n}\,
a(n),
\]
with a function $a(n) \underset{n\to\infty}{\longrightarrow} \infty$
tending {\em slowly} to infinity\,\,---\,\,think
$a(n) := \log\,\log\, n$\,\,--- and satisfying at least
$a(n) \ll n^\epsilon$ for any $\epsilon > 0$,
we want to minorize:
\[
\CR_\infty
\,:=\,
C
\big(
{\textstyle{\frac{1}{r}}}
\big)
\,=\,
C
\big(
{\textstyle{\frac{1}{\sqrt{n}\,a(n)}}}
\big).
\]

\begin{Lemma}
\label{Lemma-minoration-C-1-over-r}
One has:
\[
C
\big(
{\textstyle{\frac{1}{\sqrt{n}\,a(n)}}}
\big)
\,\,\geqslant\,\,
e^{\,\frac{1}{2}\,\frac{\sqrt{n}}{a(n)}}
\underset{n\to\infty}{\,\,\,\longrightarrow\,\,}
\infty.
\]
\end{Lemma}

\proof
Take logarithm:
\[
\!\!\!\!\!\!\!\!\!\!\!\!\!\!\!
\aligned
\log\,
C\,\big(
{\textstyle{\frac{1}{r}}}
\big)
&
\,\,=\,\,
\sum_{k=1}^{n-1}\,
\bigg(
-\,\log\,
\Big(
1
-
\frac{2}{r^k}
\Big)
+
\big(n-k\big)\,
\log\,
\Big(
1
-
\frac{1}{r^k}
\Big)
-
\big(
n-k-1
\big)\,\log\,
\Big(
1
-
\Big(
\frac{2}{r^k}
-
\frac{1}{r^{k+1}}
\Big)
\Big)
\bigg)
\\
&
\,\,=\,\,
-\,\log\,
\Big(
1
-
\frac{2}{r}
\Big)
+
\big(n-1\big)\,
\log\,
\Big(
1
-
\frac{1}{r}
\Big)
-
\big(n-2\big)\,
\log\,
\Big(
1
-
\Big(
\frac{2}{r}
-
\frac{1}{r^2}
\Big)
\Big)
\,-
\\
&
\ \ \ \ \
-\,\log\,
\Big(
1
-
\frac{2}{r^2}
\Big)
+
\big(n-2\big)\,
\log\,
\Big(
1
-
\frac{1}{r^2}
\Big)
-
\big(n-3\big)\,
\log\,
\Big(
1
-
\Big(
\frac{2}{r^2}
-
\frac{1}{r^3}
\Big)
\Big)
\,+
\\
&
\ \ \ \ \
+
\sum_{k=3}^{n-1}\,
\bigg\{
-\,\log\,
\Big(
1
-
\frac{2}{r^k}
\Big)
+
\big(
n-k
\big)\,
\log\,
\Big(
1
-
\frac{1}{r^k}
\Big)
-
\big(n-1-k\big)\,
\log\,
\Big(
1
-
\Big(
\frac{2}{r^k}
-
\frac{1}{r^{k+1}}
\Big)
\Big)
\bigg\},
\endaligned
\]
use the minorations:
\[
\aligned
-\,\log\,\big(1-\varepsilon\big)
&
\,\,\geqslant\,\,
\varepsilon,
\\
\log\,\big(1-\delta\big)
&
\,\,\geqslant\,\,
-\,\delta
-
\delta^2,
\endaligned
\]
to get:
\[
\aligned
\log\,
C\big(
{\textstyle{\frac{1}{r}}}
\big)
&
\,\,\geqslant\,\,
\frac{2}{r}
+
\big(n-1\big)\,
\bigg[
-\frac{1}{r}
-
\frac{1}{r^2}
\bigg]
+
\big(n-2\big)\,
\bigg[
\frac{2}{r}
-
\frac{1}{r^2}
\bigg]
\,+
\\
&
\ \ \ \ \
+
\frac{2}{r^2}
+
\big(n-2\big)\,
\bigg[
\underbrace{
-\frac{1}{r^2}
-
\frac{1}{r^4}}_{\geqslant\,-\,\frac{2}{r^2}}
\bigg]
+
\big(n-3\big)\,
\bigg[
\underbrace{
\frac{2}{r^2}
-
\frac{1}{r^3}}_{\geqslant\,\frac{1}{r^2}}
\bigg]
\,+
\\
&
\ \ \ \ \
+
\sum_{k=3}^{n-1}\,
\bigg\{
\frac{2}{r^k}
+
\big(n-k\big)\,
\bigg[
\underbrace{
-\frac{1}{r^k}
-
\frac{1}{r^{2k}}}_{\geqslant\,-\frac{2}{r^k}}
\bigg]
+
\big(n-k-1\big)\,
\bigg[
\underbrace{
\frac{2}{r^k}
-
\frac{1}{r^{k+1}}}_{\geqslant\,\frac{1}{r^k}}
\bigg]
\bigg\}
\\
&
\,\,\geqslant\,\,
\frac{1}{r}\,
\big[
2-n+1+2\,n-4
\big]
+
\frac{1}{r^2}\,
\big[
-n+1-n+2
+
2-2\,n+4+n-3
\big]
\,+
\\
&
\ \ \ \ \
+
\sum_{k=3}^{n-1}\,
\frac{1}{r^k}\,
\big[
2-2\,n+2\,k+n-k-1
\big]
\\
&
\,\,=\,\,
\frac{1}{r}\,\big[n-1\big]
+
\frac{1}{r^2}\,
\big[-3\,n+6\big]
+
\sum_{k=3}^{n-1}\,
\frac{1}{r^k}\,
\big[
\underbrace{
-n+k+1}_{\geqslant\,-n+1}
\big]
\\
&
\,\,\geqslant\,\,
\big(n-1\big)\,
\bigg[
\frac{1}{r}
-
\frac{3}{r^2}
-
\sum_{k=3}^\infty\,
\frac{1}{r^k}
\bigg]
\\
&
\,\,\geqslant\,\,
\big(n-1\big)\,
\bigg[
\frac{1}{r}
-
\frac{4}{r^2}
\bigg].
\endaligned
\]
Lastly, again with $a(n) \ll n^\epsilon$:
\begin{align}
\log\,
C
\bigg(
\frac{1}{\sqrt{n}\,a(n)}
\bigg)
&
\,\,\geqslant\,\,
\frac{n-1}{\sqrt{n}\,a(n)}
\underbrace{
-
\frac{4n-4}{n\,a(n)^2}}_{\geqslant\,-4}
\notag
\\
&
\,\,\geqslant\,\,
\frac{1}{2}\,
\frac{\sqrt{n}}{a(n)}.
\qedhere
\end{align}
\endproof

\Section{\bf Final minorations}
\label{final-minorations}
\HEAD{{\ref{final-minorations}}.~{\sf Final minorations}
}{
Jo\"el {\sc Merker} and The-Anh Ta,
D\'epartement de Math\'ematiques d'Orsay, 
Universit\'e Paris-Sud, France}

As explained at the end of
Section~{\ref{coordinates-t-coordinates-w}}, with a suitable
choice of $r$, the goal is to show:
\[
\aligned
1
&
\overset{\text{\bf ?}}{\,\,\leqslant\,\,}
\sum_{
\substack{
0\leqslant k_2\leqslant n
\\
\ \ \ \ \ 
0\leqslant k_3\leqslant n+k_2
\\
\ \ \ \ \ \
\cdots\cdots\cdots\cdots\cdots\cdots\cdots
\\
\ \ \ \ \ \ \ \ \ \ \ 
0\leqslant k_{n-1}\leqslant n+k_{n-2}
\\
\ \ \ \ \ \ \ \ \ \ \
0\leqslant\,k_n\,\,\,\,\,\leqslant n+k_{n-1}
}}\,\,
C_{k_2,k_3,\dots,k_{n-1},k_n}\,
M_{k_2,k_3,\dots,k_{n-1},k_n}\,
\frac{1}{r^{k_2+k_3+\cdots+k_{n-1}+k_n}}
\\
&
\,\,=:\,\,
\CMR.
\endaligned
\]
Abbreviate this domain range as:
\[
\aligned
\domaink
\,:=\,
\Big\{
&
\big(k_2,k_3,\dots,k_{n-1},k_n\big)
\in
\N^n
\colon\,\,
0
\,\leqslant\,
k_2
\,\leqslant\,
n,
\notag
\\
&
\ \ \ \ \ \ \ \ \ \ \ \ \ \ \ \ \ \ \ \ \ \ \ \ \ \ \ \ \ \ \ \ \ \ \
\ \ \ \ \ \ \ \ \ \ \ \ \ 
0
\,\leqslant\,
k_3
\,\leqslant\,
n+k_2,
\notag
\\
&
\ \ \ \ \ \ \ \ \ \ \ \ \ \ \ \ \ \ \ \ \ \ \ \ \ \ \ \ \ \ \ \ \ \ \
\ \ \ \ \ \ \ \ \ \ \ \ \ 
\cdots\cdots\cdots\cdots\cdots\cdots
\notag
\\
&
\ \ \ \ \ \ \ \ \ \ \ \ \ \ \ \ \ \ \ \ \ \ \ \ \ \ \ \ \ \ \ \ \ \ \
\ \ \ \ \ \ \ \ \ \ \ \ \ 
0
\,\leqslant\,
k_{n-1}
\,\leqslant\,
n+k_{n-2},
\notag
\\
&
\ \ \ \ \ \ \ \ \ \ \ \ \ \ \ \ \ \ \ \ \ \ \ \ \ \ \ \ \ \ \ \ \ \ \
\ \ \ \ \ \ \ \ \ \ \ \ \ 
0
\,\leqslant\,
k_n
\ \ \
\,\leqslant\,
n+k_{n-1}
\Big\}.
\endaligned
\]
Observe that:
\[
\Big\{
k_2+k_3+\cdots+k_{n-1}+k_n
\,\leqslant\,
n
\Big\}
\,\,\subset\,\,
\domaink,
\]
hence {\em a fortiori} with a function $c(n) 
\underset{n\to\infty}{\longrightarrow} \infty$ tending
slowly to infinity to be chosen later:
\[
\Big\{
k_2+k_3+\cdots+k_{n-1}+k_n
\,\leqslant\,
{\textstyle{\frac{\sqrt{n}}{c(n)}}}
\Big\}
\,\,\subset\,\,
\domaink.
\]

Introduce:
\[
\CMR_\TT
\,:=\,
\sum_{k_2+\cdots+k_n\leqslant\frac{\sqrt{n}}{c(n)}}\,
C_{k_2,\dots,k_n}\,
M_{k_2,\dots,k_n}\,
\frac{1}{r^{k_2+\cdots+k_n}},
\]
the letter `$\TT$' standing for `$\TT$runcated', 
with the $\RR$emainder:
\[
\aligned
\CMR
-
\CMR_\TT
&
\,=\,
\sum_{(k_2,\dots,k_n)\in\domaink
\atop
k_2+\cdots+k_n\geqslant1+\frac{\sqrt{n}}{c(n)}}\,
C_{k_2,\dots,k_n}\,
M_{k_2,\dots,k_n}\,
\frac{1}{r^{k_2+\cdots+k_n}}
\\
&
\,=:\,
\CMR_\RR.
\endaligned
\]

Along with these quantities, introduce also:
\reqnomode\usetagform{EngelLie}
\begin{align}
\CMR_\TT^+
&
\,:=\,
\sum_{
k_2+\cdots+k_n\leqslant\frac{\sqrt{n}}{c(n)}
\atop
C_{k_2,\dots,k_n}\geqslant 0}\,
C_{k_2,\dots,k_n}\,
M_{k_2,\dots,k_n}\,
\frac{1}{r^{k_2+\cdots+k_n}}
\tag{(\geqslant 0),}
\\
\CMR_\TT^-
&
\,:=\,
\sum_{
k_2+\cdots+k_n\leqslant\frac{\sqrt{n}}{c(n)}
\atop
\,C_{k_2,\dots,k_n}<0}\,
\big(-C_{k_2,\dots,k_n}\big)\,
M_{k_2,\dots,k_n}\,
\frac{1}{r^{k_2+\cdots+k_n}}
\tag{(\geqslant 0),}
\end{align}
two nonnegative quantities which decompose:
\[
\CMR_\TT
\,=\,
\CMR_\TT^+
-
\CMR_\TT^-.
\]

In addition, without the multinomial-quotient coefficients,
introduce:
\[
\aligned
\CR
&
\,:=\,
\sum_{(k_2,\dots,k_n)\in\domaink}\,
C_{k_2,\dots,k_n}
\cdot 1\cdot
\frac{1}{r^{k_2+\cdots+k_n}},
\\
\CR_\TT
&
\,:=\,
\sum_{k_2+\cdots+k_n\leqslant\frac{\sqrt{n}}{c(n)}}\,
C_{k_2,\dots,k_n}
\cdot 1\cdot
\frac{1}{r^{k_2+\cdots+k_n}},
\\
\CR_\RR
&
\,:=\,
\sum_{(k_2,\dots,k_n)\in\domaink
\atop
k_2+\cdots+k_n\geqslant1+\frac{\sqrt{n}}{c(n)}}
C_{k_2,\dots,k_n}
\cdot 1\cdot
\frac{1}{r^{k_2+\cdots+k_n}},
\endaligned
\]
and similarly also:
\reqnomode\usetagform{EngelLie}
\begin{align}
\CR_\TT^+
&
\,:=\,
\sum_{k_2+\cdots+k_n\leqslant\frac{\sqrt{n}}{c(n)}
\atop
C_{k_2,\dots,k_n}\geqslant 0}\,
C_{k_2,\dots,k_n}
\cdot 1\cdot
\frac{1}{r^{k_2+\cdots+k_n}}
\tag{(\geqslant\,0),}
\\
\CR_\TT^-
&
\,:=\,
\sum_{k_2+\cdots+k_n\leqslant\frac{\sqrt{n}}{c(n)}
\atop
C_{k_2,\dots,k_n}<0}\,
\big(-C_{k_2,\dots,k_n}\big)
\cdot 1\cdot
\frac{1}{r^{k_2+\cdots+k_n}}
\tag{(\geqslant\,0).}
\end{align}

Recall that we are choosing: 
\[
r(n)
\,=\,
\sqrt{n}\,a(n), 
\]
and we now endeavor to 
find a condition guaranteeing that the remainder
$\big\vert \CMR_\RR\big\vert$ be small in absolute value.

To this aim, choose in the Cauchy inequalities 
$\rho := \frac{1}{\sqrt{n}}$, apply 
Lemma~{\ref{Lemma-majoration-C-hat}}:
\[
\widehat{C}
\big(
{\textstyle{\frac{1}{\sqrt{n}}}}
\big)
\,\,\leqslant\,\,
e^{12}\,
e^{\sqrt{n}},
\]
so that Observation~{\ref{Observation-Cauchy-inequalities}} gives:
\[
\widehat{C}_h
\,\,\leqslant\,\,
\frac{1}{\big(\frac{1}{\sqrt{n}}\big)^h}\,
e^{12}\,e^{\sqrt{n}}
\eqno
{\scriptstyle{(\forall\,h\,\geqslant\,0)}}.
\]

Now, majorize the remainder:
\[
\aligned
\big\vert
\CMR_\RR
\big\vert
&
\,\,=\,\,
\Bigg\vert
\sum_{(k_2,\dots,k_n)\in\domaink
\atop
k_2+\cdots+k_n\geqslant1+\frac{\sqrt{n}}{c(n)}}\,
C_{k_2,\dots,k_n}\,
M_{k_2,\dots,k_n}\,
\frac{1}{r^{k_2+\cdots+k_n}}
\Bigg\vert
\\
&
\,\,\leqslant\,\,
\sum_{(k_2,\dots,k_n)\in\domaink
\atop
k_2+\cdots+k_n\geqslant1+\frac{\sqrt{n}}{c(n)}}\,
\big\vert
C_{k_2,\dots,k_n}
\big\vert
\cdot 1\cdot
\frac{1}{r^{k_2+\cdots+k_n}}
\\
&
\,\,\leqslant\,\,
\sum_{k_2+\cdots+k_n\geqslant1+\frac{\sqrt{n}}{c(n)}}\,
\big\vert 
C_{k_2,\dots,k_n}
\big\vert\,
\frac{1}{r^{k_2+\cdots+k_n}}
\\
&
\,\,\leqslant\,\,
\sum_{k_2+\cdots+k_n\geqslant1+\frac{\sqrt{n}}{c(n)}}\,
\widehat{C}_{k_2,\dots,k_n}\,
\frac{1}{r^{k_2+\cdots+k_n}}
\\
&
\,\,=\,\,
\sum_{h=1+\frac{\sqrt{n}}{c(n)}}^\infty\,
\frac{1}{r^h}\,
\sum_{k_2+\cdots+k_n=h}\,
\widehat{C}_{k_2,\dots,k_n}
\\
&
\,\,=\,\,
\sum_{h=1+\frac{\sqrt{n}}{{c(n)}}}^\infty\,
\frac{1}{r^h}\,
\widehat{C}_h,
\endaligned
\]
and hence, thanks to what precedes:
\[
\aligned
\big\vert
\CMR_\RR
\big\vert
&
\,\,\leqslant\,\,
\sum_{h=1+\frac{\sqrt{n}}{c(n)}}^\infty\,
\frac{1}{r^h}\,
\frac{1}{\big(\frac{1}{\sqrt{n}}\big)^h}\,
e^{12}\,e^{\sqrt{n}}
\\
&
\,\,=\,\,
e^{12}\,e^{\sqrt{n}}\,
\sum_{h=1+\frac{\sqrt{n}}{c(n)}}^\infty\,
\frac{1}{
\Big(
\zero{\sqrt{n}}\,a(n)\,\frac{1}{\zero{\sqrt{n}}}
\Big)^h}
\\
&
\,\,=\,\,
e^{12}\,e^{\sqrt{n}}\,
\frac{1}{\big(
a(n)
\big)^{1+\frac{\sqrt{n}}{c(n)}}}\,
\sum_{h=0}^\infty\,
\Big(
\frac{1}{a(n)}
\Big)^h
\\
&
\,\,=\,\,
e^{12}\,e^{\sqrt{n}}\,
e^{\,-\big(1+\frac{\sqrt{n}}{c(n)}\big)\,\logsmall\,a(n)}\,
\underbrace{
\frac{1}{1-\frac{1}{a(n)}}}_{\leqslant\,2}
\\
&
\,\,\leqslant\,\,
2\,e^{12}\,
e^{-\logsmall\,a(n)}\,
e^{\,\sqrt{n}\,\big[1-\frac{\logsmall\,a(n)}{c(n)}\big]}.
\endaligned
\]
In order to insure that the right-hand side is small, since
$e^{-\logsmall\, a(n)} \underset{n\to\infty}{\longrightarrow} 0$,
it suffices to choose:
\[
c(n)
\,:=\,
\log\,a(n)
\underset{n\to\infty}{\,\,\longrightarrow\,\,}
\infty,
\]
to obtain:
\[
\big\vert
\CMR_\RR
\big\vert
\,\,\leqslant\,\,
2\,e^{12}\,
e^{-\logsmall\,a(n)}
\underset{n\to\infty}{\,\,\longrightarrow\,\,}
0.
\]

\begin{Lemma}
\label{Lemma-CMR-R-small}
With $c(n) = \log\, a(n)$, it holds:
\reqnomode\usetagform{EngelLie}
\begin{align}
\big\vert
\CMR_\RR
\big\vert
&
\,\,\leqslant\,\,
\sum_{k_2+\cdots+k_n\geqslant1+\frac{\sqrt{n}}{c(n)}}\,
\widehat{C}_{k_2,\dots,k_n}\,
\frac{1}{r^{k_2+\cdots+k_n}}
\notag
\\
&
\,\,\leqslant\,\,
2\,e^{12}\,
e^{-\logsmall\,a(n)}.
\tag{\qed}
\end{align}
\end{Lemma}

Further, it is now necessary to estimate the size of the first
terms $\CMR_\TT$, and to show that they are large. It will be useful
that:
\[
\frac{\sqrt{n}}{\log\,a(n)}
\,=\,
{\rm o}
\big(\sqrt{n}\big)
\eqno
{\scriptstyle{(n\,\longrightarrow\,\infty)}}.
\]
Introduce the quantities:
\[
\aligned
\CR_\infty
&
\,:=\,
\sum_{k_2,\dots,k_n\geqslant 0}\,
C_{k_2,\dots,k_n}\,
\frac{1}{r^{k_2+\cdots+k_n}}
\,\,=\,\,
C^{n-1}
\big(
{\textstyle{\frac{1}{r}}}
\big),
\\
\CR_\infty^+
&
\,:=\,
\sum_{k_2,\dots,k_n\geqslant 0
\atop
C_{k_2,\dots,k_n\geqslant 0}}\,
C_{k_2,\dots,k_n}\,
\frac{1}{r^{k_2+\cdots+k_n}},
\\
\CR_\infty^-
&
\,:=\,
\sum_{k_2,\dots,k_n\geqslant 0
\atop
C_{k_2,\dots,k_n<0}}\,
\big(-C_{k_2,\dots,k_n}\big)\,
\frac{1}{r^{k_2+\cdots+k_n}},
\\
\widehat{\CR}_\infty
&
\,:=\,
\sum_{k_2,\dots,k_n\geqslant 0}\,
\widehat{C}_{k_2,\dots,k_n}\,
\frac{1}{r^{k_2+\cdots+k_n}}
\,\,=\,\,
\widehat{C}^{n-1}
\big(
{\textstyle{\frac{1}{r}}}
\big),
\endaligned
\]
for which it is clear that:
\[
\CR_\infty^+
+
\CR_\infty^-
\,\,\leqslant\,\,
\widehat{\CR}_\infty.
\]

By Lemma~{\ref{Lemma-majoration-C-hat-over-C}}:
\[
\frac{\widehat{\CR}_\infty}{\CR_\infty}
\,\,\leqslant\,\,
e^{\,\frac{17}{a(n)^2}},
\]
and next:
\[
\CR_\infty^+
+
\CR_\infty^-
\,\,\leqslant\,\,
\widehat{\CR}_\infty
\,\,\leqslant\,\,
e^{\,\frac{17}{a(n)^2}}\,
\CR_\infty
\,\,=\,\,
e^{\,\frac{17}{a(n)^2}}\,
\Big(
\CR_\infty^+
-
\CR_\infty^-
\Big),
\]
from which it comes:
\[
\CR_\infty^-\,
\Big(
\underbrace{
1
+
e^{\,\frac{17}{a(n)^2}}}_{\geqslant\,2}
\Big)
\,\,\leqslant\,\,
\Big(
e^{\,\frac{17}{a(n)^2}}
-
1
\Big)\,
\CR_\infty^+,
\]
whence:
\leqnomode\usetagform{default}
\begin{align}
\label{majoration-CR-minus-infty}
\CR_\infty^-
\,\,\leqslant\,\,
{\textstyle{\frac{1}{2}}}\,
\Big(
e^{\,\frac{17}{a(n)^2}}
-
1
\Big)\,
\CR_\infty^+.
\end{align}

Next, we want to minorize $\CMR_\TT$ in order to show it is large:
\[
\aligned
\CMR_\TT
&
\,\,=\,\,
\CMR_\TT^+
-
\CMR_\TT^-
\\
&
\,\,=\,\,
\sum_{k_2+\cdots+k_n\leqslant\frac{\sqrt{n}}{c(n)}
\atop
C_{k_2,\dots,k_n}\geqslant 0}\,
C_{k_2,\dots,k_n}\,
M_{k_2,\dots,k_n}\,
\frac{1}{r^{k_2+\cdots+k_n}}
\,-
\\
&
\ \ \ \ \
-
\sum_{k_2+\cdots+k_n\leqslant\frac{\sqrt{n}}{c(n)}
\atop
C_{k_2,\dots,k_n}<0}\,
\big(
-\,
C_{k_2,\dots,k_n}
\big)\,
M_{k_2,\dots,k_n}\,
\frac{1}{r^{k_2+\cdots+k_n}}
\\
\explicationtext{Proposition~{\ref{Proposition-minoration-M}}}
\ \ \ \ \ \ \ \ \ \ \ \ \ \ \ \ \ \ \ \ \ \ \ \ \ \
&
\,\,\geqslant\,\,
\sum_{k_2+\cdots+k_n\leqslant\frac{\sqrt{n}}{c(n)}
\atop
C_{k_2,\dots,k_n}\geqslant0}\,
C_{k_2,\dots,k_n}\,
e^{-\frac{2}{c(n)^2}}\,
\frac{1}{r^{k_2+\cdots+k_n}}
\,-
\\
\explicationtext{Lemma~{\ref{Lemma-M-less-1}}}
\ \ \ \ \ \ \ \ \ \ \ \ \ \ \ \ \ \ \ \ \ \ \ \ \ \
&
\ \ \ \ \
-
\sum_{k_2+\cdots+k_n\leqslant\frac{\sqrt{n}}{c(n)}
\atop
C_{k_2,\dots,k_n}<0}\,
\big(
-C_{k_2,\dots,k_n}
\big)
\cdot
1
\cdot
\frac{1}{r^{k_2+\cdots+k_n}}
\\
&
\,\,=\,\,
e^{-\frac{2}{c(n)^2}}\,
\CR_\TT^+
-
\CR_\TT^-,
\endaligned
\]
but we yet need to compare these to the quantities 
$\CR_\infty^\pm$. Hence we estimate:
\[
\aligned
\big\vert
\CR_\infty^+
-
\CR_\TT^+
\big\vert
\,\,=\,\,
\CR_\infty^+
-
\CR_\TT^+
&
\,\,=\,\,
\sum_{k_2+\cdots+k_n\geqslant 1+\frac{\sqrt{n}}{c(n)}
\atop
C_{k_2,\dots,k_n\geqslant 0}}\,
C_{k_2,\dots,k_n}\,
\frac{1}{r^{k_2+\cdots+k_n}}
\\
&
\,\,\leqslant\,\,
\sum_{k_2+\cdots+k_n\geqslant1+\frac{\sqrt{n}}{c(n)}}\,
\widehat{C}_{k_2,\dots,k_n}\,
\frac{1}{r^{k_2+\cdots+k_n}}
\\
&
\,\,\leqslant\,\,
2\,e^{12}\,e^{-\logsmall\,a(n)}
\endaligned
\]
and more simply:
\[
-\,\CR_\TT^-
\,\,\geqslant\,\,
-\,
\CR_\infty^-,
\]
since:
\[
0
\,\,\leqslant\,\,
\CR_\infty^-
-
\CR_\TT^-
\,\,=\,\,
\sum_{k_2+\cdots+k_n\geqslant 1+\frac{\sqrt{n}}{c(n)}
\atop
C_{k_2,\dots,k_n<0}}\,
\big(-C_{k_2,\dots,k_n}\big)\,
\frac{1}{r^{k_2+\cdots+k_n}}.
\]

Thanks to all this:
\[
\aligned
\CMR_\TT
&
\,\,\geqslant\,\,
e^{-\frac{2}{c(n)^2}}\,
\CR_\TT^+
-
\CR_\TT^-
\\
&
\,\,\geqslant\,\,
e^{-\frac{2}{c(n)^2}}\,
\Big[
\CR_\infty^+
-
2\,e^{12}\,
e^{-\logsmall\,a(n)}
\Big]
-
\CR_\infty^-
\endaligned
\]
hence applying the minoration~({\ref{majoration-CR-minus-infty}})
for $-\, \CR_\infty^-$:
\[
\aligned
\CMR_\TT
&
\,\,\geqslant\,\,
e^{-\frac{2}{c(n)^2}}\,
\Big[
\CR_\infty^+
-
2\,e^{12}\,
e^{-\logsmall\,a(n)}
\Big]
-
{\textstyle{\frac{1}{2}}}\,
\Big(
e^{\frac{17}{a(n)^2}}
-
1
\Big)\,
\CR_\infty^+
\\
&
\,\,=\,\,
\CR_\infty^+\,
\Big[
\underbrace{
e^{-\frac{2}{(\logsmall\,a(n))^2}}}_{
\underset{n\to\infty}{\longrightarrow}\,1}
-
{\textstyle{\frac{1}{2}}}\,
\big(
\underbrace{
e^{\frac{17}{a(n)^2}}
-
1}_{
\underset{n\to\infty}{\longrightarrow}\,0}
\big)
\Big]
-
\underbrace{
2\,e^{12}\,
e^{-\logsmall\,a(n)}\,
e^{-\frac{2}{(\logsmall\,a(n))^2}}}_{
\underset{n\to\infty}{\longrightarrow}\,0}.
\endaligned
\]

Since trivially:
\[
\aligned
\CR_\infty^+
&
\,=\,
\CR_\infty
+
\CR_\infty^-
\\
&
\,\geqslant\,
\CR_\infty,
\endaligned
\]
it comes:
\[
\CMR_\TT
\,\,\geqslant\,\,
\CR_\infty\,
\Big[
e^{-\frac{2}{(\logsmall\,a(n))^2}}
-
{\textstyle{\frac{1}{2}}}\,
\big(
e^{\frac{17}{a(n)^2}}
-
1
\big)
\Big]
-
2\,e^{12}\,
e^{-\logsmall\,a(n)}\,
e^{-\frac{2}{(\logsmall\,a(n))^2}},
\]
whence using Lemma~{\ref{Lemma-minoration-C-1-over-r}}:
\[
\CMR_\TT
\,\,\geqslant\,\,
e^{\frac{1}{2}\,\frac{\sqrt{n}}{a(n)}}\,
\Big[
e^{-\frac{2}{(\logsmall\,a(n))^2}}
-
{\textstyle{\frac{1}{2}}}\,
\big(
e^{\frac{17}{a(n)^2}}
-
1
\big)
\Big]
-
2\,
e^{12}\,
e^{-\logsmall\,a(n)}\,
e^{-\frac{2}{(\logsmall\,a(n))^2}}.
\]
Coming back to:
\[
\aligned
\CMR
&
\,\geqslant\,
\CMR_\TT
-
\big\vert
\CMR_\RR
\big\vert
\\
&
\,\geqslant\,
\CMR_\TT
-
2\,e^{12}\,
e^{-\logsmall\,a(n)}\,
e^{-\frac{2}{(\logsmall\,a(n))^2}},
\endaligned
\]
we obtain finally:
\[
\CMR
\,\,\geqslant\,\,
e^{\frac{1}{2}\,\frac{\sqrt{n}}{a(n)}}\,
\Big[
e^{-\frac{2}{(\logsmall\,a(n))^2}}
-
{\textstyle{\frac{1}{2}}}\,
\big(
e^{\frac{17}{a(n)^2}}
-
1
\big)
\Big]
-
2\,
e^{12}\,
e^{-\logsmall\,a(n)}\,
\Big(
1
+
e^{-\frac{2}{(\logsmall\,a(n))^2}}
\Big).
\]
This minorant is $\geqslant 1$ for all $n \gg \NN$ large enough.

\proof[Proof of Theorem~{\ref{Theorem-GG-sqrt-n-log-log-n}}]
Choose\,\,---\,\,think integer value\,\,---:
\[
r
\,:=\,
\sqrt{n}\,
\frac{\logsmall\,n}{2}.
\]
Proposition~{\ref{Proposition-r-3-n}} concludes for $\X^{n-1}
\subset \P^n$ of degree at least:
\begin{align}
d_{\GG\GG}(n)
\,:=\,
&\,
\Big(
\sqrt{n}\,
\frac{\logsmall\,n}{2}
+
3
\Big)^n\,
25\,n^2
\notag
\\
\,=\,
&\,
\Big(
\sqrt{n}\,\log\,n
\Big)^n\,\,
\frac{1}{2^n}\,\,
\bigg(
1
+
\frac{6}{\sqrt{n}\,\logsmall\,n}
\bigg)^n\,\,
25\,n^2
\notag
\\
\,=\,
&\,
\Big(
\sqrt{n}\,\log\,n
\Big)^n\,\,
e^{-n\,\logsmall\,2}\,
e^{n\,\logsmall\,\big(1+\frac{6}{\sqrt{n}\,\logsmall\,n}\big)}\,
25\,n^2
\notag
\\
\,\leqslant\,
&\,
\Big(
\sqrt{n}\,\log\,n
\Big)^n\,\,
\underbrace{
e^{-n\,\logsmall\,2
+\frac{6\,\sqrt{n}}{\logsmall\,n}}\,\,
25\,n^2}_{<\,1\,\,\text{\rm when}\,\,n\,\geqslant\,\NN_{\GG\GG}}.
\qedhere
\end{align}
\endproof

\proof[Proof of Theorem~{\ref{Theorem-Kb-n-log-n}}]
Thanks to~{\cite{Riedl-Yang-2018}},
Proposition~{\ref{Proposition-r-3-n}} concludes for $\X^{n-1}
\subset \P^n$ of degree at least:
\begin{align}
d_\KK(n)
\,:=\,
&\,
\bigg(
\sqrt{2\,n}\,\,
\frac{\logsmall\,\logsmall\,(2\,n)}{2}
+
3
\bigg)^{2n}\,
25\,\big(2\,n\big)^2
\notag
\\
\,=\,
&\,
\bigg(
\frac{\sqrt{n}}{\sqrt{2}}\,
\log\,\log\,(2\,n)
\bigg)^{2n}\,
\bigg(
1
+
\frac{3\,\sqrt{2}}{\sqrt{n}\,\logsmall\,\logsmall\,(2\,n)}
\bigg)^{2n}\,
100\,n^2
\notag
\\
\,=\,
&\,
\Big(
n\,\log\,n
\Big)^n\,\,
\frac{\big(\logsmall\,\logsmall\,(2\,n)\big)^{2n}}{
\big(\logsmall\,n\big)^n}\,\,
\frac{1}{\sqrt{2}^{\,2n}}\,\,
\bigg(
1
+
\frac{3\,\sqrt{2}}{\sqrt{n}\,\logsmall\,\logsmall\,(2\,n)}
\bigg)^{2n}\,\,
100\,n^2
\notag
\\
\,=\,
&\,
\Big(
n\,\log\,n
\Big)^n\,\,
e^{\,
2n\,\logsmall\,\logsmall\,\logsmall\,(2\,n)
-
n\,\logsmall\,\logsmall\,n
-
n\,\logsmall\,2
+
2n\,\logsmall\,
\big(
1
+
\frac{3\,\sqrt{2}}{\sqrt{n}\,\logsmall\,\logsmall\,(2\,n)}
\big)}\,\,
100\,n^2
\notag
\\
\,\leqslant\,
&\,
\Big(
n\,\log\,n
\Big)^n\,\,
\underbrace{
e^{\,
2n\,\logsmall\,\logsmall\,\logsmall\,(2\,n)
-
n\,\logsmall\,\logsmall\,n
-
n\,\logsmall\,2
+
\frac{2\,\sqrt{n}\,3\,\sqrt{2}}{\logsmall\,\logsmall\,(2\,n)}}\,\,
100\,n^2}_{
<\,1\,\,\text{\rm when}\,\,n\,\geqslant\,\NN_\KK}.
\qedhere
\end{align}
\endproof

\Section{\bf Some inequalities on circles 
$\{ \vert z\vert = \rho\}$}
\label{vert-z-vert-rho}
\HEAD{{\ref{vert-z-vert-rho}}.~Some inequalities on circles 
$\{ \vert z\vert = \rho\}$
}{
Jo\"el {\sc Merker} and The-Anh Ta,
D\'epartement de Math\'ematiques d'Orsay, 
Universit\'e Paris-Sud, France}

\begin{Proposition}
\label{Proposition-G-k-z-H-ell-z}
For every radius $0 \leqslant \rho \leqslant 
\frac{1}{4}$, the functions:
\reqnomode\usetagform{EngelLie}
\begin{align}
G_k(z)
&
\,:=\,
\frac{1-z^k}{1-2\,z^k}
\tag{(k\,\geqslant\,1),}
\\
H_\ell(z)
&
\,:=\,
\frac{1-z^\ell}{1-2\,z^\ell-z^{\ell+1}}
\tag{(\ell\,\geqslant\,1)},
\end{align}
attain, on the circle $\big\{z \in \C \colon
\vert z\vert = \rho\big\}$,
their maximum modulus at the real point $z = \rho$:
\reqnomode\usetagform{EngelLie}
\begin{align}
\underset{\vert z\vert=\rho}{\max}\,\,
\bigg\vert
\frac{1-z^k}{1-2\,z^k}
\bigg\vert
&
\,\,=\,\,
\frac{1-\rho^k}{1-2\,\rho^k}
\tag{(\forall\,k\,\geqslant\,1),}
\\
\underset{\vert z\vert=\rho}{\max}\,\,
\bigg\vert
\frac{1-z^\ell}{1-2\,z^\ell-z^{\ell+1}}
\bigg\vert
&
\,\,=\,\,
\frac{1-\rho^\ell}{1-2\,\rho^\ell-\rho^{\ell+1}}
\tag{(\forall\,\ell\,\geqslant\,1),}
\end{align}
and with the choice $\rho := 0.25$, the graphs 
on the unit circle of the two quotient functions:
\[
\theta
\,\longmapsto\,
\frac{\big\vert G_k(\rho\,e^{i\theta})\big\vert}{G_k(\rho)}
\ \ \ \ \ \ \ \ \ \ \ \ \ \ \ \ \ \
\text{and}
\ \ \ \ \ \ \ \ \ \ \ \ \ \ \ \ \ \
\theta
\,\longmapsto\,
\frac{\vert H_\ell(\rho\,e^{i\theta})\big\vert}{H_\ell(\rho)}
\eqno
{\scriptstyle{(-\pi\,\leqslant\,\theta\,\leqslant\,\pi)}}
\]
show up, respectively, for the three choices $k = 2, 5, 10$
and the three choices $\ell = 2, 5, 10$, as:

\medskip

\includegraphics[scale=0.30]{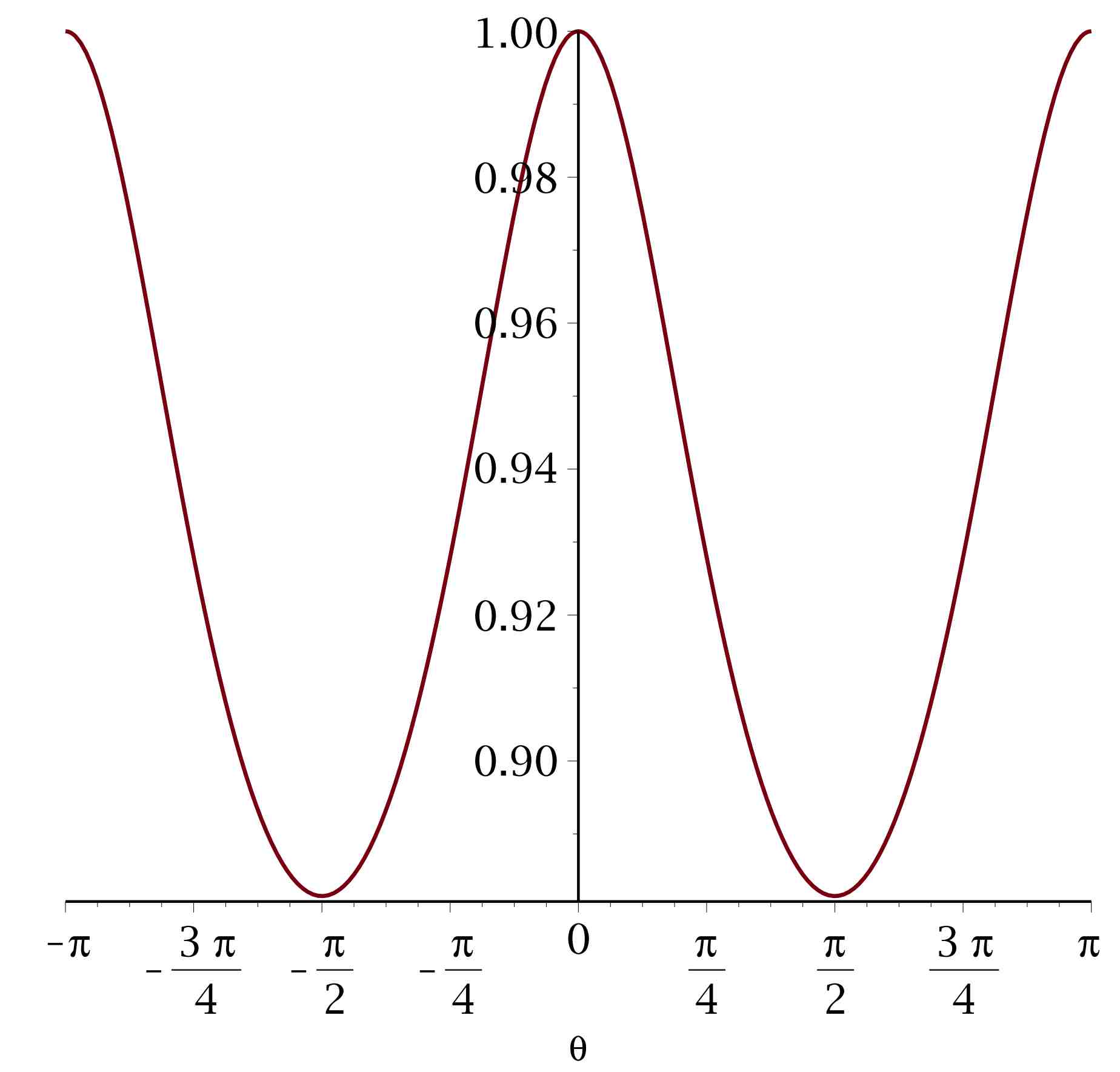}
\ \ \ \ \ \ \ 
\includegraphics[scale=0.30]{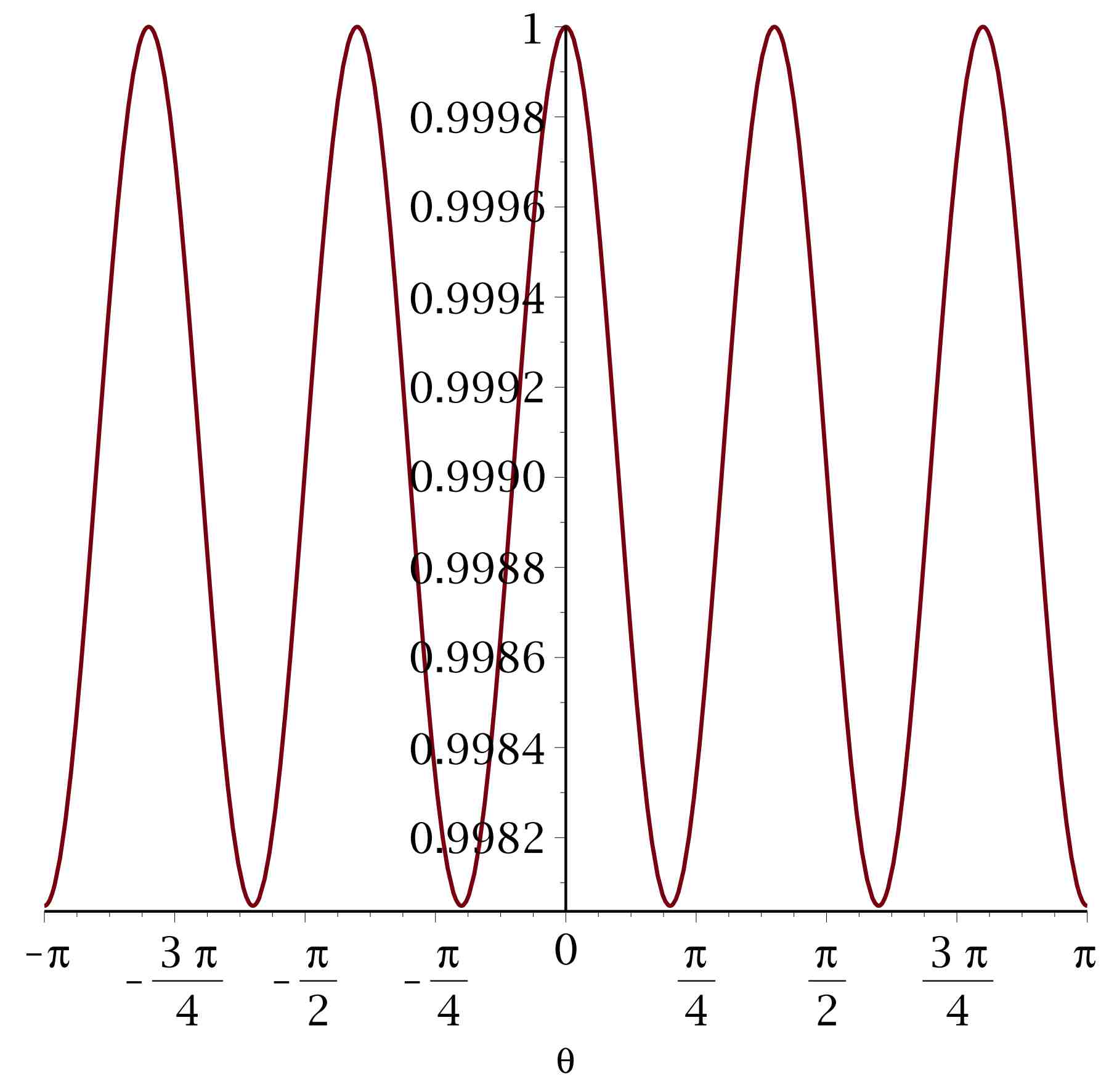}
\ \ \ \ \ \ \ 
\includegraphics[scale=0.30]{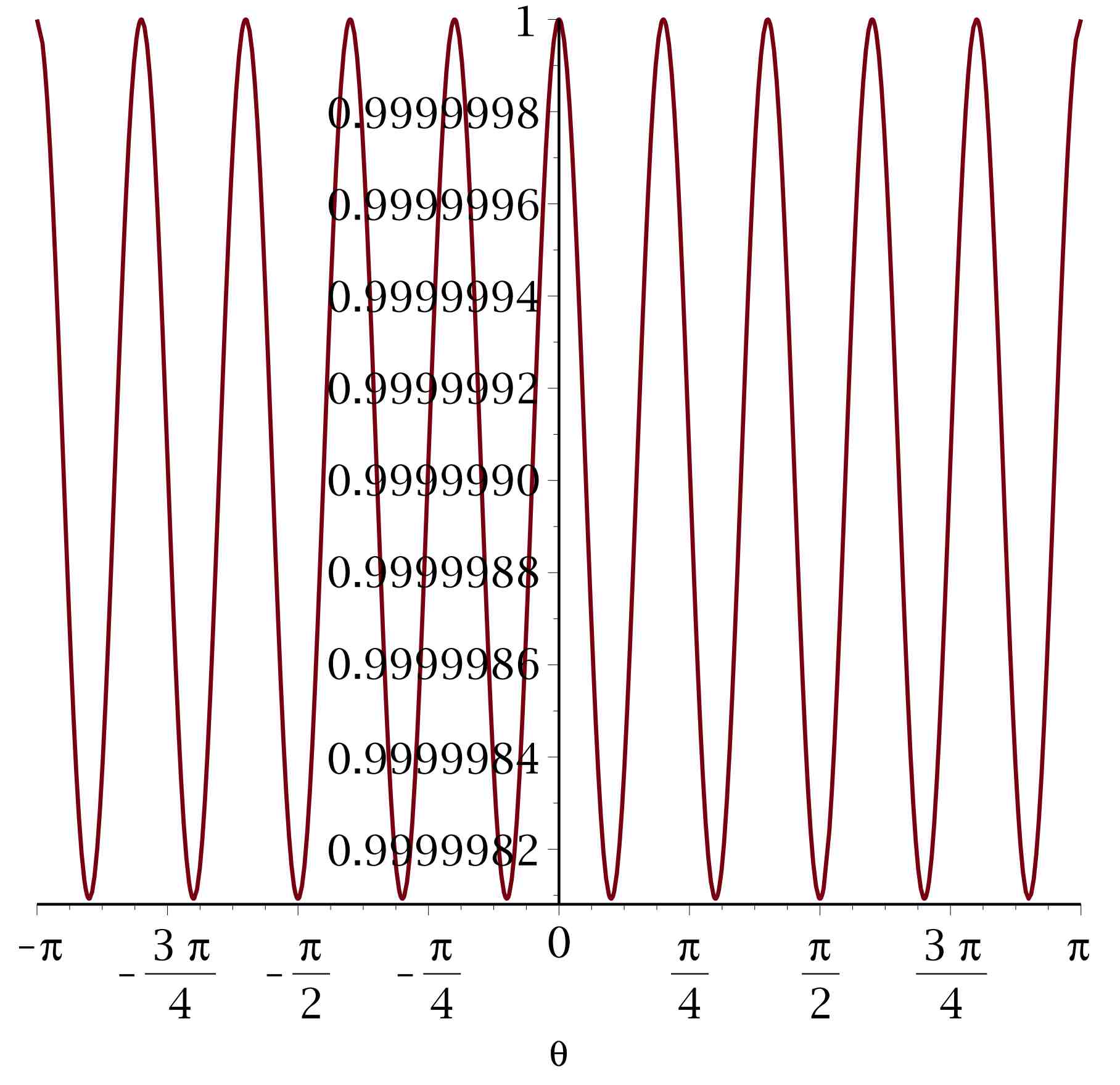}

\includegraphics[scale=0.30]{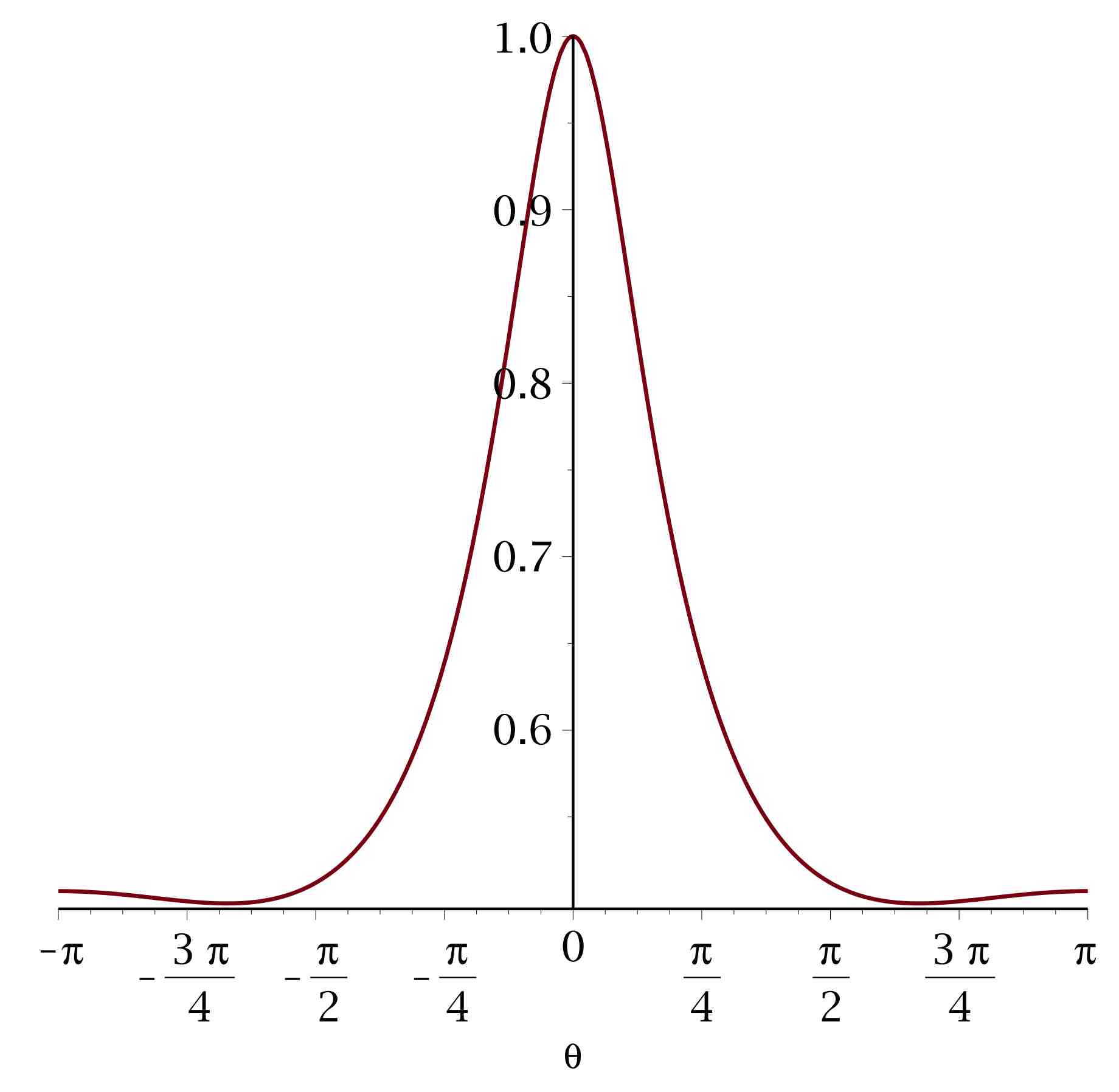}
\ \ \ \ \ \ \ 
\includegraphics[scale=0.30]{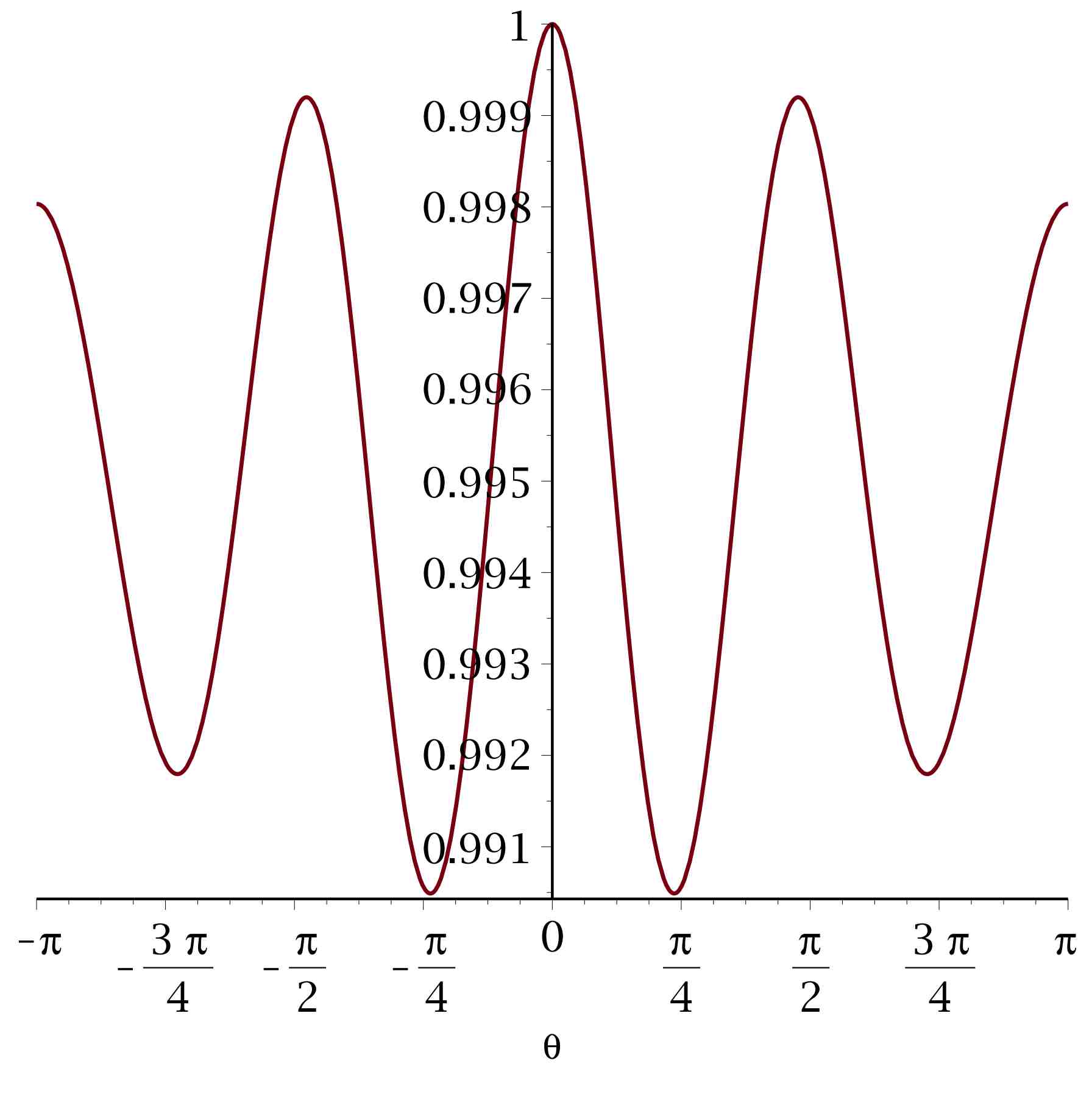}
\ \ \ \ \ \ \ 
\includegraphics[scale=0.30]{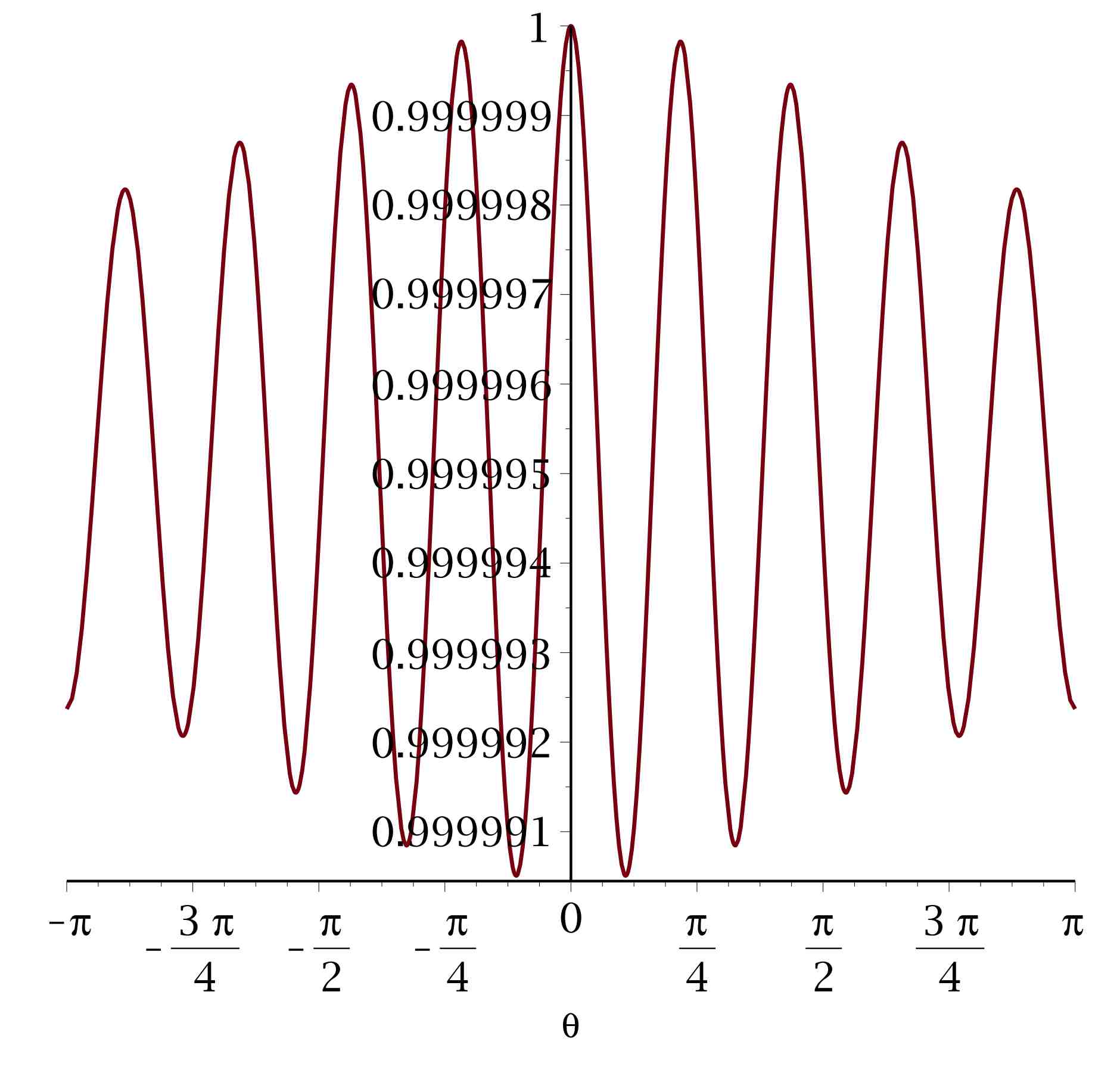}

\end{Proposition}

\proof
Treat at first the $G_k\big( \rho\, e^{i\theta} \big)$
with $\theta \in \R$, by squaring:
\[
\frac{\big\vert 1-\rho^k\,e^{ik\theta}\big\vert^2}{
\big\vert 1-2\,\rho^k\,e^{ik\theta}\big\vert^2}
\overset{\text{\bf ?}}{\,\,\leqslant\,\,}
\frac{\big(1-\rho^k\big)^2}{\big(1-2\,\rho^k\big)^2}
\eqno
{\scriptstyle{(\forall\,\theta\,\in\,\R)}},
\]
that is to say:
\[
\frac{
\big(1-\rho^k\,\cos\,k\theta\big)^2
+
\big(-\,\rho^k\,\sin\,k\theta\big)^2
}{
\big(1-2\,\rho^k\,\cos\,k\theta\big)^2
+
\big(-\,2\,\rho^k\,\sin\,k\theta\big)^2}
\overset{\text{\bf ?}}{\,\,\leqslant\,\,}
\frac{1-2\,\rho^k+\rho^{2k}}{1-4\,\rho^k+4\,\rho^{2k}},
\]
or equivalently, after crossing/clearing the fractions:
\[
\aligned
0
&
\overset{\text{\bf ?}}{\,\,\leqslant\,\,}
\big(
1-2\,\rho^k+\rho^{2k}
\big)\,
\big(
1-4\,\rho^k\,\cos\,k\theta+4\,\rho^{2k}
\big)
-
\big(
1-4\,\rho^k+4\,\rho^{2k}
\big)\,
\big(
1-2\,\rho^k\,\cos\,k\theta+\rho^{2k}
\big)
\\
&
\,\,=\,\,
\zero{1}
-
4\,\rho^k\,\cos\,k\theta
+
\zero{4\,\rho^{2k}}
\\
&
\ \ \ \ \ \ \ \ \ \ \ 
-\,2\,\rho^k
\ \ \ \ \ \ \ \ \ \ \
+\zero{8\,\rho^{2k}\,\cos\,k\theta}
-
8\,\rho^{3k}
\\
&
\ \ \ \ \ \ \ \ \ \ \ \ \ \ \ \ \ \ \ \ \ \ \ \ \ \ \ \ \ \ \ \ \ \ 
+\zero{\rho^{2k}}
\ \ \ \ \ \ \ \ \ \ \ \ \
-\,4\,\rho^{3k}\,\cos\,k\theta
+
\zero{4\,\rho^{4k}}
\\
&
\ \ \ 
-\,\zero{1}
+2\,\rho^k\,\cos\,k\theta
-\zero{\rho^{2k}}
\\
&
\ \ \ \ \ \ \ \ \ \ \ \ \,
+4\,\rho^k
\ \ \ \ \ \ \ \ \ \
-\zero{8\,\rho^{2k}\,\cos\,k\theta}
+
4\,\rho^{3k}
\\
&
\ \ \ \ \ \ \ \ \ \ \ \ \ \ \ \ \ \ \ \ \ \ \ \ \ \ \ \ \ \ \ \ \ \
-\,\zero{4\,\rho^{2k}}
\ \ \ \ \ \ \ \ \ \
+8\,\rho^{3k}\,\cos\,k\theta
-\zero{4\,\rho^{4k}}.
\endaligned
\]
Visibly, $5 \cdot 2 = 10$ underlined terms annihilate by pairs:
\[
\aligned
0
&
\overset{\text{\bf ?}}{\,\,\leqslant\,\,}
2\,\rho^k
-
2\,\rho^k\,\cos\,k\theta
-
4\,\rho^{3k}
+
4\,\rho^{3k}\,\cos\,k\theta
\\
&
\,\,=\,\,
2\,\rho^k\,
\Big[
1
-
\cos\,k\theta
-
2\,\rho^{2k}
+
2\,\rho^{2k}\,\cos\,k\theta
\Big],
\endaligned
\]
and by luck, the obtained expression factorizes under a form
which shows well that it takes only nonnegative values
because $0 \leqslant \rho \leqslant 0.25$:
\[
0
\overset{\text{\sf yes}}{\,\,\leqslant\,\,}
2\,\rho^k\,
\big(
1
-
2\,\rho^{2k}
\big)\,
\big(
1
-
\cos\,k\theta
\big)
\eqno
{\scriptstyle{(\forall\,k\,\geqslant\,1,\,\,
\forall\,\theta\,\in\,\R)}}.
\]

Secondly, for the functions $\big( H_\ell(z) \big)_{\ell \geqslant 1}$,
no such pleasant factorization is available. One can then 
view these $H_\ell(z)$ as `perturbations' of the $G_\ell(z)$, with
the addition of $-\, z^{\ell+1}$ at the denominator. More precisely,
starting from the desired inequality of which we take the 
squared modulus:
\[
\frac{\big\vert 1-\rho^\ell\,e^{i\ell\theta}\big\vert^2}{
\big\vert1-2\,\rho^\ell\,e^{i\ell\theta}
-\underline{\rho^{\ell+1}\,e^{i(\ell+1)\,\theta}}\big\vert^2}
\overset{\text{\bf ?}}{\,\,\leqslant\,\,}
\frac{\big(1-\rho^\ell\big)^2}{
\big(1-2\,\rho^\ell-\underline{\rho^{\ell+1}}\big)^2},
\]
the `perturbing terms' being underlined, after
crossing/clearing the denominators, we are led to establish
an inequality which is a `perturbation' of the one just
done above:
\[
\aligned
0
&
\overset{\text{\bf ?}}{\,\,\leqslant\,\,}
\big(
1-2\,\rho^\ell+\rho^{2\ell}
\big)\,
\Big[
\big(
1-2\,\rho^\ell\,\cos\,\ell\theta-
\underline{\rho^{\ell+1}\,\cos\,(\ell+1)\theta}
\big)^2
+
\big(
2\,\rho^\ell\,\sin\,\ell\theta
+
\underline{\rho^{\ell+1}\,\sin\,(\ell+1)\theta}
\big)^2
\Big]
\\
&
-\,
\big(
1-4\,\rho^\ell+4\,\rho^{2\ell}-
\underline{2\,\rho^{\ell+1}+4\,\rho^{2\ell+1}+\rho^{2\ell+2}}
\big)
\Big[
\big(
1-\rho^\ell\,\cos\,\ell\theta
\big)^2
+
\big(
\rho^\ell\,\sin\,\ell\theta
\big)^2
\Big],
\endaligned
\]
the perturbing terms being still underlined, that is to say:
\[
\aligned
0
&
\overset{\text{\bf ?}}{\,\,\leqslant\,\,}
\big(
1-2\,\rho^\ell+\rho^{2\ell}
\big)\,
\Big[
1
-
4\,\rho^\ell\,\cos\,\ell\theta
+
4\,\rho^{2\ell}\,\cos^2\ell\theta
+
4\,\rho^{2\ell}\,\sin^2\ell\theta
\\
&
\ \ \ \ \ \ \ \ \ \ \ \ \ \ \ \ \ \ \ \ 
-\,
\underline{
2\,\rho^{\ell+1}\,\cos\,(\ell+1)\theta
+
4\,\rho^{2\ell+1}\,\cos\,\ell\theta\,\cos\,(\ell+1)\theta
+
\rho^{2\ell+2}\cos^2(\ell+1)\theta}
\\
&
\ \ \ \ \ \ \ \ \ \ \ \ \ \ \ \ \ \ \ \ \ \ \ \ \ \ \ \ \ 
+
\underline{
4\,\rho^{2\ell+1}\,\sin\,\ell\theta\,\sin\,(\ell+1)\theta
+
\rho^{2\ell+2}\,\sin^2(\ell+1)\theta}
\Big]
\\
&
-\,
\big(
1-4\,\rho^\ell+4\,\rho^{2\ell}-
\underline{2\,\rho^{\ell+1}+4\,\rho^{2\ell+1}+\rho^{2\ell+2}}
\big)\,
\Big[
1
-
2\,\rho^\ell\,\cos\,\ell\theta
+
\rho^{2\ell}\,\cos^2\ell\theta
+
\rho^{2\ell}\,\sin^2\ell\theta
\Big].
\endaligned
\]
Without redoing the calculation concerning the (principal, 
not underlined) terms, and using:
\[
\cos\,(\ell+1)\theta\,\cos\,(-\ell\theta)
-
\sin\,(\ell+1)\theta\,\sin\,(-\ell\theta)
\,\,=\,\,
\cos\,\big(
(\ell+1-\ell)\,\theta
\big),
\]
we obtain:
\[
\aligned
0
&
\overset{\text{\bf ?}}{\,\,\leqslant\,\,}
2\,\rho^\ell\,
\big(
1-2\,\rho^{2\ell}
\big)\,
\big(
1
-
\cos\,\ell\theta
\big)
\\
&
\ \ \ \ \
+
\big(
1-2\,\rho^\ell+\rho^{2\ell}
\big)\,
\Big[
-2\,\rho^{\ell+1}\,
\cos\,(\ell+1)\theta
+
4\,\rho^{2\ell+1}\,\cos\,\theta
+
\rho^{2\ell+2}
\Big]
\\
&
\ \ \ \ \
+
\big(
2\,\rho^{\ell+1}-4\,\rho^{2\ell+1}-\rho^{2\ell+2}
\big)\,
\Big[
1
-
2\,\rho^\ell\,\cos\,\ell\theta
+
\rho^{2\ell}
\Big].
\endaligned
\]
Now, organize the expansion of lines 2 and 3 in a convenient 
synoptic way:
\[
\footnotesize
\aligned
&
-\,2\,\rho^{\ell+1}\,\cos\,(\ell+1)\theta
+4\,\rho^{2\ell+1}\,\cos\,\theta 
\ \ \ \ \ \ \ \
+\zero{\rho^{2\ell+2}}
\\
&
\ \ \ \ \ \ \ \ \ \ \ \ \ \ \ \ \ \ \ \ \ \ \ \ \ \ \ \ \ \ \ \ \ \ \
+4\,\rho^{2\ell+1}\,\cos\,(\ell+1)\theta 
\ \ \ \ \ \ \ 
-8\,\rho^{3\ell+1}\,\cos\,\theta
-
2\,\rho^{3\ell+2}
\\
&
\ \ \ \ \ \ \ \ \ \ \ \ \ \ \ \ \ \ \ \ \ \ \ \ \ \ \ \ \ \ \ \ \ \ \ 
\ \ \ \ \ \ \ \ \ \ \ \ \ \ \ \ \ \ \ \ \ \ \ \ \ \ \ \ \ \ \ \ \ \ \ 
\ \ \ \ \ \ \ \ \ \ 
-\,2\,\rho^{3\ell+1}\,\cos\,(\ell+1)\theta
+4\,\rho^{4\ell+1}\,\cos\,\theta
+
\zero{\rho^{4\ell+2}}
\\
&
+2\,\rho^{\ell+1}
\ \ \ \ \ \ \ \ \ \ \ \ \ \ \ \ \ \ \ 
-4\,\rho^{2\ell+1}\,\cos\,\ell\theta
\ \ \ \ \ \ \ \ \ \ \ \ \ \ \ \ \ \ 
+2\,\rho^{3\ell+1}
\\
&
\ \ \ \ \ \ \ \ \ \ \ \ \ \ \ \ \ \ \ \ \ \ \ \ \ \ \ \ \ \ \ \ \ \ \
-4\,\rho^{2\ell+1}
\ \ \ \ \ \ \ \ \ \ \ \ \ \ \ \ \ \ \ \ \ \ \ \ \ \ \ \ \,
+8\,\rho^{3\ell+1}\,\cos\,\ell\theta 
\ \ \ \ \ \ \ \ \ \ \
-4\,\rho^{4\ell+1}
\\
&
\ \ \ \ \ \ \ \ \ \ \ \ \ \ \ \ \ \ \ \ \ \ \ \ \ \ \ \ \ \ \ \ \ \ \
\ \ \ \ \ \ \ \ \ \ \ \ \ \ \ \ \ \ \ \ \ \ \ \ \ \ \ \ \ \ \ \ \ \ \
\ \ 
-\,\zero{\rho^{2\ell+2}}
\ \ \ \ \ \ \ \ \ \ \ \ \ \ \ \ \ \
+2\,\rho^{3\ell+2}\,\cos\,\ell\theta
\ \ \ \ \ \ \ \ \ \ 
-\zero{\rho^{4\ell+2}}.
\endaligned
\]
Only $4 = 2 \cdot 2$ terms annihilate by pairs, and the question 
is reduced to determine whether there is nonnegativity:
\[
\aligned
0
&
\overset{\text{\bf ?}}{\,\,\leqslant\,\,}
2\,\rho^\ell\,
\big(
1-2\,\rho^{2\ell}
\big)\,
\big(
1
-
\cos\,\ell\theta
\big)
+
2\,
\rho^{\ell+1}\,
\big(
1
-
\cos\,(\ell+1)\theta
\big)
\\
&
\ \ \ \ \
+
\rho^{2\ell+1}\,
\Big[
4\,\cos\,\theta
+
4\,\cos\,(\ell+1)\theta
-
4\,\cos\,\ell\theta
-
4
\Big]
\\
&
\ \ \ \ \
+
\rho^{3\ell+1}\,
\Big[
-\,8\,\cos\,\theta
-
2\,\cos\,(\ell+1)\theta
+
8\,\cos\,\ell\theta
+
2
\Big]
\\
&
\ \ \ \ \
+
\rho^{3\ell+2}\,
\Big[
2\,\cos\,\ell\theta
-
2
\Big]
+
\rho^{4\ell+1}\,
\Big[
4\,\cos\,\theta
-
4
\Big]
\\
&
\,\,=:\,\,
f_{\ell,\rho}(\theta),
\endaligned
\] 
for a certain family $\big( f_{\ell,\rho} 
\big)_{1\leqslant \ell}^{0 \leqslant \rho \leqslant 1/4}$ of
$2\pi$-periodic functions. Since this is trivially satisfied
when $\rho = 0$, we shall from now on assume that:
\[
0
\,<\,
\rho
\,\leqslant\,
{\textstyle{\frac{1}{4}}}
\,=\,
0.25
\eqno
{\scriptstyle{(\text{\rm assumption})}}.
\]

For the first term of $f_{\ell,\rho}$, since we have:
\[
1
\,\,<\,\,
2\,\big(1-2\cdot0.25^2\big)
\,\,\leqslant\,\,
2\,\big(1-2\,\rho^{2\ell}\big)
\eqno
{\scriptstyle{(\forall\,\ell\,\geqslant\,1)}},
\]
after division by $\rho^\ell$, it would suffice to have,
with certain new minorinzing functions:
\[
g_{\ell, \rho}
\,\leqslant\,
{\textstyle{\frac{1}{\rho^\ell}}}\,
f_{\ell,\rho}, 
\]
the nonnegativity:
\[
\aligned
0
&
\overset{\text{\bf ?}}{\,\,\leqslant\,\,}
g_{\ell,\rho}(\theta)
\\
&
\,\,:=\,\,
1
-
\cos\,\ell\theta
+
2\,\rho\,
\big(
1
-
\cos\,(\ell+1)\theta
\big)
\\
&
\ \ \ \ \
+
\rho^{\ell+1}\,
\Big[
4\,\cos\,\theta
+
4\,\cos\,(\ell+1)\theta
-
4\,\cos\,\ell\theta
-
4
\Big]
\\
&
\ \ \ \ \
+
\rho^{2\ell+1}\,
\Big[
-\,8\,\cos\,\theta
-
2\,\cos\,(\ell+1)\theta
+
8\,\cos\,\ell\theta
+
2
\Big]
\\
&
\ \ \ \ \
+
\rho^{2\ell+2}\,
\Big[
2\,\cos\,\ell\theta
-
2
\Big]
+
\rho^{3\ell+1}\,
\Big[
4\,\cos\,\theta
-
4
\Big].
\endaligned
\]

For instance, again with the choice $\rho := 0.25$,
the graphs on the unit circle of the functions
\[
\theta
\,\longmapsto\,
g_{\ell,\rho}(\theta)
\eqno
{\scriptstyle{(-\pi\,\leqslant\,\theta\,\leqslant\,\pi)}}
\]
show up, respectively, for the three choices $\ell = 2, 5, 10$, as:

\medskip

\includegraphics[scale=0.30]{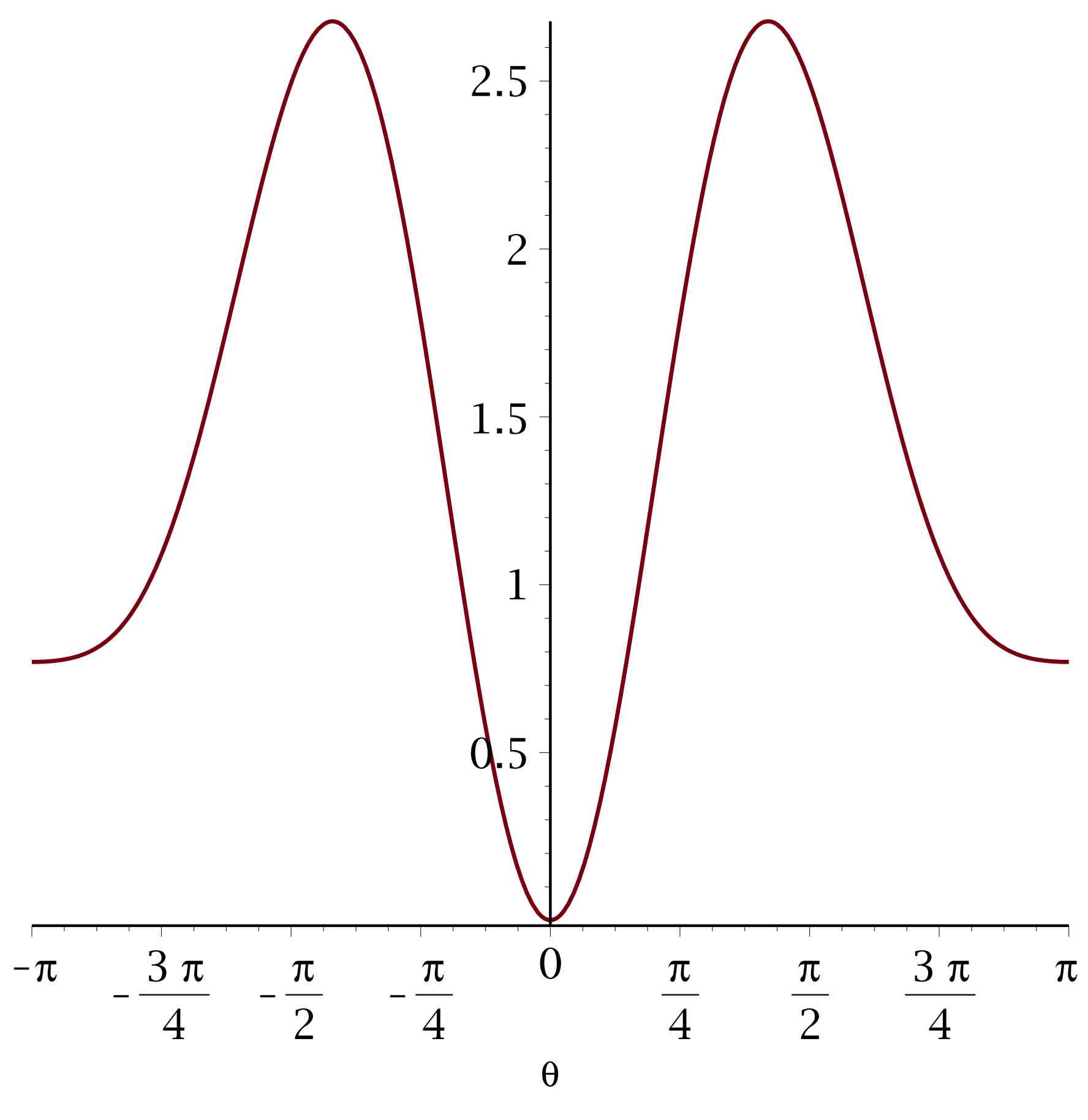}
\ \ \ \ \ \ \ 
\includegraphics[scale=0.30]{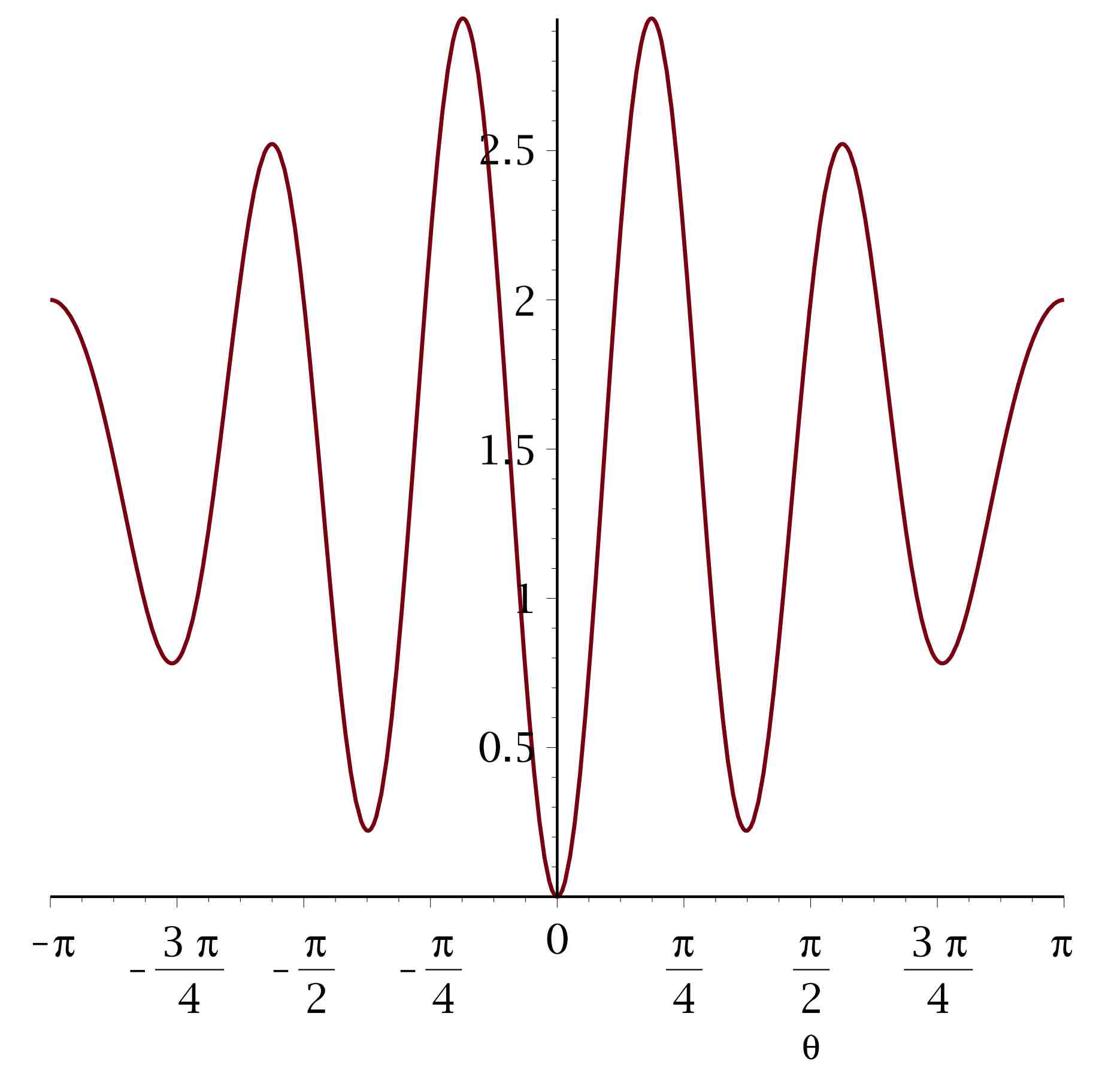}
\ \ \ \ \ \ \ 
\includegraphics[scale=0.30]{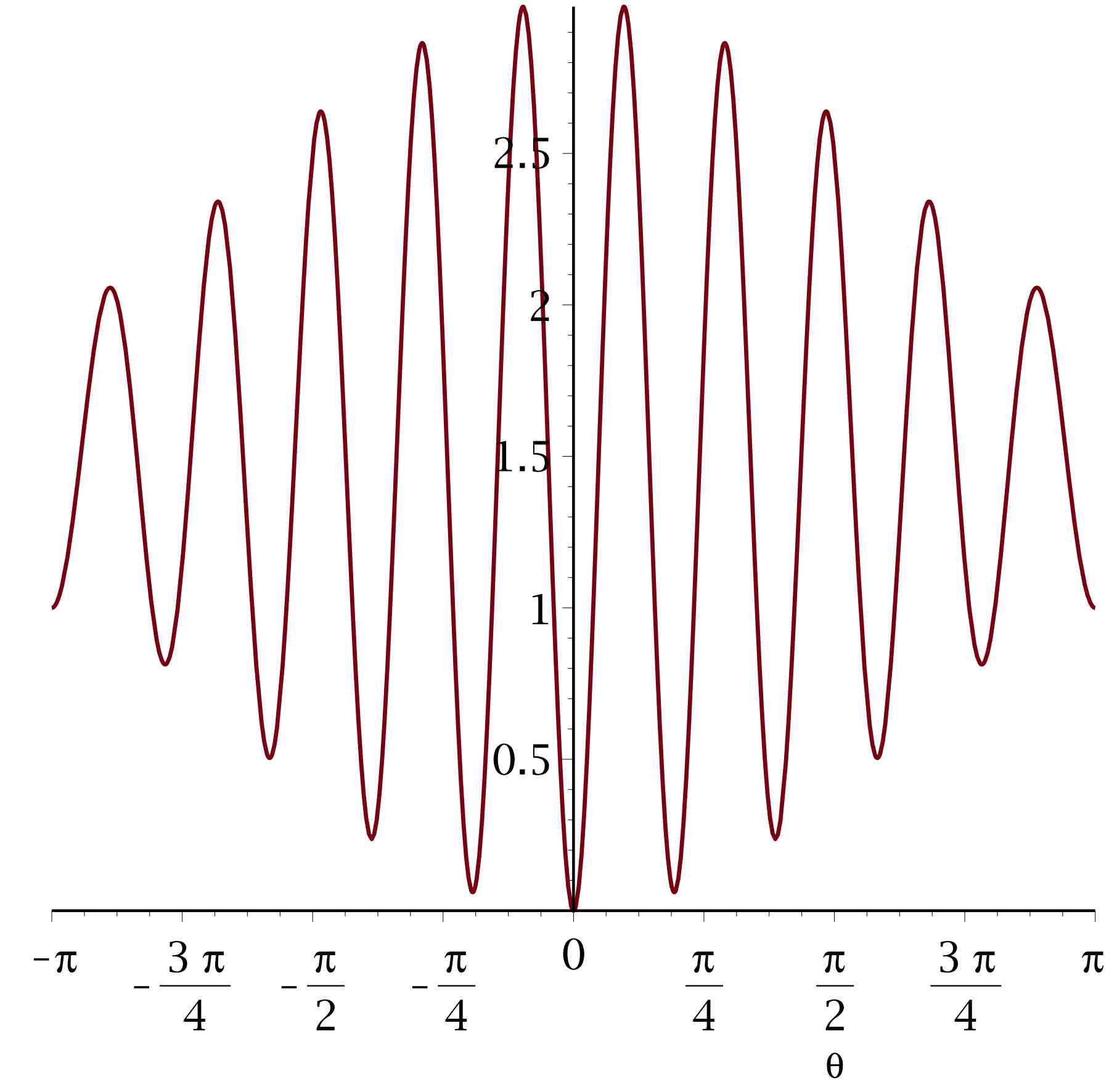}

Since these functions $g_{\ell,\rho}$ 
are even, it suffices to establish
their nonnegativity on $[0, \pi]$. Let us begin with 
examining their behavior in a right
half-neighborhood of $0$.
Starting from:
\[
0
\,=\,
g_{\ell,\rho}(0),
\]
a positivity of the first derivatives
of the $g_{\ell, \rho}$ would be welcome, 
at least on a small interval
like $]0,\, \frac{\pi}{4\ell} \big]$. 

\begin{Lemma}
For all real $0 < \rho \leqslant 0.25$
and for every integer $\ell \geqslant 1$, one has:
\[
g_{\ell,\rho}'(\theta)
\,\,>\,\,
0
\eqno
{\scriptstyle{(\forall\,\,
0\,<\,\theta\,\leqslant\,\frac{\pi}{4\ell})}}.
\]
\end{Lemma}

\proof
Observe that this is true even when $\rho = 0$, since
the function $g_{\ell, 0}(\theta) = 1 - \cos\,\ell
\theta$ has derivarive $\ell\, \sin\, \ell\theta > 0$ on
$]0, \frac{\pi}{4\ell} \big]$. Anyway, we assume
$0 < \rho \leqslant 0.25$.

Our aim is to minorize this derivative:
\[
\aligned
g_{\ell,\rho}'(\theta)
&
\,\,=\,\,
\ell\,\sin\,\ell\theta
+
2\,\rho\,(\ell+1)\,\sin\,(\ell+1)\theta
\\
&
\ \ \ \ \
+
\rho^{\ell+1}\,
\Big[
-4\,\sin\,\theta
-
4(\ell+1)\,\sin\,(\ell+1)\theta
+
4\ell\,\sin\,\ell\theta
\Big]
\\
&
\ \ \ \ \
+
\rho^{2\ell+1}\,
\Big[
8\,\sin\,\theta
+
2(\ell+1)\,\sin\,(\ell+1)\theta
-
8\ell\,\sin\,\ell\theta
\Big]
\\
&
\ \ \ \ \
+
\rho^{2\ell+2}\,
\Big[
-2\ell\,\sin\,\ell\theta
\Big]
+
\rho^{3\ell+1}\,
\Big[
-4\,\sin\,\theta
\Big],
\endaligned
\]
by a quantity which can be seen to be positive.
However, we have to treat the special
case $\ell = 1$ separately, namely for all
$0 < \rho \leqslant \frac{1}{4}$ and for all $0 < \theta 
\leqslant \frac{\pi}{4}$, we first check that:
\[
\aligned
g_{1,\rho}'(\theta)
&
\,\,=\,\,
\sin\,\theta
+
4\rho\,\sin\,2\theta
+
\rho^2\,
\big[
-\,8\,\sin\,2\theta
\big]
+
\rho^3\,
\big[
4\,\sin\,2\theta
\big]
+
\rho^4\,
\big[
-6\,\sin\,\theta
\big]
\\
&
\,\,=\,\,
\sin\,\theta\,
\Big\{
1
+
8\rho\,\cos\,\theta
-
\rho^2\,16\,\cos\,\theta
+
\rho^3\,8\,\cos\,\theta
-
\rho^4\,6
\Big\}
\\
&
\,\,\geqslant\,\,
\sin\,\theta\,
\Big\{
\zero{1}
+
8\,\rho\,
{\textstyle{\frac{1}{\sqrt{2}}}}
-
\zero{{\textstyle{\frac{1}{4^2}}}\,16}
+
\rho^3\,
{\textstyle{\frac{8}{\sqrt{2}}}}
-
\rho^4\,6
\Big\}
\\
&
\,\,\geqslant\,\,
\rho\,
\sin\,\theta
\cdot
\Big\{
{\textstyle{\frac{8}{\sqrt{2}}}}
+
0^2\,
{\textstyle{\frac{8}{\sqrt{2}}}}
-
{\textstyle{\frac{1}{4^3}}}\,6
\Big\}
\\
&
\,\,=\,\,
\rho\,\sin\,\theta
\cdot
5,563\cdots
\\
&
\,\,>\,\,
0.
\endaligned
\]

So we may assume $\ell \geqslant 2$.
If we use the classical inequalities valid for $\varphi \in 
[0,\pi]$:
\[
\sin\,\varphi
\,\geqslant\,
\varphi
-
{\textstyle{\frac{1}{6}}}\,
\varphi^3
\ \ \ \ \ \ \ \ \ \ \ \ \
\text{and}
\ \ \ \ \ \ \ \ \ \ \ \ \
-\,\sin\,\varphi
\,\geqslant\,
-\varphi,
\]
we are conducted to ask ourselves whether:
\reqnomode\usetagform{EngelLie}
\begin{align}
g_{\ell,\rho}'(\theta)
&
\,\,\geqslant\,\,
\ell\,
\big(
\ell\theta
-
{\textstyle{\frac{1}{6}}}
(\ell\theta)^3
\big)
+
2\,\rho\,(\ell+1)\,
\big(
(\ell+1)\theta
-
{\textstyle{\frac{1}{6}}}
((\ell+1)\theta)^3
\big)
\notag
\\
&
\ \ \ \ \
+
\rho^{\ell+1}\,
\Big[
-4\theta
-
4(\ell+1)\,(\ell+1)\theta
+
4\ell\,\big(
\ell\theta
-
{\textstyle{\frac{1}{6}}}
(\ell\theta)^3
\big)
\Big]
\notag
\\
&
\ \ \ \ \
+
\rho^{2\ell+1}\,
\Big[
8\,\big(
\theta
-
{\textstyle{\frac{1}{6}}}
\theta^3
\big)
+
2\,(\ell+1)\,
\big(
(\ell+1)\theta
-
{\textstyle{\frac{1}{6}}}
((\ell+1)\theta)^3
\big)
-
8\ell\,\ell\theta
\Big]
\notag
\\
&
\ \ \ \ \
+
\rho^{2\ell+2}\,
\big[
-2\ell\,\ell\theta
\big]
+
\rho^{3\ell+1}\,
\big[
-4\,\theta
\big]
\notag
\\
&
\overset{\text{\bf ?}}{\,\,>\,\,}
0
\tag{(\forall\,\,0\,<\,\theta\,\leqslant\,\frac{\pi}{4\ell}).}
\end{align}
To have a better view, let us set:
\[
t
\,:=\,
\ell\,\theta,
\ \ \ \ \ \ \ \ \ \ \ \ \
\text{whence}
\ \ \ \ \ \ \ \ \ \ \ \ \
0
\,<\,
t
\,\leqslant\,
{\textstyle{\frac{\pi}{4}}}
\,<\,
1,
\]
and let us simplify this minorant by writing $(\ell+1)\theta
= \frac{\ell+1}{\ell}\, \ell\theta = \frac{\ell+1}{\ell}\, t$:
\reqnomode\usetagform{EngelLie}
\begin{align}
g_{\ell,\rho}'(\theta)
&
\,\,\geqslant\,\,
\ell\,t\,
\big(
1
-
{\textstyle{\frac{1}{6}}}\,
t^2
\big)
+
2\,\rho\,(\ell+1)\,
{\textstyle{\frac{\ell+1}{\ell}}}\,
t
\big(
1
-
{\textstyle{\frac{1}{6}}}
\big({\textstyle{\frac{\ell+1}{\ell}}}\big)^2\,
t^2
\big)
\notag
\\
&
\ \ \ \ \
+
\rho^{\ell+1}\,
\Big[
-{\textstyle{\frac{4}{\ell}}}\,t
-
4(\ell+1)\,
{\textstyle{\frac{\ell+1}{\ell}}}\,
t
+
4\ell\,t\,
\big(
1
-
{\textstyle{\frac{1}{6}}}\,
t^2
\big)
\Big]
\notag
\\
&
\ \ \ \ \
+
\rho^{2\ell+1}\,
\Big[
8\,
{\textstyle{\frac{t}{\ell}}}\,
\big(
1
-
{\textstyle{\frac{1}{6}}}\,
\big(
{\textstyle{\frac{t}{\ell}}}
\big)^2
\big)
+
2\,(\ell+1)\,
{\textstyle{\frac{\ell+1}{\ell}}}\,
t\,
\big(
1
-
{\textstyle{\frac{1}{6}}}
\big({\textstyle{\frac{\ell+1}{\ell}}}\big)^2\,
t^2
\big)
-
8\ell\,t
\Big]
\notag
\\
&
\ \ \ \ \
+
\rho^{2\ell+2}\,
\big[
-2\ell\,t
\big]
+
\rho^{3\ell+1}\,
\big[
-\,{\textstyle{\frac{4}{\ell}}}\,t
\big]
\notag
\\
&
\overset{\text{\bf ?}}{\,\,>\,\,}
0
\tag{(\forall\,\,0\,<\,\theta\,\leqslant\,\frac{\pi}{4\ell}).}
\end{align}

In order to minorize this by an even simpler quantity, we can
use, since $\ell \geqslant 2$:
\[
-\,\big(
{\textstyle{\frac{\ell+1}{\ell}}}
\big)^2
\,\geqslant\,
-\,\big({\textstyle{\frac{3}{2}}}\big)^2
\ \ \ \ \ \ \ \ \ \ \ \ \
\text{and also}
\ \ \ \ \ \ \ \ \ \ \ \ \
-\,t^2
\,\geqslant\,
-1,
\]
so that:
\[
\aligned
{\sf preceding}\,\,{\sf minorant}
&
\,\,\geqslant\,\,
\ell\,t\,
{\textstyle{\frac{5}{6}}}
+
\underline{
2\,\rho\,(\ell+1)\,
{\textstyle{\frac{\ell+1}{\ell}}}\,
t\,
{\textstyle{\frac{5}{8}}}}
\\
&
\ \ \ \ \
+
\rho^{\ell+1}\,
\Big[
-2\,t
-
4\,(\ell+1)\,
{\textstyle{\frac{3}{2}}}\,
t
+
4\,\ell\,t\,
{\textstyle{\frac{5}{6}}}
\Big]
\\
&
\ \ \ \ \
+
\rho^{2\ell+1}\,
\Big[
8\cdot0
+
2\,(\ell+1)\,1\,t\,
{\textstyle{\frac{5}{8}}}
-
8\,\ell\,t
\Big]
\\
&
\ \ \ \ \
+
\rho^{2\ell+2}\,
\big[
-2\,\ell\,t
\big]
+
\rho^{3\ell+1}\,
\big[
-2\,t
\big]
\ \ \ \ \ \ \ \ \ \ \ \ \ \ 
\overset{\text{\bf ?}}{\,\,>\,\,}
0,
\endaligned
\]
and we even once more minorize this intermediate 
minorant by neglecting the
term underlined and summing the expressions in brackets:
\[
\aligned
g_{\ell,\rho}'(\theta)
&
\,\,\geqslant\,\,
\ell\,t\,
{\textstyle{\frac{5}{6}}}
+
0
\\
&
\ \ \ \ \
+
\rho^{\ell+1}\,
\Big[
-8\,t
+
{\textstyle{\frac{16}{6}}}\,
\ell\,t
\Big]
\\
&
\ \ \ \ \
+
\rho^{2\ell+1}\,
\Big[
{\textstyle{\frac{5}{4}}}\,
t
-
{\textstyle{\frac{27}{4}}}\,
\ell\,t
\Big]
\\
&
\ \ \ \ \
+
\rho^{2\ell+2}\,
\big[
-2\,\ell\,t
\big]
+
\rho^{3\ell+1}\,
\big[
-2\,t
\big]
\ \ \ \ \ \ \ \ \ \ \ \ \ \ 
\overset{\text{\bf ?}}{\,\,>\,\,}
0.
\endaligned
\]
We conclude by a factorization and by a final computer 
check, still for all $\ell \geqslant 2$: 
\begin{align}
g_{\ell,\rho}'(\theta)
&
\,\,\geqslant\,\,
t\,
\Big\{
{\textstyle{\frac{5}{6}}}\,\ell
-
\rho^{\ell+1}\,
\big[
8
+
{\textstyle{\frac{16}{6}}}\,\ell
\big]
-
\rho^{2\ell+1}\,
\big[
{\textstyle{\frac{5}{4}}}
+
{\textstyle{\frac{27}{4}}}\,
\ell
\big]
-
\rho^{2\ell+2}\,
\big[2\,\ell\big]
-
\rho^{3\ell+1}\,
\big[2\big]
\Big\}
\notag
\\
&
\,\,\geqslant\,\,
t\,
\Big\{
{\textstyle{\frac{5}{6}}}\,\ell
-
0.25^{\ell+1}\,
\big[
8
+
{\textstyle{\frac{16}{6}}}\,\ell
\big]
-
0.25^{2\ell+1}\,
\big[
{\textstyle{\frac{5}{4}}}
+
{\textstyle{\frac{27}{4}}}\,
\ell
\big]
-
0.25^{2\ell+2}\,
\big[2\,\ell\big]
-
0.25^{3\ell+1}\,
\big[2\big]
\Big\}
\notag
\\
&
\,\,\geqslant\,\,
t
\cdot
1,442\cdots.
\qedhere
\end{align}
\endproof

In summary, we have established for all $\ell \geqslant 1$
the positivity on a starting interval:
\[
0
\,<\,
g_{\ell,\rho}(\theta)
\eqno
{\scriptstyle{(\forall\,\theta\,\in\,
]0,\,\frac{\pi}{4\ell}])}},
\]
and our next goal is to establish the positivity of this minoring
function $g_{\ell,\rho}$ 
on the remaining (large) subinterval of $[0, \pi]$:
\[
0
\overset{\text{\bf ?}}{\,\,<\,\,}
g_{\ell,\rho}(\theta)
\eqno
{\scriptstyle{(\forall\,\,\theta\,\in\,
[\frac{\pi}{4\ell},\pi])}}.
\]
We first finish the case $\ell = 1$.

\begin{Lemma}
For all $0 < \rho \leqslant 0.25$, the function $g_{1, \rho}(\theta)$
is positive on $\big[ \frac{\pi}{4\cdot 1}, \pi]$.
\end{Lemma}

\proof
Indeed:
\begin{align}
g_{1,\rho}(\theta)
&
\,=\,
1
-
\cos\,\theta
+
2\rho\,\big(1-\cos\,2\theta\big)
\notag
\\
&
\ \ \ \ \
+
\rho^2\,
\big[
4\,\cos\,2\theta
-
4
\big]
\notag
\\
&
\ \ \ \ \
+
\rho^3\,
\big[
-2\,\cos\,2\theta+2
\big]
\notag
\\
&
\ \ \ \ \
+
\rho^4\,
\big[6\,\cos\,\theta-6\big]
\notag
\\
&
\,=\,
\big(1-\cos\,\theta\big)\,
\big[
1
-
6\,\rho^4
\big]
+
\big(1-\cos\,2\theta\big)\,
\big[
2\rho\,
\big(1-2\,\rho+\rho^2\big)
\big]
\notag
\\
&
\,\geqslant\,
\big(
1
-
{\textstyle{\frac{1}{\sqrt{2}}}}
\big)\,
\big[
1
-
6\cdot 0.25^4
\big]
+
{\sf nonnegative}
\notag
\\
&
\,=\,
0.286\cdots
\notag
\\
&
\,>\,
0.
\qedhere
\end{align}
\endproof

From now one, when we work on $\big[ \frac{\pi}{4\ell}, \pi]$,
we can therefore assume that:
\[
\ell
\,\geqslant\,
2.
\]

\begin{Lemma}
For all $0 < \rho < 0.25$ and every integer $\ell \geqslant 2$,
there is on $\big[ \frac{\pi}{4\ell}, \pi]$ a minoration:
\[
0
\overset{\text{\bf ?}}{\,\,<\,\,}
h_{\ell,\rho}(\theta)
\,\,\leqslant\,\,
g_{\ell,\rho}(\theta)
\,\,\leqslant\,\,
{\textstyle{\frac{1}{\rho^\ell}}}\,
f_{\ell,\rho}(\theta),
\]
with the new:
\[
h_{\ell,\rho}(\theta)
\,:=\,
1
-
\cos\,\ell\theta
+
2\rho\,
\big(
1
-
\cos\,(\ell+1)\theta
\big)
-
18\,\rho^{\ell+1}.
\]
\end{Lemma}

Again with the choice $\rho := 0.25$,
the graphs on the unit circle of the functions
\[
\theta
\,\longmapsto\,
h_{\ell,\rho}(\theta)
\eqno
{\scriptstyle{(-\pi\,\leqslant\,\theta\,\leqslant\,\pi)}}
\]
show up, respectively, for the three choices $\ell = 2, 5, 10$, as:

\medskip

\includegraphics[scale=0.30]{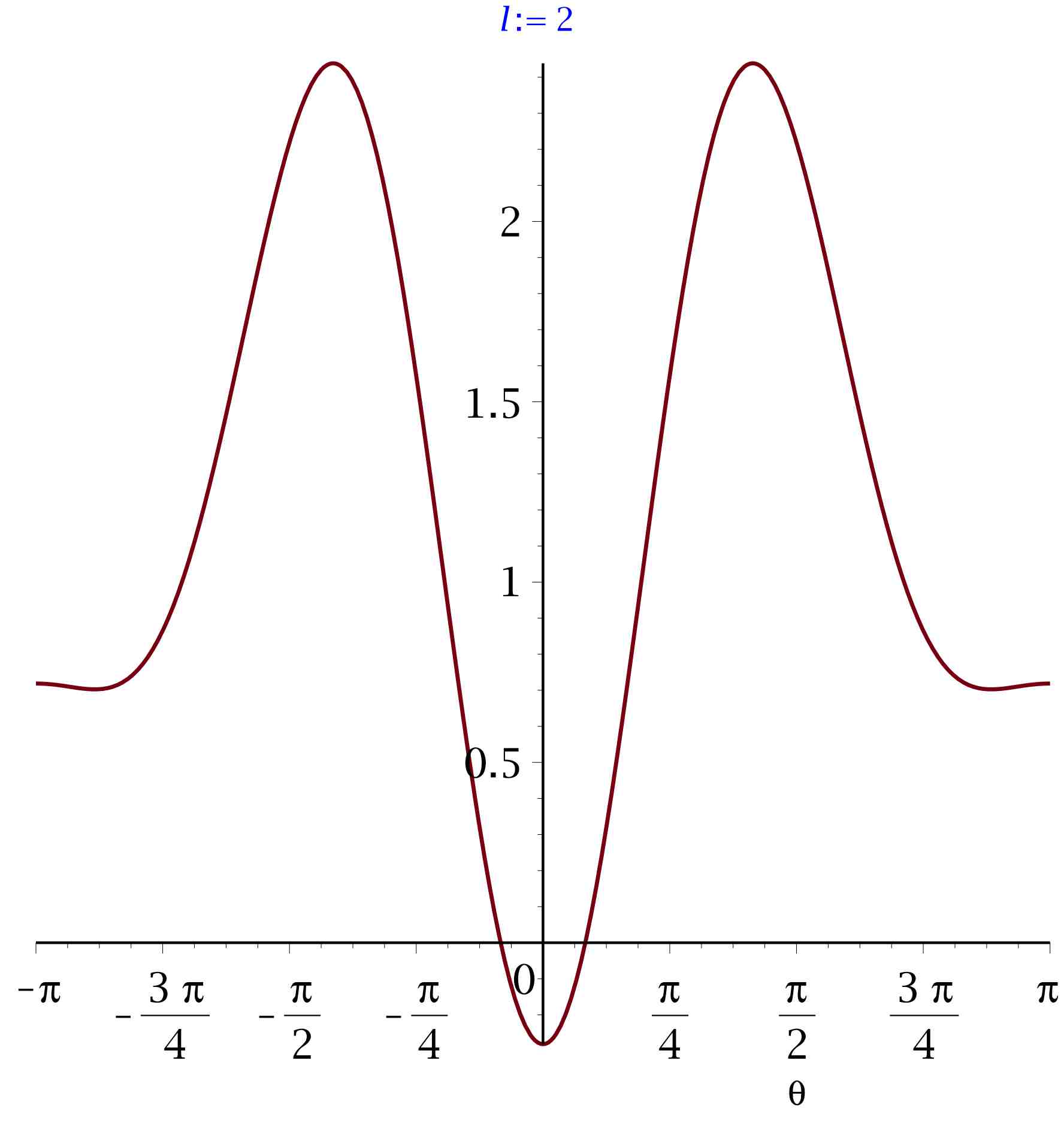}
\ \ \ \ \ \ \ 
\includegraphics[scale=0.30]{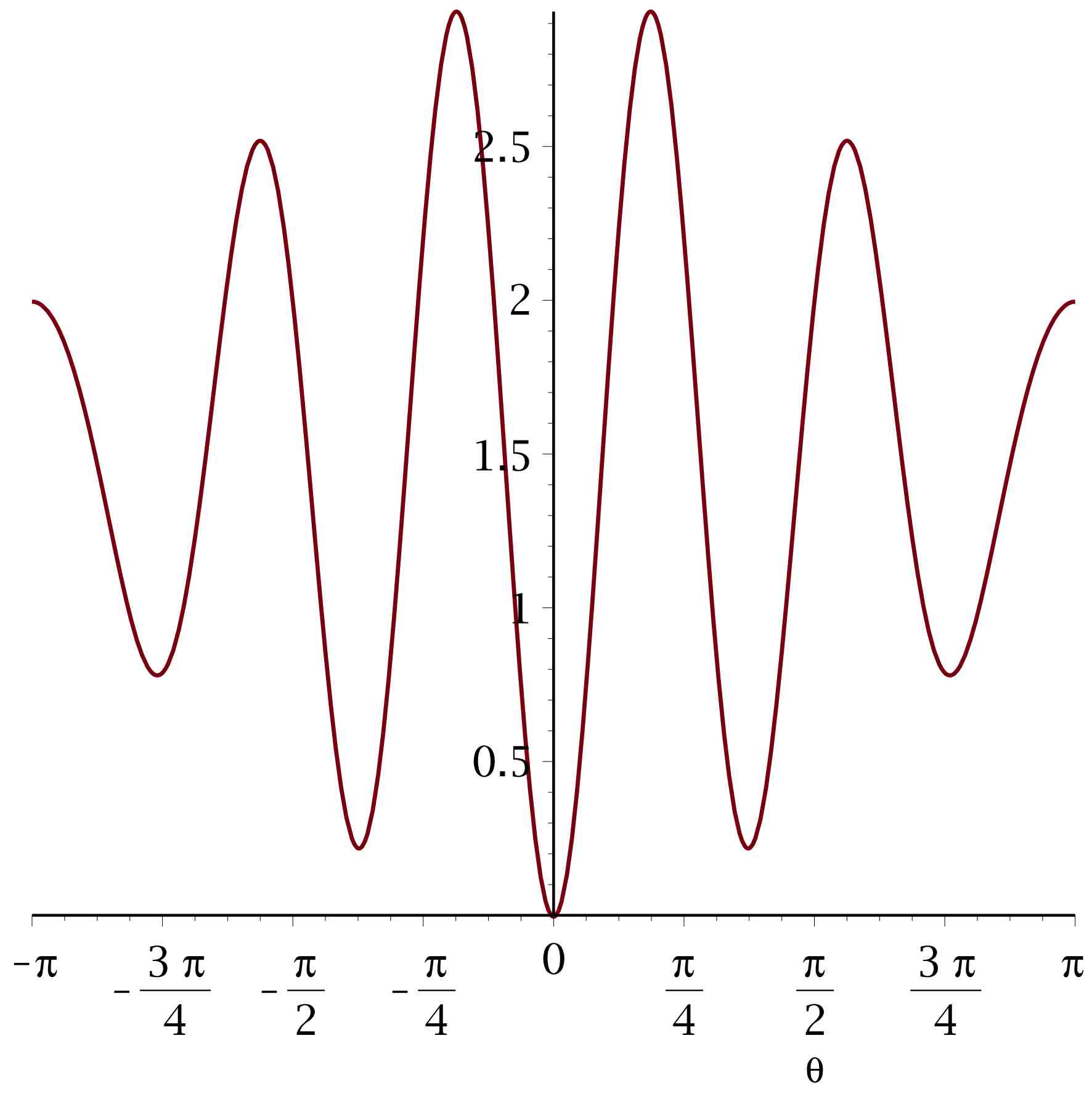}
\ \ \ \ \ \ \ 
\includegraphics[scale=0.30]{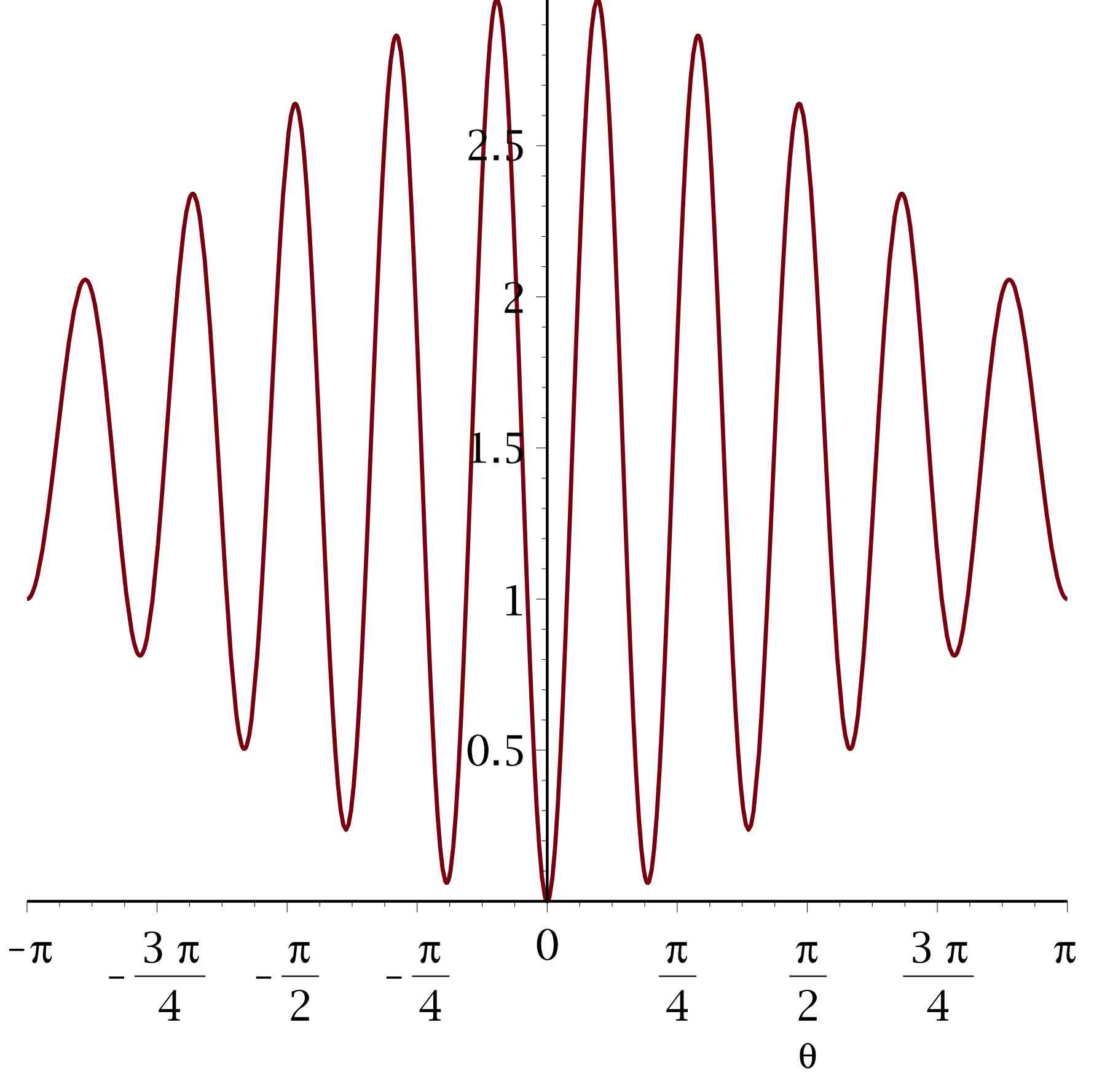}

\proof
Indeed, we minorize simply the reminders:
\[
\aligned
g_{\ell,\rho}(\theta)
&
\,\,\geqslant\,\,
1
-
\cos\,\ell\theta
+
2\rho\,
\big(
1
-
\cos\,(\ell+1)\theta
\big)
\\
&
\ \ \ \ \
-\,
\rho^{\ell+1}\,
\big[
4+4+4+4
\big]
-
\rho^{2\ell+1}\,
\big[8+2+8+2\big]
-
\rho^{2\ell+2}\,
\big[2+2\big]
-
\rho^{3\ell+1}\,
\big[4+4\big],
\endaligned
\]
and to even simplify the second line by replacing it by
$-18\, \rho^{\ell+1}$ as announced, we assert that:\[
-16\,\rho^{\ell+1}
-
20\,\rho^{2\ell+1}
-
4\,\rho^{2\ell+2}
-
8\,\rho^{3\ell+1}
\,\,\geqslant\,\,
-\,18\,\rho^{\ell+1},
\]
simply since:
\begin{align}
\rho^{\ell+1}\,
\big(
1
-
10\,\rho^\ell
-
2\,\rho^{\ell+1}
-
4\,\rho^{2\ell}
\big)
&
\,\,\geqslant\,\,
\rho^{\ell+1}\,
\big(
1
-
10\cdot0.25^2
-
2\cdot0.25^3
-
4\cdot 0.25^4
\big)
\notag
\\
&
\,\,>\,\,
\rho^{\ell+1}\,
\big(
0.328\cdots
\big)
\notag
\\
&
\,>\,
0.
\qedhere
\end{align}
\endproof

It therefore remains to treat all the cases $\ell \geqslant 2$.
We start by looking at the subintervals
$\big[\frac{\pi}{4\ell},\, \frac{7\pi}{4\ell} \big]
\subset \big[ \frac{\pi}{4\ell}, \pi]$.  

\begin{Lemma}
For every real $0 < \rho \leqslant 0.25$ and
every integer $\ell \geqslant 2$:
\[
0
\,<\,
h_{\ell,\rho}(\theta)
\eqno
{\scriptstyle{(\forall\,\,\theta\,\in\,
[\frac{\pi}{4\ell},\frac{7\pi}{4\ell}])}}.
\]
\end{Lemma}

\proof
Since:
\[
{\textstyle{\frac{\pi}{4}}}
\,\leqslant\,
\ell\,\theta
\,\leqslant\,
{\textstyle{\frac{7\pi}{4}}}
\]
it comes:
\[
1
-
\cos\,\ell\theta
\,\geqslant\,
1
-
\frac{1}{\sqrt{2}},
\]
and since $1 - \cos\,(\ell+1)\theta \geqslant 0$ anyway, 
we can minorize:
\begin{align}
h_{\ell,\rho}(\theta)
&
\,\geqslant\,
1
-
\cos\,\ell\theta
+
0
-
18\,\rho^{\ell+1}
\notag
\\
&
\,\geqslant\,
1
-
{\textstyle{\frac{1}{\sqrt{2}}}}
-
18\cdot
0.25^{2+1}
\notag
\\
&
\,=\,
0.01164\cdots.
\qedhere
\end{align}
\endproof

We can now finish the case $\ell = 2$. It remains to show
positivity of $h_{2,\rho}(\theta)$ on $\big[ \frac{7\pi}{8}, \pi]$.
Since $\frac{21\pi}{8} \leqslant 3\,\theta \leqslant 3\,\pi$,
or equivalently $\frac{5\pi}{8} \leqslant 3\,\theta - \pi \leqslant 
\pi$, we can minorize:
\[
\aligned
h_{2,\rho}(\theta)
&
\,=\,
1
-
\cos\,2\theta
+
2\rho\,
\big(
1-\cos\,3\theta
\big)
-
18\,\rho^3
\\
&
\,\geqslant\,
0
+
2\rho\,
\big[
\big(
1
-
\cos\,
{\textstyle{\frac{5\pi}{8}}}
\big)
-
9\cdot 0.25^2
\big]
\\
&
\,=\,
2\rho
\cdot
0.820\cdots
\\
&
\,>\,
0.
\endaligned
\] 

It still remains to treat all the cases $\ell \geqslant 3$.

\begin{Lemma}
\label{Lemma-ell-3-h-ell-rho}
For every real $0 \leqslant \rho \leqslant 0.25$ and
every integer $\ell \geqslant 3$, the function:
\[
h_{\ell,\rho}(\theta)
\,:=\,
1
-
\cos\,\ell\theta
+
2\,\rho\,
\big(
1
-
\cos\,(\ell+1)\theta
\big)
-
18\,\rho^{\ell+1}
\]
takes only positive values in the interval 
$\big[\frac{7\pi}{4\ell}, \pi]$:
\[
h_{\ell,\rho}(\theta)
\,>\,
0
\eqno
{\scriptstyle{(\forall\,\,
\frac{7\pi}{4\ell}\,\leqslant\,\theta\,\leqslant\,\pi)}}.
\]
\end{Lemma}

\proof
Using $1 - \cos\, \varphi = 2\, \sin^2 \frac{\varphi}{2}$,
let us rewrite:
\[
h_{\ell,\rho}(\theta)
\,=\,
2\,\sin^2\,
{\textstyle{\frac{\ell\,\theta}{2}}}
+
4\,\rho\,\sin^2\,
{\textstyle{\frac{(\ell+1)\,\theta}{2}}}
-
18\,\rho^{\ell+1}.
\]
At a point $\theta \in \big[\frac{7\pi}{4\ell}, \pi \big]$,
if we have either:
\[
2\,\sin^2\,
{\textstyle{\frac{\ell\,\theta}{2}}}
-
18\,\rho^{\ell+1}
\,\,>\,\,
0
\ \ \ \ \ \ \ \ \ \ \ \ \
\text{or}
\ \ \ \ \ \ \ \ \ \ \ \ \
4\,\rho\,\sin^2\,
{\textstyle{\frac{(\ell+1)\,\theta}{2}}}
-
18\,\rho^{\ell+1}
\,\,>\,\,
0,
\]
then there is nothing to prove. We claim that the opposite 
inequalities cannot hold. 

\begin{Assertion}
For every $\ell \geqslant 3$, there is no 
$\theta \in \big[ \frac{7\pi}{4\ell}, \pi]$ at which:
\[
\sin^2\,
{\textstyle{\frac{\ell\,\theta}{2}}}
\,\,\leqslant\,\,
9\,\rho^{\ell+1}
\ \ \ \ \ \ \ \ \ \ \ \ \
\text{and}
\ \ \ \ \ \ \ \ \ \ \ \ \
\sin^2\,
{\textstyle{\frac{(\ell+1)\,\theta}{2}}}
\,\,\leqslant\,\,
{\textstyle{\frac{9}{2}}}\,
\rho^\ell.
\]
\end{Assertion}

\proof
Suppose nevertheless that such a $\theta \in \big[ 
\frac{7\pi}{4\ell}, \pi]$ exists. Modulo $\pi$, there exist
two unique representatives $-\frac{\pi}{2} < \alpha, \beta 
\leqslant \frac{\pi}{2}$ of:
\[
{\textstyle{\frac{\ell\,\theta}{2}}}
-
p\,\pi
\,=\,
\alpha
\ \ \ \ \ \ \ \ \ \ \ \ \
\text{and}
\ \ \ \ \ \ \ \ \ \ \ \ \
{\textstyle{\frac{(\ell+1)\,\theta}{2}}}
-
q\,\pi
\,=\,
\beta,
\]
with certain unique integers $p, q \in \Z$, whence:
\[
\sin^2\alpha
\,=\,
\sin^2\,
{\textstyle{\frac{\ell\,\theta}{2}}}
\,\,\leqslant\,\,
9\,\rho^{\ell+1}
\ \ \ \ \ \ \ \ \ \ \ \ \
\text{and}
\ \ \ \ \ \ \ \ \ \ \ \ \
\sin^2\beta
\,=\,
\sin^2\,
{\textstyle{\frac{(\ell+1)\,\theta}{2}}}
\,\,\leqslant\,\,
{\textstyle{\frac{9}{2}}}\,
\rho^\ell.
\]

\begin{Lemma}
For all $0 \leqslant \vert \gamma \vert \leqslant \frac{\pi}{2}$,
one has the classical inequality $\frac{\vert \gamma\vert}{2}
\leqslant \vert \sin\,\gamma \vert \leqslant \vert \gamma \vert$.\qed
\end{Lemma}

Consequently, using $\rho^{1/2} \leqslant 0.25^{1/2} = \frac{1}{2}$,
it comes:
\[
{\textstyle{\frac{\vert\alpha\vert}{2}}}
\,\leqslant\,
\vert\sin\,\alpha\vert
\,\leqslant\,
3\,\rho^{\frac{\ell+1}{2}}
\,\leqslant\,
3\,
{\textstyle{\frac{1}{2^{\ell+1}}}}
\ \ \ \ \ \ \ \ \ \ \ \ \
\text{and}
\ \ \ \ \ \ \ \ \ \ \ \ \
{\textstyle{\frac{\vert\beta\vert}{2}}}
\,\leqslant\,
\vert\sin\,\beta\vert
\,\leqslant\,
{\textstyle{\frac{3}{\sqrt{2}}}}\,
\rho^{\frac{\ell}{2}}
\,\leqslant\,
{\textstyle{\frac{3}{\sqrt{2}}}}\,
{\textstyle{\frac{1}{2^{\ell}}}},
\]
that is to say:
\[
\vert\alpha\vert
\,\leqslant\,
3\,
{\textstyle{\frac{1}{2^{\ell}}}}
\ \ \ \ \ \ \ \ \ \ \ \ \
\text{and}
\ \ \ \ \ \ \ \ \ \ \ \ \
\vert\beta\vert
\,\leqslant\,
3\,\sqrt{2}\,
{\textstyle{\frac{1}{2^{\ell}}}}.
\]
From a chain of estimations:
\[
\aligned
3\,
{\textstyle{\frac{1+\sqrt{2}}{2^{\ell}}}}
\,\,\geqslant\,\,
\vert\beta-\alpha\vert
&
\,\,=\,\,
\big\vert
{\textstyle{\frac{\theta}{2}}}
+
(q-p)\,\pi
\big\vert
\\
&
\,\,\geqslant\,\,
\big\vert
{\textstyle{\frac{1}{2}}}\,
{\textstyle{\frac{7\pi}{4\ell}}}
+
(q-p)\,\pi
\big\vert
\\
&
\,\,\geqslant\,\,
\underset{r\in\Z}{\min}\,
\big\vert
{\textstyle{\frac{7\pi}{8\ell}}}
-
r\,\pi
\big\vert
\\
\explicationmath{{\textstyle{\frac{7\pi}{8\ell}}}
\leqslant{\textstyle{\frac{7\pi}{64}}}}
\ \ \ \ \ \ \ \ \ \ \ \ \ \ \ \ \ \ \ \ \ \ \ \ \ \
&
\,\,=\,\,
{\textstyle{\frac{7\pi}{8\ell}}}
\endaligned
\]
we obtain the inequality:
\[
3\,
{\textstyle{\frac{1+\sqrt{2}}{2^{\ell}}}}
\,\,\geqslant\,\,
{\textstyle{\frac{7\pi}{8\ell}}},
\]
which is visibly false for large $\ell$, and which begins
to be false when $\ell \geqslant 3$: 
\[
3\,
{\textstyle{\frac{1+\sqrt{2}}{2^3}}}
\,=\,
0.905\cdots
\overset{\text{\bf no!}}{\,\,\geqslant\,\,}
0.916\cdots
\,=\,
{\textstyle{\frac{7\pi}{8\cdot 3}}},
\]
This contradiction
proves the assertion.
\endproof

The proof of 
Lemma~{\ref{Lemma-ell-3-h-ell-rho}} 
is complete.
\endproof

The proof of 
Proposition~{\ref{Proposition-G-k-z-H-ell-z}} is complete.
\endproof


\linestop


\vfill\end{document}